\documentclass[pdflatex,sn-mathphys-num]{sn-jnl}

\usepackage{graphicx}%
\usepackage{multirow}%
\usepackage{amsmath,amssymb,amsfonts}%
\usepackage{amsthm}%
\usepackage{mathrsfs}%
\usepackage[title]{appendix}%
\usepackage{xcolor}%
\usepackage{textcomp}%
\usepackage{manyfoot}%
\usepackage{booktabs}%
\usepackage{algorithm}%
\usepackage{algorithmicx}%
\usepackage{algpseudocode}%
\usepackage{listings}%
\usepackage{bm}
\usepackage{enumitem}
\usepackage{arydshln}
\usepackage{multirow}
\usepackage{graphicx}

\newcommand{\abs}[1]{\left\vert#1\right\vert}

\def\RR{{\mathbb{R}}}

\theoremstyle{thmstyleone}%
\newtheorem{theorem}{Theorem}
\theoremstyle{thmstyletwo}%
\newtheorem{remark}{Remark}%

\theoremstyle{thmstylethree}%
\newtheorem{definition}{Definition}%

\begin{document}

\title[Error  and long-term analysis of two-step symmetric  methods for relativistic charged-particle dynamics]{Error  and long-term analysis of two-step symmetric  methods for relativistic charged-particle dynamics}


\author*[1,2]{\fnm{Ting} \sur{Li}}\email{1009587520@qq.com}

\author[1]{\fnm{Bin} \sur{Wang}}\email{wangbinmaths@xjtu.edu.cn}
\author[3]{\fnm{Ruili} \sur{Zhang}}\email{zhangrl@bjtu.edu.cn}
\equalcont{These authors contributed equally to this work.}

\affil[1]{School of Mathematics and Statistics, Xi'an Jiaotong University, 710049 Xi'an, China}
\affil[2]{Mathematisches Institut, University of T\"{u}bingen, 72076 T\"{u}bingen, Germany}
\affil[3]{School of Mathematics and Statistics, Beijing Jiaotong University, 100044 Beijing, China}


\abstract{In this work, we consider the error estimates and the long-time conservation or near-conservation of geometric structures, including energy, mass shell and phase-space volume, for four two-step symmetric methods applied to relativistic charged-particle dynamics.  We begin by introducing a two-step symmetric numerical method based on a splitting scheme that exactly preserves the mass shell and the phase-space volume of the relativistic system. Building on this formulation, we develop three additional two-step symmetric methods with further modifications, for which the long-time near-conservation of energy and mass shell can be rigorously established through the backward error analysis. All methods are shown to achieve second-order accuracy.  The theoretical results are illustrated and complemented by numerical experiments.
}




\maketitle

\section{Introduction}
Relativistic charged-particle dynamics in time-dependent electromagnetic fields is of fundamental theoretical interest and has significant applications in plasma physics \cite{Eriksson2003,Guan2010,RohrlichF2008}. The dynamics often involves multi-scale processes and long-term simulations. 
A comprehensive understanding of such dynamics is essential for both theoretical insights and practical applications.
The trajectory of a charged particle is obtained by solving the differential equation
\begin{equation}\label{rela-CPD}
	\begin{aligned}	\frac{dx}{d\bar{t}}=\frac{{\mu}}{\gamma},\quad \ 
		\frac{d{\mu}}{d\bar{t}}=\frac{{\mu}}{\gamma}\times{B(x)}+E(x),\ \ 
		\gamma=\sqrt{1+\mu^2},
	\end{aligned}
\end{equation}
where $x(\bar{t})\in \RR^3$ represents the position of a particle, ${\mu}(\bar{t})\in \RR^3$ is the momentum, and $\gamma$ is the relativistic factor. The vector potential $A{(x)} \in \RR^3 $
determines the magnetic field via $B(x) = \nabla \times A{( x)}$, while the electric field is given by $E(x)=-\nabla U(x)$ with a scalar potential $U(x)\in \RR$. Throughout this paper, we assume that $A(x)$ and $U(x)$ are smooth functions, and that all physical quantities, including the rest mass, speed of light, and charge of the particle, are normalized to 1 for convenience.

The relativistic charged-particle dynamics \eqref{rela-CPD} can equivalently be expressed as the Euler-Lagrange equations $\frac{d}{d\bar{t}}\nabla_{\nu}L
=\nabla_{x}L$
for the Lagrangian function
$L(x,\nu)=-\gamma^{-1}+A(x) \cdot \nu-U(x),$
where $\nu=\mu\gamma^{-1}$ and thus $\gamma^{-1}=\sqrt{1-\nu^2}$. Here, the dot $(\cdot)$ denotes the three-dimensional Euclidean inner product. By introducing the conjugate momentum $p=\nabla_{\nu}L=\gamma \nu+A{(x)}={\mu}+A{(x)}$, we apply the Legendre transformation $H=p \cdot \nu-L$ to derive the Hamiltonian (energy of the charged particle)
\begin{equation}\label{E}
	H(x,p)=\gamma+ U(x).
\end{equation}
To emphasize the dependence of the energy on $x$ and $\gamma$, we often write $H(x, \gamma)$ instead of using the canonical variables $(x, p)$.

The numerical methods investigated in this paper are based on a different formulation designed to work within four-dimensional (4D) spacetime. The equations \eqref{rela-CPD} can be written in terms of the proper time $\tau$
\begin{equation*}
	\frac{d x}{d \tau}={\mu}, \ \ \ \ 
	\frac{d {\bar{t}}}{d \tau}=\gamma,\ \ \ \ 
	\frac{d{\mu}}{d \tau}= {\mu} \times B(x)+\gamma E(x),  \ \ \ \
	\frac{d \gamma}{d\tau}=E(x) \cdot {\mu}.
\end{equation*}
We reformulate the above system in Minkowski spacetime by setting 
$v={\mu}$ and introducing the imaginary variables 
$w=i\gamma$ and $t=i\bar{t}$, thereby yielding the following formulation
\begin{equation}\label{rela-CPD-4d}
	\frac{d x}{d \tau}=v, \ \ \ \ 
	\frac{d t}{d \tau}=w,\ \ \ \ 
	\frac{dv}{d \tau}=v \times B(x)-iw E(x),  \ \ \ \
	\frac{d w}{d\tau}=i E(x)\cdot v.
\end{equation}
The four-dimensional position and momentum vectors are denoted by 
${\bm{y}}=(x;t)=(x^{\intercal},t)^{\intercal}$ and ${\bm{u}}=(v;w)=(v^{\intercal},w)^{\intercal}$, respectively, with boldface letters used to represent 4D vectors and matrices in what follows. The Lagrangian for the above four-dimensional system is $\mathcal{L}(x,t,v,w)=(w^2+v^2)/2+A(x) \cdot v+iwU(x)$
and the conjugate momentum is $\bm{p}=\nabla_{\bm{u}}\mathcal L=\big(\frac{\partial \mathcal L(x,t,v,w)}{\partial v}; \frac{\partial \mathcal L(x,t,v,w)}{\partial w}\big)=\big(v+A(x); w+iU(x)\big)$.  Then the corresponding Hamiltonian  $\mathcal{H}=\bm{p} \cdot \bm{u}-\mathcal L$, which characterizes the mass shell of the charged particle described in \eqref{rela-CPD-4d}, is thus given by
\begin{equation}\label{mass shell}
	\mathcal{H}(v,w)=(w^2+v^2)/2.
\end{equation}

The Hamiltonians $H$ (energy)
and $\mathcal{H}$ (mass shell)
are conserved quantities along the solutions of the equivalent systems \eqref{rela-CPD} and \eqref{rela-CPD-4d}. Moreover, for the system \eqref{rela-CPD-4d}, we denote the phase-space variable by $r=(x^{\intercal},t,v^{\intercal},w)^{\intercal}$, and the associated vector field is $f(x,t,v,w)=\big(v^{\intercal}, w, (v \times B(x)-iw E(x))^{\intercal},i E(x)\cdot v\big)^{\intercal}:=(f_1,f_2,f_3,f_4)$. One readily checks that this vector field is divergence-free, i.e., $\operatorname{div} f(t,x,w,v)=\frac{\partial f_1}{\partial t}+\frac{\partial f_2}{\partial x}+\frac{\partial f_3}{\partial w}+\frac{\partial f_4}{\partial v}=0$. Hence, by Liouville’s theorem (see \cite{Hairer2002}), the phase-space volume is preserved along the solution flow of the system \eqref{rela-CPD-4d}.
The conservation of geometric structures plays a crucial role in ensuring the physical validity of numerical solutions for simulating the motion of charged particles in electromagnetic fields. In this paper, we shall analyze the error estimates and structure-preserving properties of the proposed numerical methods, focusing on the long-term conservation or near-conservation of energy and mass shell, as well as the preservation of phase-space volume.

Various methods for relativistic charged-particle dynamics have been successfully developed in recent years within the field of plasma physics \cite{Y.Wang2016, Matsuyama2017, Y.Wang2021, J.Xiao2019}. The Boris method \cite{Boris 1970}, a classic second-order accurate leapfrog scheme, is one of the most commonly used methods in relativistic and non-relativistic charged-particle dynamics. In addition to the Boris integrator, several other integrators have been developed. An implicit scheme for solving the relativistic equation of motion in an external electromagnetic field was introduced in \cite{Petri2017}, while \cite{Qiang2017} developed a fast Runge-Kutta relativistic integrator to address force balance issues. The study in \cite{Ripperda2018} conducted a comparison of multiple numerical methods and formulated a new fully implicit  particle integrator. An alternative integrator designed to better resolve the $E \times B$ drift at relativistic speeds was presented  in \cite{Vay2008}. 

Structure-preserving algorithms find many applications in different areas of physics, and extensive research has also been carried out in plasma physics \cite{Feng2010,Gauckler2018,Hairer2002,YHe2015,Ricketson2020}.
By preserving important geometric structures such as symplecticity, energy, and phase-space volume, these methods are well known for their long-term accuracy and stability, and they have been shown to be highly effective in addressing multiscale and nonlinear dynamics. Several structure-preserving methods have been investigated for relativistic charged-particle dynamics. For example, volume-preserving integrators have been discussed in \cite{Finn2005,Y.He2016,A.V.Higuera2017,R.Zhang2015}, and symplectic methods have been explored in \cite{Wolski2018,Wu2003,R.Zhang2018}.
Furthermore, a recent study \cite{Hairer2023} constructed some leapfrog methods and analyzed their properties in terms of geometric structure preservation. 

In addition, splitting methods \cite{McLachlan2002} enhance numerical efficiency in geometric integration by decomposing systems into simpler subsystems, and have been extensively applied in plasma physics \cite{Chartier2016,Crouseilles2015,Lubich2008,WangB2021}.  A structure-preserving algorithm based on a splitting scheme was recently proposed \cite{zhang-wang 2024}, which attains first-order uniform accuracy for relativistic charged-particle dynamics under the scaling $B(x)=B(\epsilon x)/\epsilon$. Concerning the numerical methods for solving relativistic charged-particle dynamics \eqref{rela-CPD}, to the best of our knowledge, the second-order splitting scheme proposed in this work has not yet been analyzed in terms of its geometric structure-preserving properties.
The objective of this article is to present four novel two-step symmetric methods based on this splitting scheme, which achieve:
\begin{description}
	\item[--] Second-order accuracy;
	\item[--] Exact conservation of the mass shell and phase-space volume;
	\item[--] Long-time near-conservation of energy and mass shell.
\end{description}

The organization of this paper is as follows.
In Section \ref{sec:methods}, we introduce four two-step symmetric methods based on a Strang splitting scheme.
Section \ref{results and experiments} presents the main results concerning the convergence and the long-time conservation or near-conservation of geometric structures, together with five numerical experiments. Section \ref{Error and long-term analysis} is devoted to rigorous proofs of the convergence results as well as the long-term conservation or near-conservation properties of energy, mass shell, and volume for the proposed methods.
Concluding remarks are given in Section \ref{sec:conclusions}.

\section{Numerical integrators}
\label{sec:methods}%
This section introduces four two-step symmetric methods designed to solve the four-dimensional relativistic charged-particle dynamics \eqref{rela-CPD-4d}. Notably,   a two-step method denoted by $y_{n+1}=\Phi_{h}(y_{n-1},y_n)$ is called symmetric if exchanging $y_{n-1} \leftrightarrow y_{n+1}$ and $h \leftrightarrow  -h$ leaves the method unaltered. As pointed out in \cite{Hairer2002}, symmetric methods exhibit excellent long-term behavior and play a central role in the geometric integration of differential equations. Based on the splitting approach, we first present the framework of two-step symmetric methods, as described in the following subsection.
\subsection{Splitting two-step symmetric framework}
To begin with, the relativistic dynamics \eqref{rela-CPD-4d} can be rewritten as
\begin{equation}\label{CPD-4d}
	\small{	
		\frac{d}{d\tau}\begin{array}[c]{cc}
			\left(
			\begin{array}{cc}
				x\\
				t\\
			\end{array}
			\right)=\left(
			\begin{array}{cc}
				v\\
				w\\
			\end{array}
			\right)\end{array},\qquad 
		\frac{d}{d\tau}\begin{array}{cc}
			\left(
			\begin{array}{cc}
				v\\
				w\\
			\end{array}
			\right)={\bm{F}}{(x)}\left(
			\begin{array}{cc}
				v\\
				w\\
			\end{array}
			\right),
	\end{array}}
\end{equation}
where
\begin{equation}\label{F}
		\begin{aligned}
			{\bm{F}}( x)=	\left(
			\begin{array}{ccc;{1pt/3pt}c}
				0 & B_3( x) & -B_2( x) & -iE_1( x) \\
				-B_3( x) & 0 & B_1( x) & -iE_2( x) \\
				B_2( x) & -B_1( x) & 0 & -iE_3( x) \\
				\hdashline[1pt/3pt]
				iE_1( x) & iE_2( x) & iE_3( x) & 0 \\
			\end{array}
			\right):=\left(
			\begin{array}{c;{1pt/3pt}c}
				\tilde{B}( x)& -iE( x)    \\
				\hdashline[1pt/3pt]
				iE( x) ^{\intercal} & 0\\
			\end{array}
			\right)
	\end{aligned}
\end{equation}
is a skew-symmetric matrix and $\tilde{B}( x)$ satisfies $\tilde{B}( x)\alpha = \alpha \times B( x)$ for all $\alpha \in \mathbb{R}^3$. The equations \eqref{CPD-4d} are then split into the following two subsystems
\begin{equation*}\label{}
	\small{\setlength{\arraycolsep}{2pt}
		\begin{aligned}
			&Z_1:\quad 	\frac{d}{d\tau}\begin{array}[c]{cc}
				\left(
				\begin{array}{cc}
					x\\
					t\\
				\end{array}
				\right)=0\end{array},\quad\qquad 
			\frac{d}{d\tau}\begin{array}{cc}
				\left(
				\begin{array}{cc}
					v\\
					w\\
				\end{array}
				\right)={\bm{F}}{(x)}\left(
				\begin{array}{cc}
					v\\
					w\\
				\end{array}
				\right),
			\end{array}
			\\
			&	Z_2:\quad 	\frac{d}{d\tau}\begin{array}[c]{cc}
				\left(
				\begin{array}{cc}
					x\\
					t\\
				\end{array}
				\right)=\left(
				\begin{array}{cc}
					v\\
					w\\
				\end{array}
				\right)\end{array},\quad 
			\frac{d}{d\tau}\begin{array}{cc}
				\left(
				\begin{array}{cc}
					v\\
					w\\
				\end{array}
				\right)=0.
			\end{array}
	\end{aligned}}
\end{equation*}
Both subsystems can be solved exactly, and their solution flows $\Psi_{h}^{i}$ $(i=1,2)$ have explicit representations given by
\begin{equation*}\label{flows}
	\begin{aligned}
		&\Psi_{h}^{1}:
		\begin{array}[c]{ll}\left(
			\begin{array}{cc}
				x_{n+1}\\
				{t}_{n+1}\\
			\end{array}
			\right)=\left(
			\begin{array}{cc}
				x_{n}\\
				{t}_{n}\\
			\end{array}
			\right),\quad \left(
			\begin{array}{cc}
				v_{n+1}\\
				{w}_{n+1}\\
			\end{array}
			\right)=e^{h{\bm{F}}(x_{n})}\left(
			\begin{array}{cc}
				v_{n}\\
				{w}_{n}\\
			\end{array}
			\right),
		\end{array}\\[0.3em]
		&\Psi_{h}^{2}:
		\begin{array}[c]{ll}\left(
			\begin{array}{cc}
				x_{n+1}\\
				{t}_{n+1}\\
			\end{array}
			\right)=\left(
			\begin{array}{cc}
				x_{n}\\
				{t}_{n}\\
			\end{array}
			\right)+h \left(
			\begin{array}{cc}
				v_{n}\\
				{w}_{n}\\
			\end{array}
			\right),\quad \left(
			\begin{array}{cc}
				v_{n+1}\\
				{w}_{n+1}\\
			\end{array}
			\right)=\left(
			\begin{array}{cc}
				v_{n}\\
				{w}_{n}\\
			\end{array}
			\right),
		\end{array}
	\end{aligned}
\end{equation*}
where $h$ denotes the chosen stepsize. In order to solve the relativistic dynamics \eqref{rela-CPD-4d}, we employ a Strang splitting scheme (see \cite{Hairer2002}), given by $	\Psi_{h}=\Psi_{h/2}^{1} \circ \Psi_{h}^{2} \circ \Psi_{h/2}^{1}.$
Then an exponential splitting scheme
is derived as
\begin{equation}\label{ori-exp}
	\left(
	\begin{array}{cc}
		x_{n+1}\\
		{t}_{n+1}\\
	\end{array}
	\right)=\left(
	\begin{array}{cc}
		x_{n}\\
		{t}_{n}\\
	\end{array}
	\right)+he^{\frac{h}{2}{\bm{F}}(x_{n})}\left(
	\begin{array}{cc}
		v_{n}\\
		{w}_{n}\\
	\end{array}
	\right),\ 
	\left(
	\begin{array}{cc}
		v_{n+1}\\
		{w}_{n+1}\\
	\end{array}
	\right)=e^{\frac{h}{2}{\bm{F}}(x_{n+1})}e^{\frac{h}{2}{\bm{F}}(x_{n})}\left(
	\begin{array}{cc}
		v_{n}\\
		{w}_{n}\\
	\end{array}
	\right).
\end{equation}

By introducing the notations ${\bm{y}}_{n}=(x_{n};{t}_{n})$ and  ${\bm{u}}_{n}=(v_{n};w_n)$, the above equations can be rewritten as 
\begin{equation}\label{M1-ori}
	{\bm{y}}_{n+1}={\bm{y}}_{n}+he^{\frac{h}{2}{\bm{F}}(x_{n})}{\bm{u}}_{n},\quad 
	{\bm{u}}_{n+1}=e^{\frac{h}{2}{\bm{F}}(x_{n+1})}e^{\frac{h}{2}{\bm{F}}(x_{n})}{\bm{u}}_{n}.
\end{equation}
Using the property of the symmetric method, it is easy to confirm that $\Psi_{h}$ is symmetric.
As a result, we can conclude that ${\bm{y}_{n-1}}={\bm{y}_{n}}-he^{-\frac{h}{2}{\bm{F}}(x_n)}{\bm{u}_{n}}.$
Combining  the  first equation of \eqref{M1-ori} with this expression yields
\begin{equation}\label{ori-meth}
		\begin{aligned}
			&\frac{{\bm{y}_{n+1}}-2{\bm{y}_{n}}+{\bm{y}_{n-1}}}{h^2}
			=\frac{2}{h}\frac{e^{\frac{h}{2}{\bm{F}}({x_n})}-e^{-\frac{h}{2}{\bm{F}}(x_n)}}{e^{\frac{h}{2}{\bm{F}}({x_n})}+e^{-\frac{h}{2}{\bm{F}}(x_n)}}	\frac{{\bm{y}_{n+1}}-{\bm{y}_{n-1}}}{2h}\\
			&\qquad\qquad\qquad\qquad\quad =\frac{2}{h}\tanh\Big(\frac{h}{2}{\bm{F}}({x_{n}})\Big)\frac{{\bm{y}_{n+1}}-{\bm{y}_{n-1}}}{2h},\\[0.3em]
			&{\bm{u}}_{n}=\frac{2}{e^{\frac{h}{2}{\bm{F}}({x_n})}+e^{-\frac{h}{2}{\bm{F}}(x_n)}}	\frac{{\bm{y}_{n+1}}-{\bm{y}_{n-1}}}{2h}=\operatorname{sech}\Big(\frac{h}{2}{\bm{F}}({x_{n}})\Big)\frac{{\bm{y}_{n+1}}-{\bm{y}_{n-1}}}{2h}.
		\end{aligned}
	\end{equation}
To complete the above formulation, we derive the explicit expressions of the matrix functions $\tanh(h\bm{F})$ and $\operatorname{sech}(h\bm{F})$ in the subsequent subsection.
	\subsection{Coefficient derivations for \texorpdfstring{$\tanh(h\bm{F})$}{tanh(hF)} and \texorpdfstring{$\operatorname{sech}(h\bm{F})$}{sech(hF)}}
	We now present the detailed derivation of $\tanh(h{\bm{F}})$ and $\operatorname{sech}(h{\bm{F}})$, with the argument $x_n$ omitted for clarity. As shown in \cite{zhang-wang 2024}, the definition of ${\bm{F}}$ in \eqref{F} yields the following relations:
	\begin{equation}\label{F3}
		{\bm{F}}^3=r_1{\bm{F}}+r_2\hat{{\bm{F}}}, \ \ {\bm{F}}\hat{{\bm{F}}}=r_2{\bm{I}},
	\end{equation}
	with $r_1=E^{\intercal}E-B^{\intercal}B$,  $r_{2}=-E^{\intercal} B$ and  
	\begin{equation*}\label{hatF} 
		\small{
			\begin{aligned}
				\hat{{\bm{F}}}=\left(
				\begin{array}{ccc;{1pt/3pt}c}
					0 & E_3 & -E_2 & iB_1 \\
					-E_3 & 0 & E_1 & iB_2 \\
					E_2 & -E_1 & 0 & iB_3 \\
					\hdashline[1pt/3pt]
					-iB_1 & -iB_2 & -iB_3 & 0 \\
				\end{array}
				\right):=\left(
				\begin{array}{c;{1pt/3pt}c}
					\tilde{E} & iB    \\
					\hdashline[1pt/3pt]
					-iB^{\intercal} & 0\\
				\end{array}
				\right).
		\end{aligned}}
	\end{equation*}
	Applying the Taylor expansion yields $	\tanh(h{\bm{F}})=\sum_{k=1}^{\infty}\frac{2^{2k}(2^{2k}-1)\mathcal{B}_{2k}h^{2k-1}}{(2k)!} {\bm{F}}^{2k-1}$ and $\operatorname{sech}(h{\bm{F}})={\bm{I}}+\sum_{k=1}^{\infty}\frac{\mathcal{E}_{2k}h^{2k}}{(2k)!} {\bm{F}}^{2k}$,
	where ${\bm{I}}$ denotes the 4-dimensional identity matrix, and $\mathcal{B}_{2k}$ and $\mathcal{E}_{2k}$ denote the Bernoulli and Euler numbers, respectively.  Using the relations \eqref{F3},  these series can be reformulated as
	\begin{equation*}
		\begin{aligned}
			&	\tanh(h{\bm{F}})
			=\sum_{k=1}^{\infty}\frac{2^{2k}(2^{2k}-1)\mathcal{B}_{2k}h^{2k-1}}{(2k)!}\big(a_{2k-1}{\bm{F}}+b_{2k-1}\hat{{\bm{F}}}\big),\\
			&\operatorname{sech}(h{\bm{F}})
			={\bm{I}}+\sum_{k=1}^{\infty}\frac{\mathcal{E}_{2k}h^{2k}}{(2k)!}\big(a_{2k}{\bm{I}}+b_{2k}{\bm{F}}^2\big).
		\end{aligned}
	\end{equation*}
	The coefficients $a_{2k-1}$,  $b_{2k-1}$, $a_{2k}$ and $b_{2k}$ satisfy the following recurrence relations	
	\begin{equation}\label{recursions}
		\begin{aligned}
			&	a_{2k-1}=r_1a_{2k-3}+r_2b_{2k-3},\ \
			b_{2k-1}=r_2a_{2k-3},\ \
			a_{1}=1,\ \  b_{1}=0,\\
			&  a_{2k+2}=r_2^2b_{2k},\ \ \ \
			b_{2k+2}=a_{2k}+r_1b_{2k},\ \
			a_{2}=0,\ \  b_{2}=1.
		\end{aligned}
	\end{equation}
	Define the coefficient functions by
	\begin{equation*}\label{}
		\begin{aligned}
			&K_1=\sum_{k=1}^{\infty}\frac{2^{2k}(2^{2k}-1)\mathcal{B}_{2k}h^{2k-1}}{(2k)!}a_{2k-1},\quad 
			K_2=\sum_{k=1}^{\infty}\frac{2^{2k}(2^{2k}-1)\mathcal{B}_{2k}h^{2k-1}}{(2k)!}b_{2k-1},\\
			&K_3=1+\sum_{k=1}^{\infty}\frac{\mathcal{E}_{2k}h^{2k}}{(2k)!}a_{2k},\qquad \qquad \qquad \ \
			K_4=\sum_{k=1}^{\infty}\frac{\mathcal{E}_{2k}h^{2k}}{(2k)!}b_{2k}.
		\end{aligned}
	\end{equation*}
	Then, the matrix functions $\tanh(h\bm{F})$ and $\operatorname{sech}(h\bm{F})$ can be expressed as
	\begin{equation}\label{tanhF}
		\begin{aligned}
			\tanh(h{\bm{F}})
			=K_1{\bm{F}}+K_2\hat{{\bm{F}}},\quad \
			\operatorname{sech}(h{\bm{F}})=K_3{\bm{I}}+K_4{{\bm{F}}}^2.
		\end{aligned}
	\end{equation}

Explicit derivations are presented for $K_1$ and $K_3$, while the expressions for $K_2$ and $K_4$ follow by analogy. From the recurrence relations of \eqref{recursions}, it follows that
	\begin{equation*}\label{}
		\small{\left(
			\begin{array}{cc}
				a_{2k-1}\\
				b_{2k-1}\\
			\end{array}
			\right)=\left(
			\begin{array}{cc}
				r_1 & r_2    \\
				r_2 & 0\\
			\end{array}
			\right)\left(
			\begin{array}{cc}
				a_{2k-3}\\
				b_{2k-3}\\
			\end{array}
			\right):=P_1\left(
			\begin{array}{cc}
				a_{2k-3}\\
				b_{2k-3}\\
			\end{array}
			\right),\ \ \left(
			\begin{array}{cc}
				a_{1}\\
				b_{1}\\
			\end{array}
			\right)=\left(
			\begin{array}{cc}
				1\\
				0\\
			\end{array}
			\right),}
	\end{equation*}
	and
	\begin{equation*}\label{}
		\left(
		\begin{array}{cc}
			a_{2k+2}\\
			b_{2k+2}\\
		\end{array}
		\right)=\left(
		\begin{array}{cc}
			0 & r_{2}^{2}    \\
			1 & r_1\\
		\end{array}
		\right)\left(
		\begin{array}{cc}
			a_{2k}\\
			b_{2k}\\
		\end{array}
		\right):=P_2\left(
		\begin{array}{cc}
			a_{2k}\\
			b_{2k}\\
		\end{array}
		\right),\ \ \left(
		\begin{array}{cc}
			a_{2}\\
			b_{2}\\
		\end{array}
		\right)=\left(
		\begin{array}{cc}
			0\\
			1\\
		\end{array}
		\right).
	\end{equation*}
	These relations then yield $(a_{2k-1}, b_{2k-1})^{\intercal}=P_1^{k-1}(1, 0)^{\intercal}$ and $(a_{2k+2},b_{2k+2})^{\intercal}=P_2^{k}(0, 1)^{\intercal}$.  Diagonalization of $P_i$ $(i=1,2)$ yields
	\begin{equation*}\label{}
		\begin{aligned}
			&	a_{2k-1}=(\lambda_{1}^k-\lambda_{2}^k)/(\lambda_1-\lambda_2),\quad \
			b_{2k-1}=r_2(\lambda_{1}^{k-1}-\lambda_{2}^{k-1})/(\lambda_1-\lambda_2),\\
			&a_{2k}=r_2^2(\lambda_{1}^{k-1}-\lambda_{2}^{k-1})/(\lambda_1-\lambda_2),\quad \ 
			b_{2k}=(\lambda_{1}^{k}-\lambda_{2}^{k})/(\lambda_1-\lambda_2).
		\end{aligned}
	\end{equation*}
	Both matrices $P_{1}$ and $P_{2}$ share the same eigenvalues, given by
	$\lambda_{1}=(r_1+d)/2$ and $\lambda_{2}=(r_1-d)/2$, where $d=\sqrt{r_1^2+4r_2^2}$. Substituting the expressions for $a_{2k-1}$ and $a_{2k}$ into the definitions of $K_1$ and $K_3$, respectively, we obtain
	\begin{equation*}\label{}
		\begin{aligned}
			&K_1=\sum_{k=1}^{\infty}\frac{2^{2k}(2^{2k}-1)\mathcal{B}_{2k}h^{2k-1}(\lambda_{1}^k-\lambda_{2}^k)}{(2k)!(\lambda_1-\lambda_2)},\ \
			K_3=1+\sum_{k=1}^{\infty}\frac{\mathcal{E}_{2k}h^{2k}r_2^2(\lambda_{1}^{k-1}-\lambda_{2}^{k-1})}{(2k)!(\lambda_1-\lambda_2)}.\\
		\end{aligned}
	\end{equation*}
	From the Taylor expansions of $\tanh(\cdot)$ and $\operatorname{sech}(\cdot)$, one verifies that 
	\begin{equation*}
		\begin{aligned}
			&\sum_{k=1}^{\infty}\frac{2^{2k}(2^{2k}-1)\mathcal{B}_{2k}h^{2k-1}\lambda_{1,2}^k}{(2k)!}=\sqrt{\lambda_{1,2}}\tanh\big(h\sqrt{\lambda_{1,2}}\big),\\
			&\sum_{k=1}^{\infty}\frac{\mathcal{E}_{2k}h^{2k}\lambda_{1,2}^{k-1}}{(2k)!}=\frac{1}{\lambda_{1,2}}\big(\operatorname{sech}\big(h\sqrt{\lambda_{1,2}}\big)-1\big).
		\end{aligned}
	\end{equation*}
	Thus, the explicit expressions for $K_{1}$ and $K_{3}$ are given by
	\begin{equation*}
		\begin{aligned}
			&K_{1}=\big(\sqrt{\lambda_1}\tanh\big(h\sqrt{\lambda_1}\big)-\sqrt{\lambda_2}\tanh\big(h\sqrt{\lambda_2}\big)\big)/d,\\
			&K_3=(1-r_1/d)\operatorname{sech}\big(h\sqrt{\lambda_1}\big)/2+(1+r_1/d)\operatorname{sech}\big(h\sqrt{\lambda_2}\big)/2.\\
		\end{aligned}
	\end{equation*}
	Similarly, the explicit formulas for $K_{2}$ and $K_{4}$ are obtained as
	\begin{equation*}
		\begin{aligned}
			&K_{2}=r_2\big(\tanh\big(h\sqrt{\lambda_1}\big)/\sqrt{\lambda_1}-\tanh\big(h\sqrt{\lambda_2}\big)/\sqrt{\lambda_2}\big)/d,\\
			&K_4=\big(\operatorname{sech}\big(h\sqrt{\lambda_1}\big)-\operatorname{sech}\big(h\sqrt{\lambda_2}\big)\big)/d.
		\end{aligned}
	\end{equation*}
We are now in a position to present the explicit formulations of the four two-step symmetric methods, as shown in the next subsection.
	\subsection{Formulations of four two-step symmetric methods}
Building on the coefficient derivations in the previous subsection, substituting $\tanh(h\bm{F})$ and $\operatorname{sech}(h\bm{F})$ from \eqref{tanhF} into \eqref{ori-meth} yields
	\begin{equation*}\label{two-step}
		\begin{aligned}
			&\frac{{\bm{y}_{n+1}}-2{\bm{y}_{n}}+{\bm{y}_{n-1}}}{h^2}
			=\frac{2}{h}\big(K_1{\bm{F}}+K_2\hat{{\bm{F}}}\big)\frac{{\bm{y}_{n+1}}-{\bm{y}_{n-1}}}{2h},\\
			&{\bm{u}}_{n}=\big(K_3{\bm{I}}+K_4{{\bm{F}}}^2\big)\frac{{\bm{y}_{n+1}}-{\bm{y}_{n-1}}}{2h}.
		\end{aligned}
	\end{equation*}
	Recalling the notations ${\bm{y}}_{n}$, ${\bm{u}}_{n}$, and the definitions of ${\bm{F}}$ and $\hat{{\bm{F}}}$, we have
	\begin{equation*}
			\begin{aligned}
				&\left(
				\begin{array}{cc}
					\frac{{x_{n+1}}-2{x_{n}}+{x_{n-1}}}{h^2}\\
					\frac{{t_{n+1}}-2{t_{n}}+{t_{n-1}}}{h^2}
				\end{array}
				\right)
				=\frac{2}{h}		\left(
				\begin{array}{cc}
					K_1\tilde{B}+K_2\tilde{E} & -i\big(K_1E-K_2B\big)   \\
					i\big(K_1E^{\intercal}-K_2B^{\intercal}\big) & 0\\
				\end{array}
				\right)\left(
				\begin{array}{cc}
					\frac{{x_{n+1}}-{x_{n-1}}}{2 h}\\
					\frac{{t_{n+1}}-{t_{n-1}}}{2 h}\\
				\end{array}
				\right),\\
				&\left(
				\begin{array}{cc}
					v_{n}\\
					{w}_{n}\\
				\end{array}
				\right)=\left(
				\begin{array}{cc}
					K_3I+K_4(\tilde{B}^2+EE^{\intercal}) & -iK_4\tilde{B}E  \\
					iK_4E^{\intercal}\tilde{B}& K_3+K_4E^{\intercal}E\\
				\end{array}
				\right)\left(
				\begin{array}{cc}
					\frac{{x_{n+1}}-{x_{n-1}}}{2 h}\\
					\frac{{t_{n+1}}-{t_{n-1}}}{2 h}\\
				\end{array}
				\right).
			\end{aligned}
		\end{equation*}
With these preparations, four two-step symmetric methods can be rigorously formulated. The first method is stated as follows.
		\begin{definition} \label{}
			A  two-step symmetric  method for solving the relativistic charged-particle dynamics \eqref{rela-CPD-4d} is
			defined as
			\begin{equation}\label{M1}
				\begin{aligned}
					&\frac{{x_{n+1}}-2{x_{n}}+{x_{n-1}}}{h^2}
					=\frac{2}{h}\big(K_1({x_{n}})\tilde{B}( {x_{n}})+K_2({x_{n}})\tilde{E}( {x_{n}})\big)\frac{{x_{n+1}}-{x_{n-1}}}{2 h}\\
					&\qquad\qquad\qquad\qquad\qquad\ \ -\frac{2}{h}i\big(K_1({x_{n}})E( {x_{n}})-K_2({x_{n}})B( {x_{n}})\big)\frac{t_{n+1}-t_{n-1}}{2 h},\\
					&	\frac{{t_{n+1}}-2{t_{n}}+{t_{n-1}}}{h^2}
					= \frac{2}{h}i\big(K_1({x_{n}})E( {x_{n}})^{\intercal}-K_2({x_{n}})B( {x_{n}})^{\intercal}\big)\frac{{x_{n+1}}-{x_{n-1}}}{2 h}.
				\end{aligned}
			\end{equation} 
			The momentum vector $(v;w)$ is approximated by
			$$(v_{n+1};w_{n+1})=e^{\frac{h}{2}{\bm{F}}(x_{n+1})}e^{\frac{h}{2}{\bm{F}}(x_{n})}(v_{n}; w_{n}).$$
			We shall call it as \textbf{M1}. 
		\end{definition}
		
		By modifying the $(x;t)$ components of M1 and introducing a new approximation for the momentum vector $(v;w)$, we construct three additional two-step symmetric methods, denoted by M2--M4.
		\begin{definition} \label{def:M2--M4}
			We present three two-step symmetric methods (referred to as \textbf{M2}--\textbf{M4}) for solving the relativistic charged-particle dynamics \eqref{rela-CPD-4d}. The detailed formulations are as follows:
			\begin{itemize}
				\item \textbf{M2}: The first method is defined by
				\begin{equation}\label{M2}
					\begin{aligned}
						&\frac{x_{n+1} - 2x_n + x_{n-1}}{h^2}
						= \frac{2}{h} \big(K_1(x_0) \tilde{B}(x_n) + K_2(x_0) \tilde{E}(x_n)\big) \frac{x_{n+1} - x_{n-1}}{2h} \\
						&\qquad\qquad\qquad \qquad \qquad\ - \frac{2i}{h} \big(K_1(x_0) E(x_n) - K_2(x_0) B(x_0)\big) \frac{t_{n+1} - t_{n-1}}{2h}, \\
						&\frac{t_{n+1} - 2t_n + t_{n-1}}{h^2}
						= \frac{2i}{h} \big(K_1(x_0) E(x_n)^\intercal - K_2(x_0) B(x_0)^\intercal \big) \frac{x_{n+1} - x_{n-1}}{2h}.
					\end{aligned}
				\end{equation}
				A new approximation for the momentum vector $(v;w)$ is given by
				\begin{equation*}\label{}
					\begin{aligned}
						&v_n
						=\big(K_3({x_{0}})I+K_4({x_{0}})\big(\tilde{B}( {x_{n}})^2+E( {x_{n}})E( {x_{n}})^{\intercal}\big)\big)\frac{{x_{n+1}}-{x_{n-1}}}{2 h}\\
						&\qquad\ \ -iK_4({x_{0}})\tilde{B}( {x_{n}})E( {x_{n}})\frac{t_{n+1}-t_{n-1}}{2 h},\\
						&	w_n
						= iK_4({x_{0}})E( {x_{n}})^{\intercal}\tilde{B}( {x_{n}})\frac{{x_{n+1}}-{x_{n-1}}}{2 h}+\big(K_3({x_{0}})+K_4({x_{0}})E( {x_{n}})^{\intercal}E( {x_{n}})\big)\frac{t_{n+1}-t_{n-1}}{2 h}.
					\end{aligned}
				\end{equation*} 
				
				\item \textbf{M3}: Based on the Taylor expansions 
				$K_2(x_0) = -\frac{h^3}{24} r_2(x_0) + \frac{h^5}{240} r_1(x_0) r_2(x_0) + \mathcal{O}(h^7)$ and $K_4(x_0) = -\frac{h^2}{8}+ \frac{5h^4}{384} r_1(x_0)+ \mathcal{O}(h^6)$,  we introduce the second method as follows
				\begin{equation*}\label{M3}
					\begin{aligned}
						&\frac{x_{n+1} - 2x_n + x_{n-1}}{h^2}
						= \frac{2}{h} K_1(x_0) \Big(\tilde{B}(x_n) \frac{x_{n+1} - x_{n-1}}{2h} - i E(x_n) \frac{t_{n+1} - t_{n-1}}{2h} \Big), \\
						&\frac{t_{n+1} - 2t_n + t_{n-1}}{h^2}
						= \frac{2i}{h} K_1(x_0) E(x_n)^\intercal \frac{x_{n+1} - x_{n-1}}{2h}.
					\end{aligned}
				\end{equation*} 
				The approximation for the momentum vector $(v;w)$ is given by
				$$
				v_n
				=K_3({x_{0}})(x_{n+1}-x_{n-1})/(2h),\quad \ w_n=K_3({x_{0}})(t_{n+1}-t_{n-1})/(2h).
				$$
				
				\item \textbf{M4}: Using the Taylor expansions 
				$K_1(x_0) = \frac{h}{2} - \frac{h^3}{24} r_1(x_0) +\mathcal{O}(h^5)$ and $K_3(x_0) =1+ \frac{5h^4}{384} r_2^2(x_0)+ \mathcal{O}(h^6)$, the third method is formulated as 
				\begin{equation*}\label{M4}
					\begin{aligned}
						&\frac{x_{n+1} - 2x_n + x_{n-1}}{h^2}
						= \tilde{B}(x_n) \frac{x_{n+1} - x_{n-1}}{2h} - i E(x_n) \frac{t_{n+1} - t_{n-1}}{2h}, \\
						&\frac{t_{n+1} - 2t_n + t_{n-1}}{h^2}
						= i E(x_n)^\intercal \frac{x_{n+1} - x_{n-1}}{2h}.
					\end{aligned}
				\end{equation*}
				The approximation for the momentum vector $(v;w)$ is given by
				$$
				v_n
				=(x_{n+1}-x_{n-1})/(2h),\quad \ w_n=(t_{n+1}-t_{n-1})/(2h).
				$$
			\end{itemize}
			To simplify the expressions in the following sections, we use the abbreviations
			$S_1=K_1({x_{0}})$, $S_2=K_2({x_{0}})$, $S_3=K_3({x_{0}})$, $S_4=K_4({x_{0}})$ and $B_0=B( {x_{0}})$.
		\end{definition}
		
		\begin{remark}\label{remark1}
The construction of M2--M4 is motivated by two primary considerations.
First, we aim to preserve the second-order accuracy of M1, which is confirmed by the error analysis in Section~\ref{error analysis}. 
Second, M1 exhibits energy drift, as observed in Problem~2 of Section~\ref{experiments}. To address this issue, we modify M1 by treating $K_i(x_n)$ $(i=1,2)$ and $B(x_n)$ in the $(x;t)$ components as the constants $K_i(x_0)$ $(i=1,2)$ and $B(x_0)$, and introducing a new momentum approximation, yielding M2, which achieves long-term near-conservation of energy.  This method, along with its simplified versions M3 and M4, successfully achieves long-term near-conservation of energy. The detailed proof is provided in Section~\ref{energy analysis}.
		\end{remark}

		\begin{remark}
			It is worth noting that the starting value $({x_{1}};t_1)$ is required for implementing the above two-step symmetric methods. This value can be determined by setting $n = 0$ in the first equation of \eqref{ori-exp}. Moreover, it is well known that implicit methods can be implemented by fixed-point iteration or by reformulating them as explicit schemes through collecting all terms involving $x_{n+1}$
			and inverting the resulting expression. The numerical experiments presented in this paper employ the latter approach.
		\end{remark}

		\section{Main results and numerical experiments}\label{results and experiments}
		This section presents the main results on error bounds and long-time conservation or near-conservation of key geometric structures for the proposed numerical methods M1--M4, including energy, mass shell, and phase-space volume.  The first part is concerned with theoretical results,  and the second focuses on the numerical experiments.
		
		\subsection{Main results}\label{results}
		This subsection investigates the error bounds and the long-time conservation or near-conservation of three geometric structures.
		The following theorem establishes the error bounds, demonstrating the second-order accuracy of the proposed methods. 
		\begin{theorem}\label{error bound2}\textbf{(Convergence)}
			Suppose that the relativistic charged-particle dynamics \eqref{rela-CPD-4d} has sufficiently smooth solutions and the functions $B(x)$ and $E(x)$ are sufficiently  differentiable. We further assume that $B(x)$ and $E(x)$ are locally Lipschitz continuous with a Lipschitz constant $L$.
			If the step size $h$ is chosen such that $ h \leq h_0$
			for a constant $h_0>0$, then applying the two-step symmetric methods M1--M4 to \eqref{rela-CPD-4d} yields the following global errors
			\begin{equation*}\label{R1}
				\abs{{\bm{y}}_{n+1}-{\bm{y}}(\tau_{n+1})} \leq Ch^2, \ \quad	\abs{{\bm{u}}_{n+1}-{\bm{u}}(\tau_{n+1})} \leq Ch^2,
			\end{equation*}
			for $(n+1)h\leq T$, where $C>0$ is a generic constant independent of $h$ or $n$ but depends on $L$ and $T$.
		\end{theorem}
		
		The following theorem establishes the near-conservation of total energy over long times by the two-step symmetric methods M2--M4 under a quadratic potential. 
		
		\begin{theorem}\label{conservation-H1}\textbf{(Near-conservation of energy)}
			Consider the two-step symmetric methods M2--M4 for solving  the relativistic charged-particle dynamics \eqref{rela-CPD-4d} with a quadratic potential $U(x)=\frac{1}{2}x^{\intercal}Qx+q^{\intercal}x$, where
			$Q$ is a symmetric matrix.
			Assume that the numerical solution remains within a compact set independent of  $h$. Then, the energy \eqref{E} is nearly preserved as follows:
			\begin{equation*}
				\abs{H(x_{n},\gamma_n)-	H(x_{0},\gamma_0)} \leq Ch^2 \quad \ \text{for} \quad nh \leq h^{-N+2},
			\end{equation*}			
			where  $N \geq 3$ is the truncation number and  $C$ is a constant that is independent of $n$ and $h$, but depends on the exponent $N$.
		\end{theorem}
		
		The good long-term behavior of the mass shell \eqref{mass shell} along the numerical solutions M1--M4 is established in the following theorem. 
		\begin{theorem}\label{conservation-H2}\textbf{(Conservation or near-conservation of mass shell)}
			Assume that the numerical solution belongs to a compact set independent of $h$,  the following results on the conservation or near-conservation of the mass shell \eqref{mass shell} hold:
			\begin{enumerate}[label=(\roman*)] \item \textbf{Conservation by M1:}  M1 exactly preserves the mass shell, i.e.,
				\begin{equation}\label{MS0}
					\mathcal{H}(v_{n+1},w_{n+1})=\mathcal{H}(v_n,w_n).
				\end{equation}
				
				\item \textbf{Near-conservation by M2--M4:} M2--M4 nearly preserve the mass shell as follows:
				\begin{equation}\label{MS2}
					\left| \mathcal{H}(v_n, w_n) - \mathcal{H}(v_0, w_0) \right| \leq Ch^2 \quad \ \text{for} \quad nh \leq h^{-N+2}.
				\end{equation}
			\end{enumerate}
			Here, $N \geq 3$ is the truncation number, and $C$ is a constant independent of $n$ and $h$, but depends on $N$.
		\end{theorem}	
		
		The previous two theorems have established that the proposed methods ensure long-term conservation or near-conservation of energy and mass shell. The following theorem further confirms the preservation of volume.
		\begin{theorem}\label{conservation-vol}
			\textbf{(Conservation of phase-space volume for M1)}  
			The discrete flow $\Psi_h$ of the two-step symmetric method M1 preserves the phase-space volume: for every bounded open set $\Omega \subset \mathbb{R}^8$, one has 
			$\operatorname{vol}(\Psi_h(\Omega)) = \operatorname{vol}(\Omega)$.
		\end{theorem}

		The four theorems above demonstrate that the proposed methods possess second-order accuracy and maintain the long-term conservation or near-conservation of energy, mass shell and volume. For clarity, these properties are summarized in Table \ref{tab:geo-properties} below.
		\setlength{\arrayrulewidth}{1pt}
		\renewcommand{\arraystretch}{1.2}
		\begin{table}[htbp]
			\centering
			\normalsize
			\resizebox{\textwidth}{!}{%
				\begin{tabular}{|l|l|l|l|l|}
					\hline
					Method &Error bounds& \hspace{1cm}Energy near-conservation & 
					\begin{tabular}{@{}c@{}}
						Mass Shell \\ \hline
						\begin{tabular}{@{}l l@{}}
							Conservation & Near-Conservation
						\end{tabular}
					\end{tabular} & 
					Volume Conservation \\
					\hline
					\hspace{0.3cm}M1 & \hspace{0.3cm} $\mathcal{O}(h^2)$&
					\checkmark (constant magnetic field-\textbf{\textit{numerical result}})
					& \begin{tabular}{@{}l @{}l}
						\hspace{0.65cm}	\checkmark\hspace{2.3cm} & -
					\end{tabular}  & \hspace{1.3cm}\checkmark \\ \hline
					\hspace{0.3cm}M2& \hspace{0.4cm}$\mathcal{O}(h^2)$& \hspace{1.3cm}\checkmark (quadratic potential)  & \begin{tabular}{@{}l l@{}}
						\hspace{0.8cm}- \hspace{1.85cm} & \checkmark
					\end{tabular}  &  \hspace{1.35cm}- \\ \hline
					\hspace{0.3cm}M3 &\hspace{0.4cm}$\mathcal{O}(h^2)$ &  \hspace{1.3cm}\checkmark (quadratic potential)  & \begin{tabular}{@{}l l@{}}
						\hspace{0.8cm}- \hspace{1.85cm} & \checkmark 
					\end{tabular}  & \hspace{1.25cm} - \\ \hline
					\hspace{0.3cm}M4 &\hspace{0.4cm}$\mathcal{O}(h^2)$ & \hspace{1.3cm}\checkmark (quadratic potential) & \begin{tabular}{@{}l l@{}}
						\hspace{0.7cm}	- \hspace{1.85cm} & \checkmark
					\end{tabular}  &\hspace{1.25cm} - \\
					\hline
				\end{tabular}
			}
			\caption{\centering{Properties of the methods.}}
			\label{tab:geo-properties}
		\end{table}

		\begin{remark}
			As summarized in Table \ref{tab:geo-properties}, all methods M1--M4 achieve second-order accuracy. M1 demonstrates both similarities and differences relative to M2--M4 in long-term geometric conservation. Notably, although M2, M3, and M4 are theoretically equivalent, numerical results reveal practical discrepancies in their behavior, as detailed below:
			\begin{enumerate}[label=(\roman*)]
				\item In Problems 1--4, M2 achieves slightly smaller errors in the $x$-direction and performs more favorably than M3 and M4 in the $v$-direction;
				
				\item M2 also provides better preservation of the mass shell, consistently showing smaller errors than M3 and M4 in Problems 2--4;
				
				\item Under a constant strong magnetic field in Problem 5, M2 demonstrates superior accuracy and long-term near-conservation of energy and mass shell compared to M3 and M4. This will be an important topic for further investigation in our future work.
			\end{enumerate}
			These observations motivate the development of the diverse method formulations presented in this paper.
		\end{remark}

\subsection{Numerical experiments}\label{experiments}
In this subsection, numerical experiments are conducted to verify the error estimates and demonstrate the long-term conservation or near-conservation of energy and mass shell by the proposed methods M1--M4.  
To test the performance of these methods, we evaluate the global errors
$error_{x}:=\frac{\abs{x_{n}-x(\tau_n)}}{\abs{x(\tau_n)}},error_{v}:=\frac{\abs{v_{n}-v(\tau_n)}}{\abs{v(\tau_n)}}$, 
the energy  error
$e_{H}:=\frac{\abs{H(x_{n},\gamma_n)-	H(x_{0},\gamma_0)}}{\abs{H(x_{0},\gamma_0)}}$ and
the mass shell error
$e_{M}:=\frac{\abs{\mathcal{H}({v_{n}},w_{n})-	\mathcal{H}({v_{0}},w_{0})}}{\abs{	\mathcal{H}({v_{0}},w_{0})}}$. The reference solution is generated using MATLAB's ode45 solver.

\begin{figure}[H]
 \centering\tabcolsep=0.5mm
 \begin{tabular}[c]{cc}
 	\includegraphics[width=4cm,height=2.8cm]{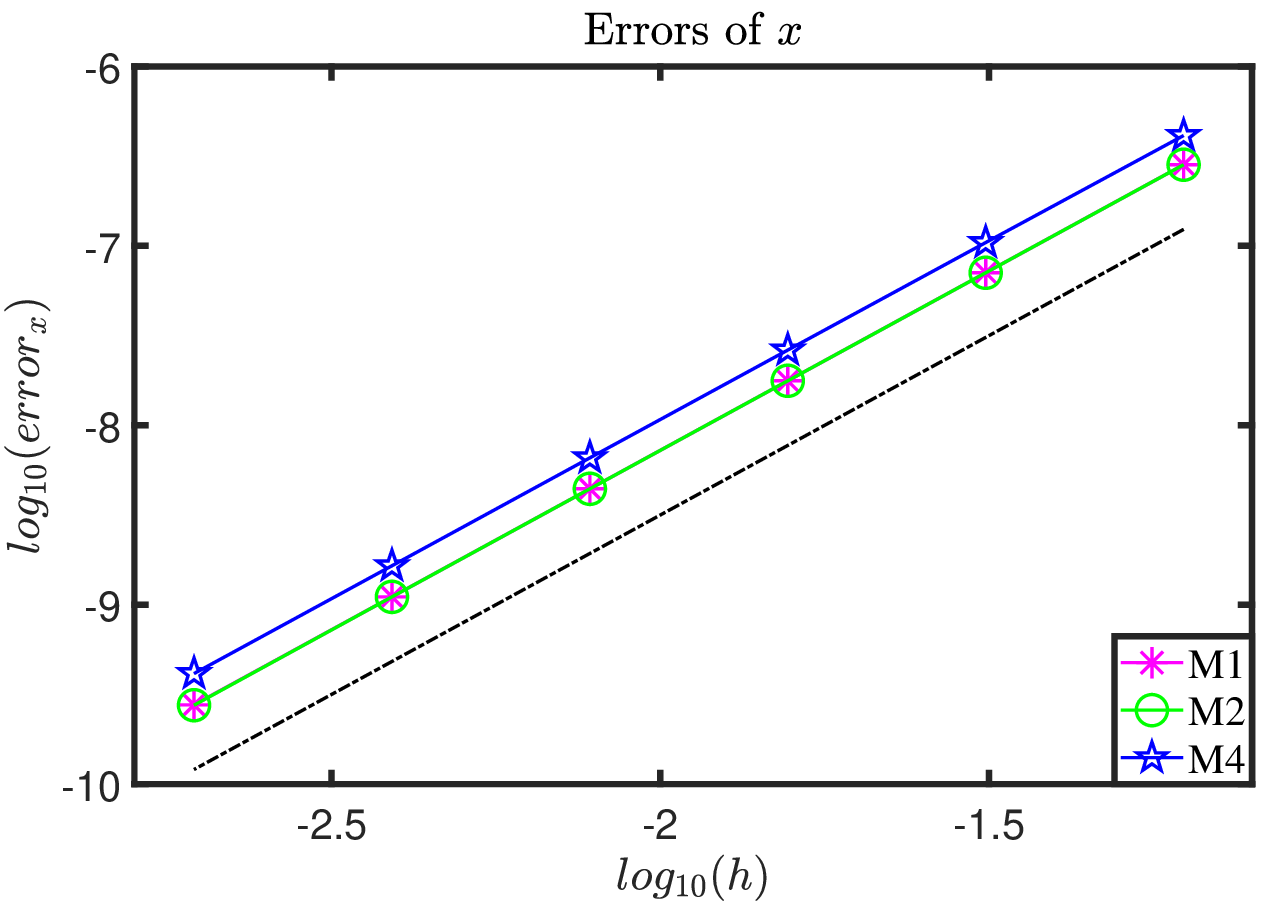}  \qquad\quad
 	\includegraphics[width=4cm,height=2.8cm]{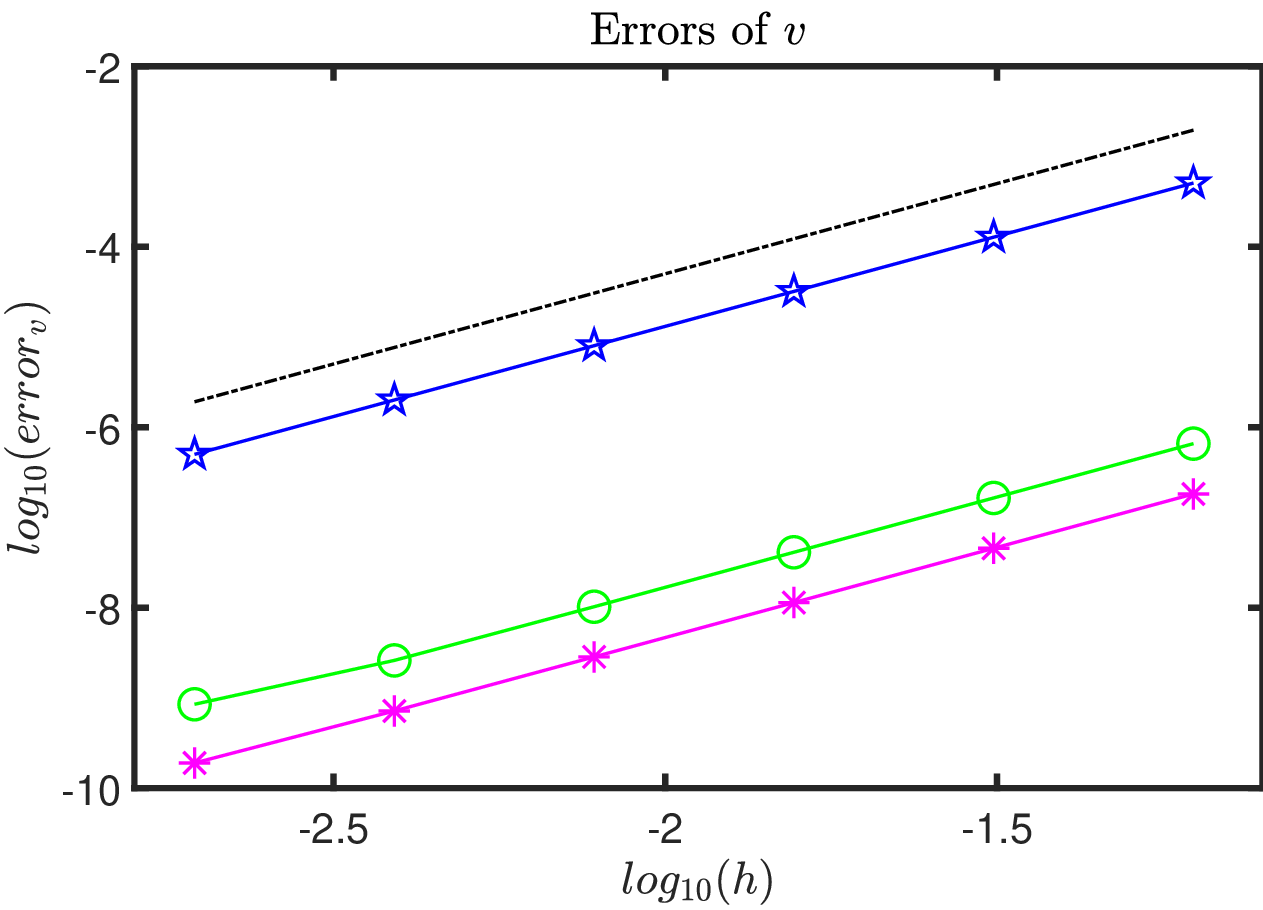}
 \end{tabular}
 \caption{Problem 1 (tokamak field).
 	The global errors $error_{x}$ and $error_{v}$ with $\tau=1$ and $h=1/2^{k}$ for $k=4,\ldots,9$ (the dash-dot line is slope two).}
 \label{fig:problem1Err}
 \end{figure}

\begin{figure}[H]
 \centering\tabcolsep=0.5mm
 \begin{tabular}	[c]{cc}
 	\includegraphics[width=5.5cm,height=2.8cm]{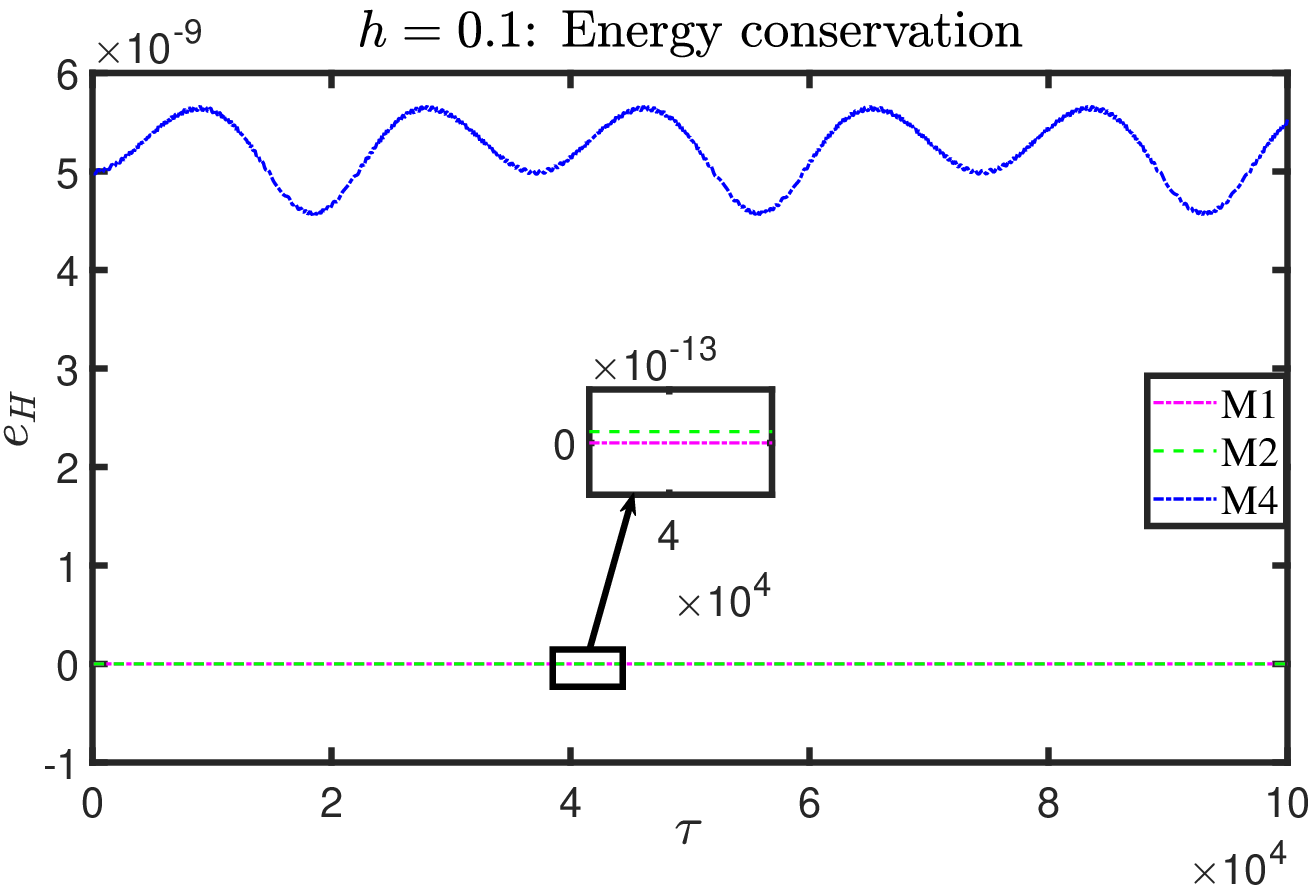}\qquad\quad \includegraphics[width=5.8cm,height=2.8cm]{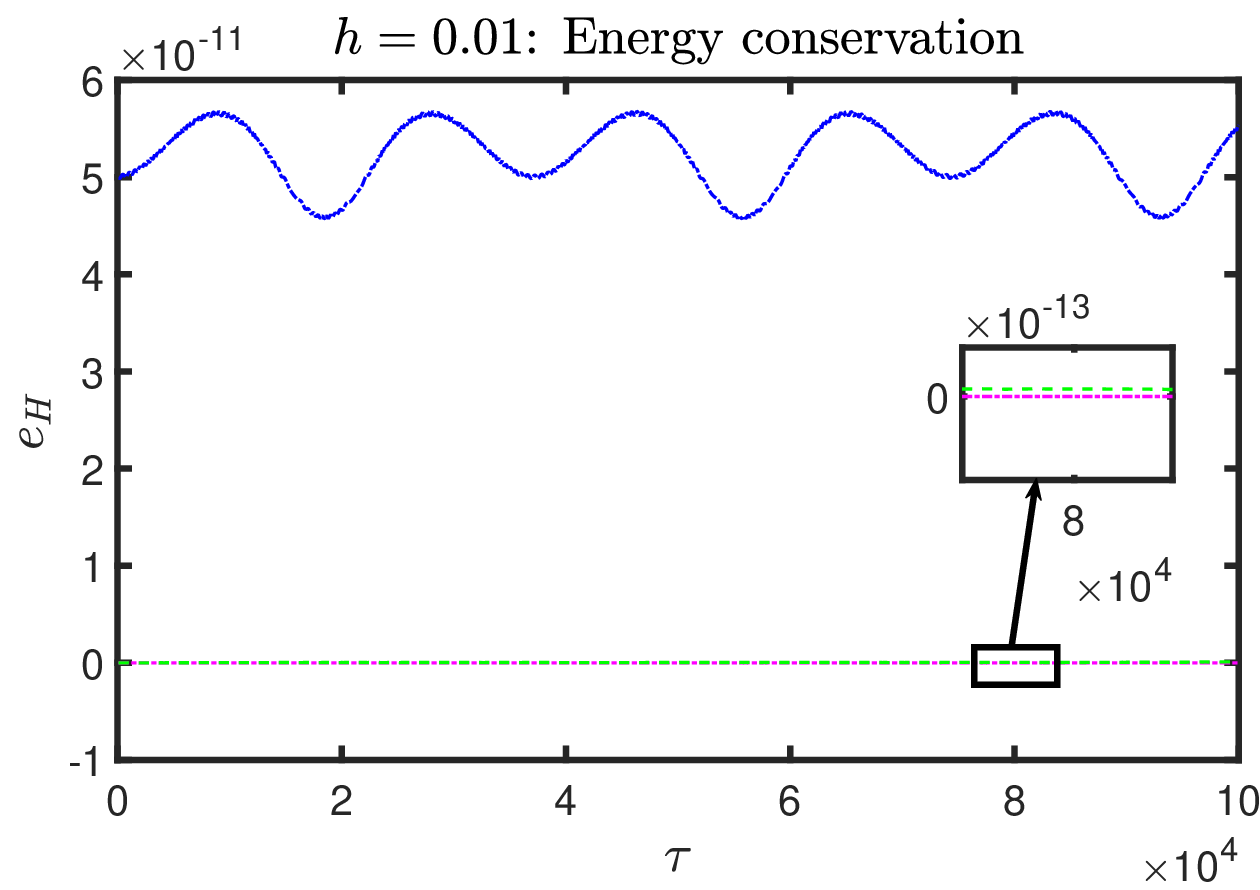}\\
 		\includegraphics[width=5.8cm,height=2.8cm]{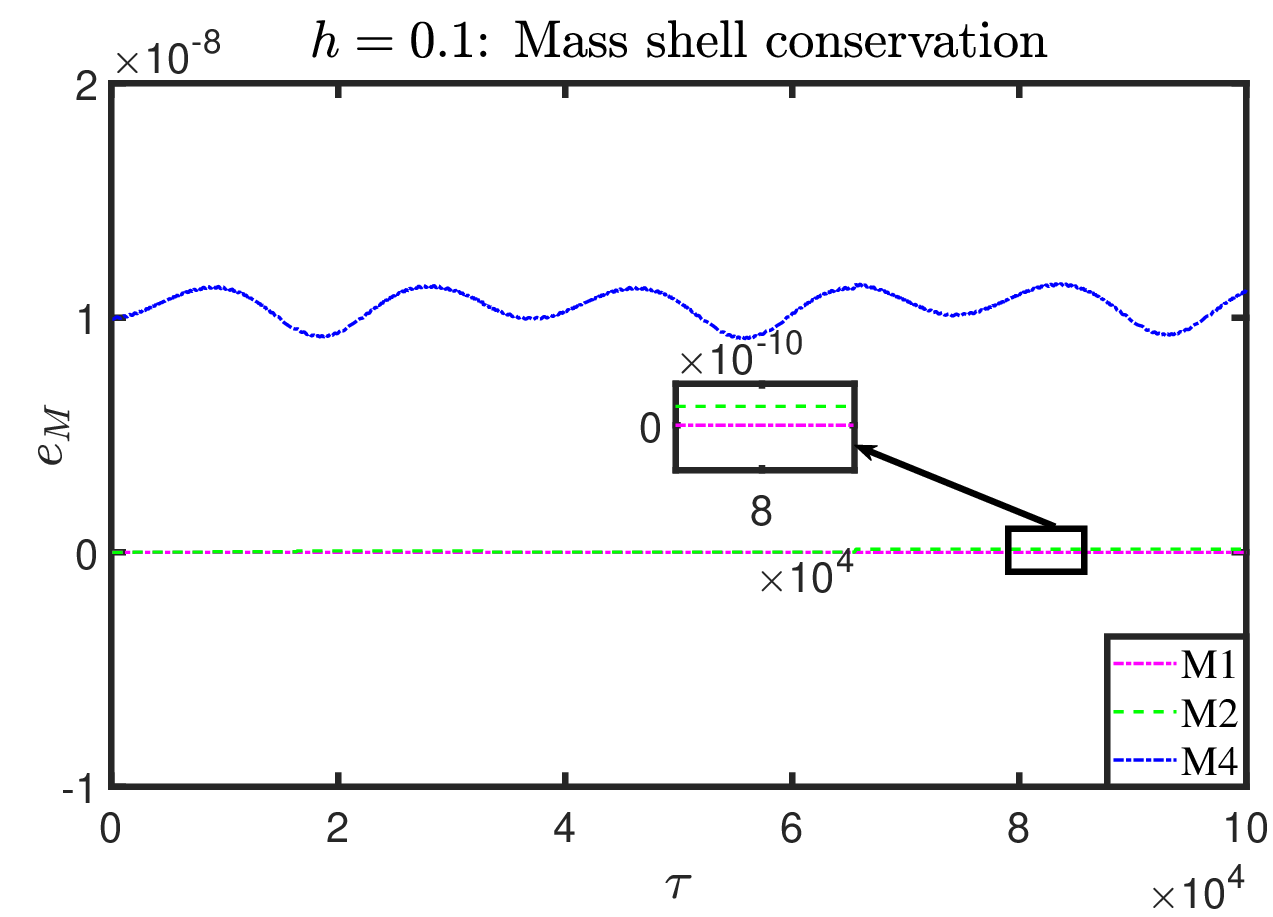}\qquad\quad
 	\includegraphics[width=5.8cm,height=2.8cm]{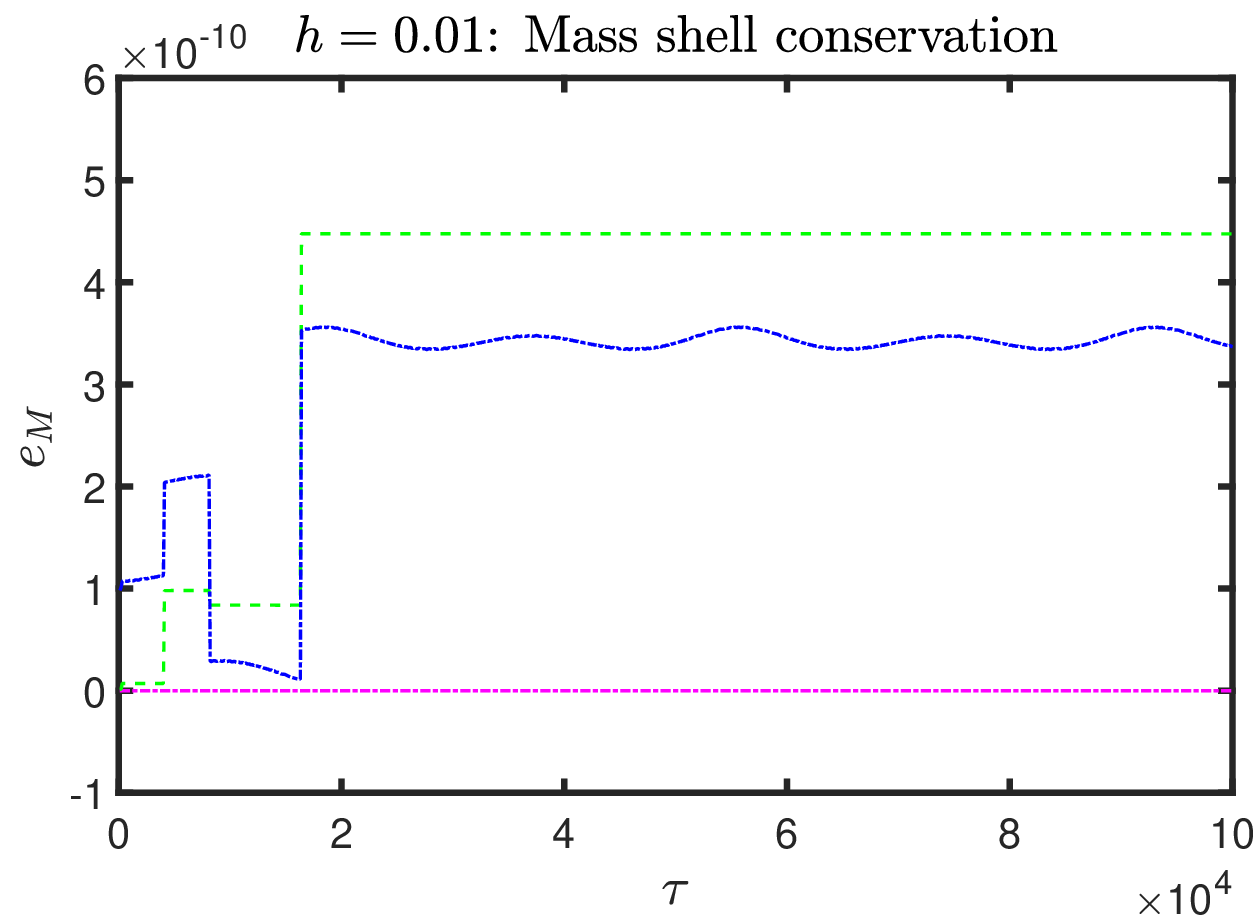}
 \end{tabular}
 \caption{Problem 1 (tokamak field). Evolution of  $e_{H}$ and $e_{M}$ with different step sizes $h$.}
 \label{fig:problem1H}
 \end{figure}

\textbf{Problem 1 (Tokamak field).}
 The first numerical experiment focuses on a charged particle moving in a Tokamak magnetic field without an electric field. 
The magnetic field is given in Cartesian coordinates by
${B}(x)=\big(-\frac{2x_2+x_1x_3}{2R^2},\frac{2x_1-x_2x_3}{2R^2},\frac{R-1}{2R}\big)^{\intercal}$, where $R=\sqrt{x_1^2+x_2^2}$.
The initial values are  defined as  $(x(0); t(0))=(1.05,0,0,0)^{\intercal}$ 
and $(v(0); w(0))=(2.1 \times 10^{-3},4.3 \times 10^{-4},0,w(0))^{\intercal}$  with $w(0)=i\sqrt{1+v(0)^{\intercal}v(0)}$. As M2 and M3 yield essentially identical results in this case, only the numerical results for M2 are presented below. Figure \ref{fig:problem1Banana} displays the projected particle orbits in the $(R, x_3)$ plane.
To illustrate the errors stated in Theorem \ref{error bound2}, we solve this problem on the interval $[0,1]$ with $h = 2^{-k}$ for $k = 4, \dots, 9$, and the global errors are presented in Figure \ref{fig:problem1Err}. We further solve the problem with two step sizes, $h=0.1$ and $h=0.01$, on the interval $[0,100000]$. Figure~\ref{fig:problem1H} shows the conservation errors of the energy $e_{H}$ and the mass shell $e_{M}$. 
 \begin{figure}[H]
 	\centering\tabcolsep=0.5mm
 	\begin{tabular}[c]{cc}
 		\includegraphics[width=3cm,height=2.8cm]{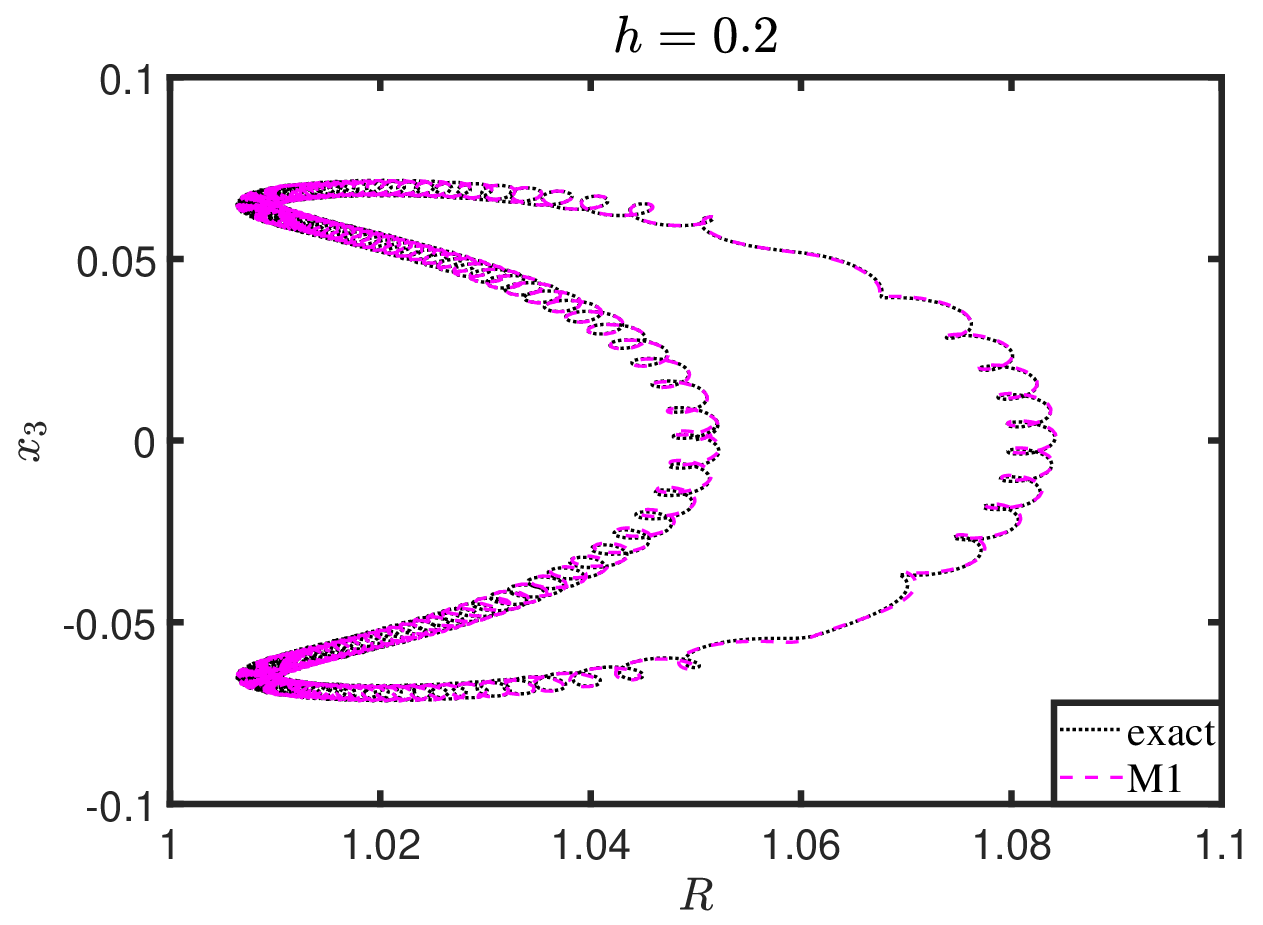}  \qquad
 		\includegraphics[width=3cm,height=2.8cm]{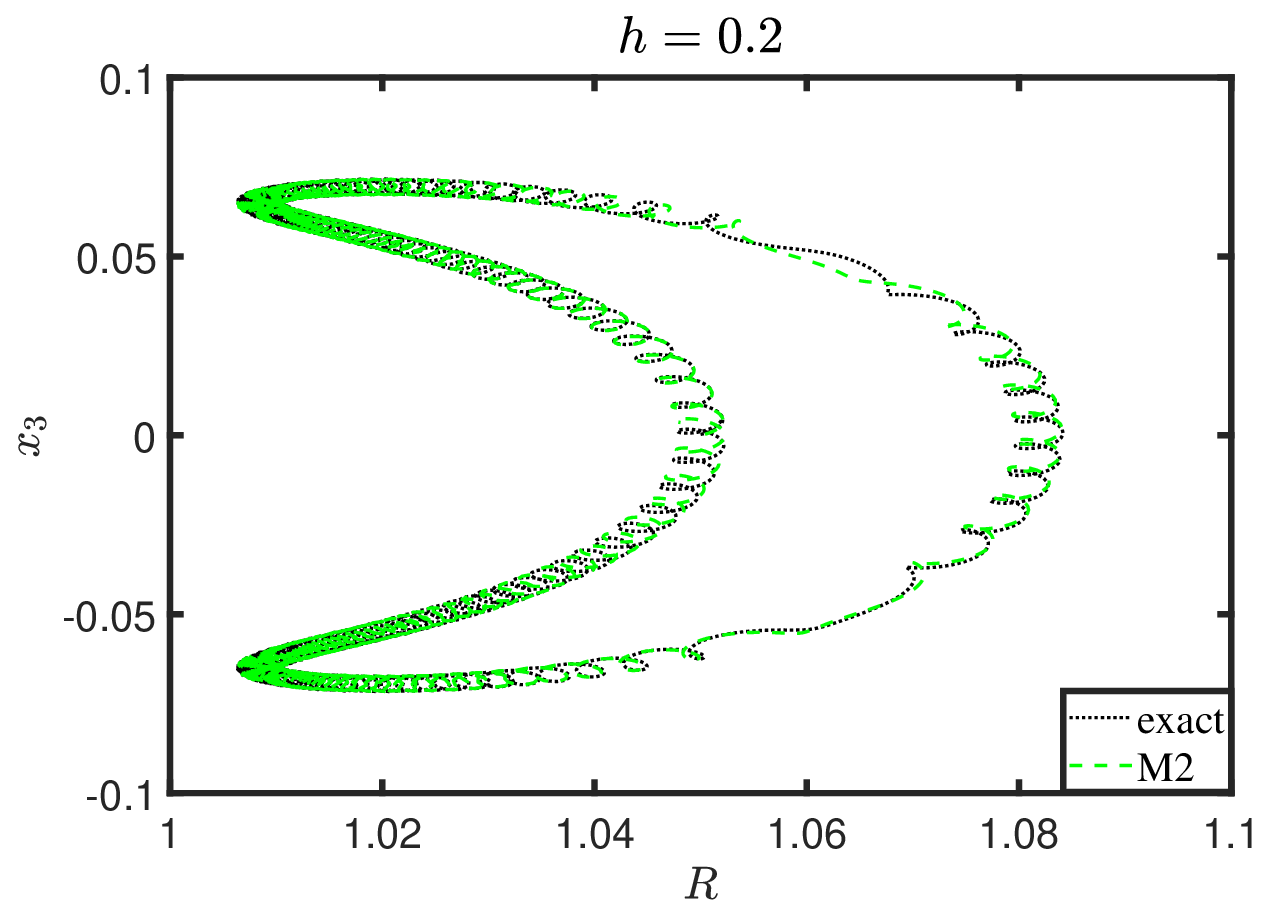}\qquad
 		\includegraphics[width=3cm,height=2.8cm]{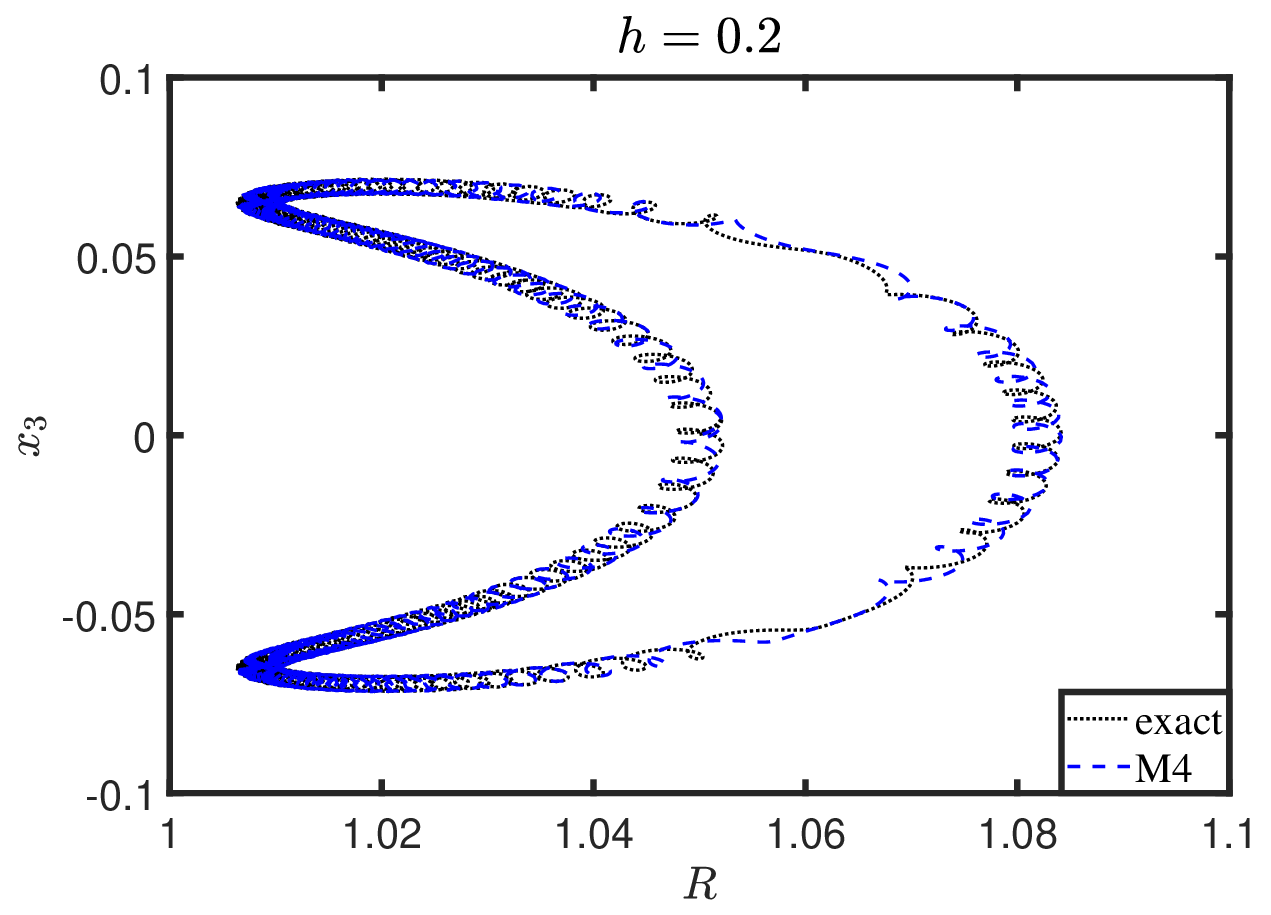}
 	\end{tabular}
 	\caption{Problem 1 (tokamak field). Particle orbit in the poloidal plane of a tokamak obtained by the proposed methods.}
 \label{fig:problem1Banana}
 \end{figure}

\begin{figure}[H]
 \centering\tabcolsep=0.5mm
 \begin{tabular}[c]{cc}
 \includegraphics[width=4cm,height=2.8cm]{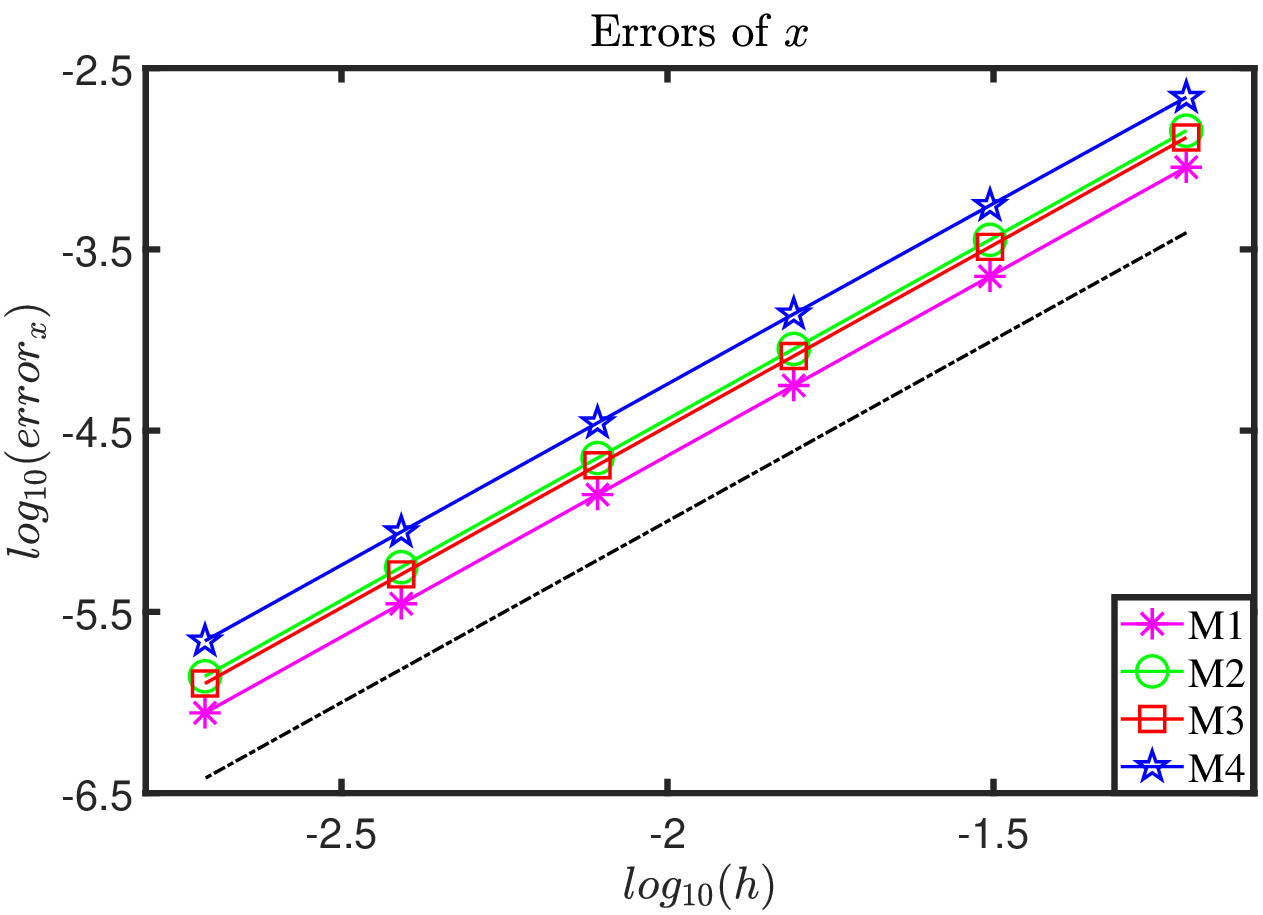}  \qquad\quad
 \includegraphics[width=4cm,height=2.8cm]{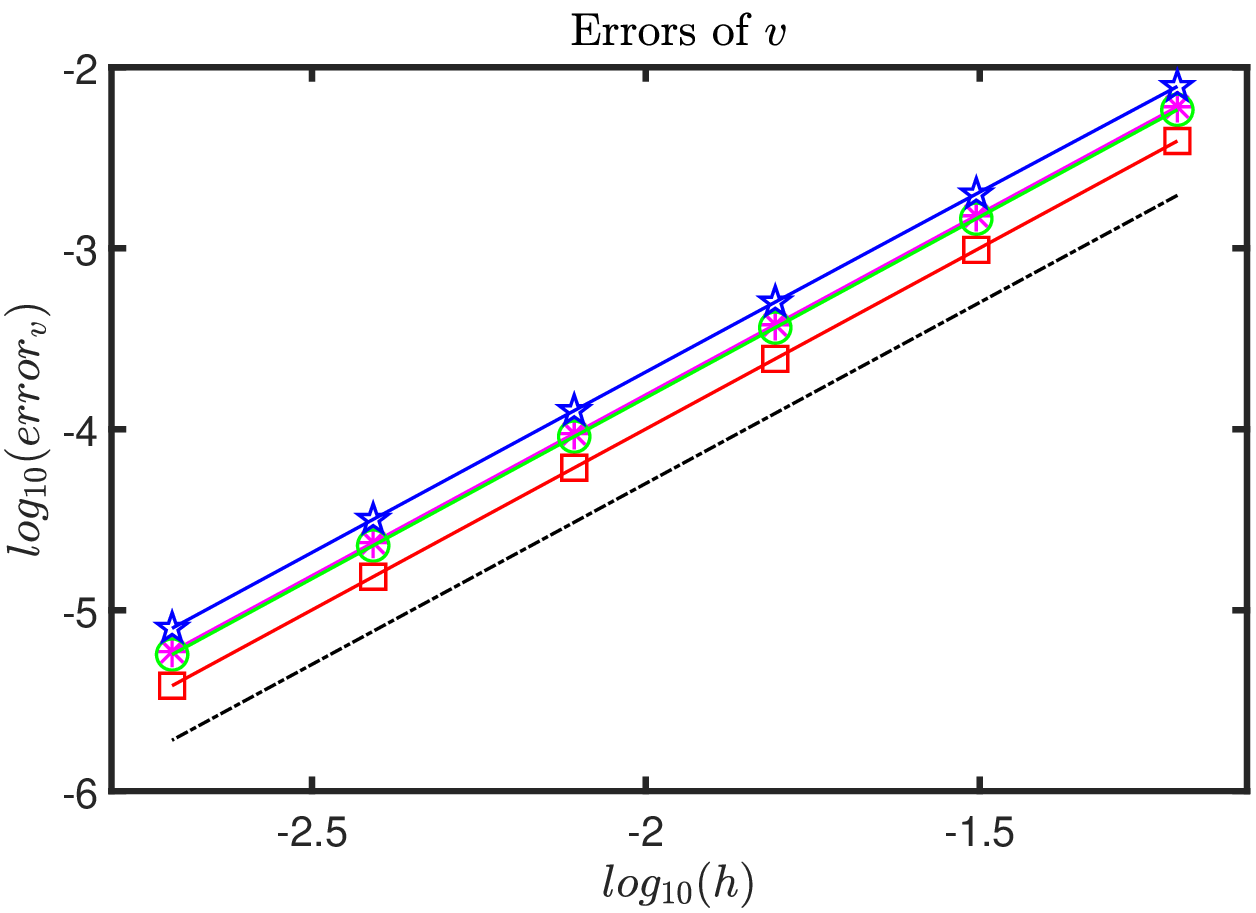}
 \end{tabular}
 \caption{Problem 2 (quadratic electric potential).
 	The global errors $error_{x}$ and $error_{v}$ with $\tau=1$ and $h=1/2^{k}$ for $k=4,\ldots,9$ (the dash-dot line is slope two).}
 \label{fig:problem2Err}
 \end{figure}

\begin{figure}[H]
 \centering\tabcolsep=0.5mm
 \begin{tabular}	[c]{cc}
 \includegraphics[width=5.8cm,height=2.6cm]{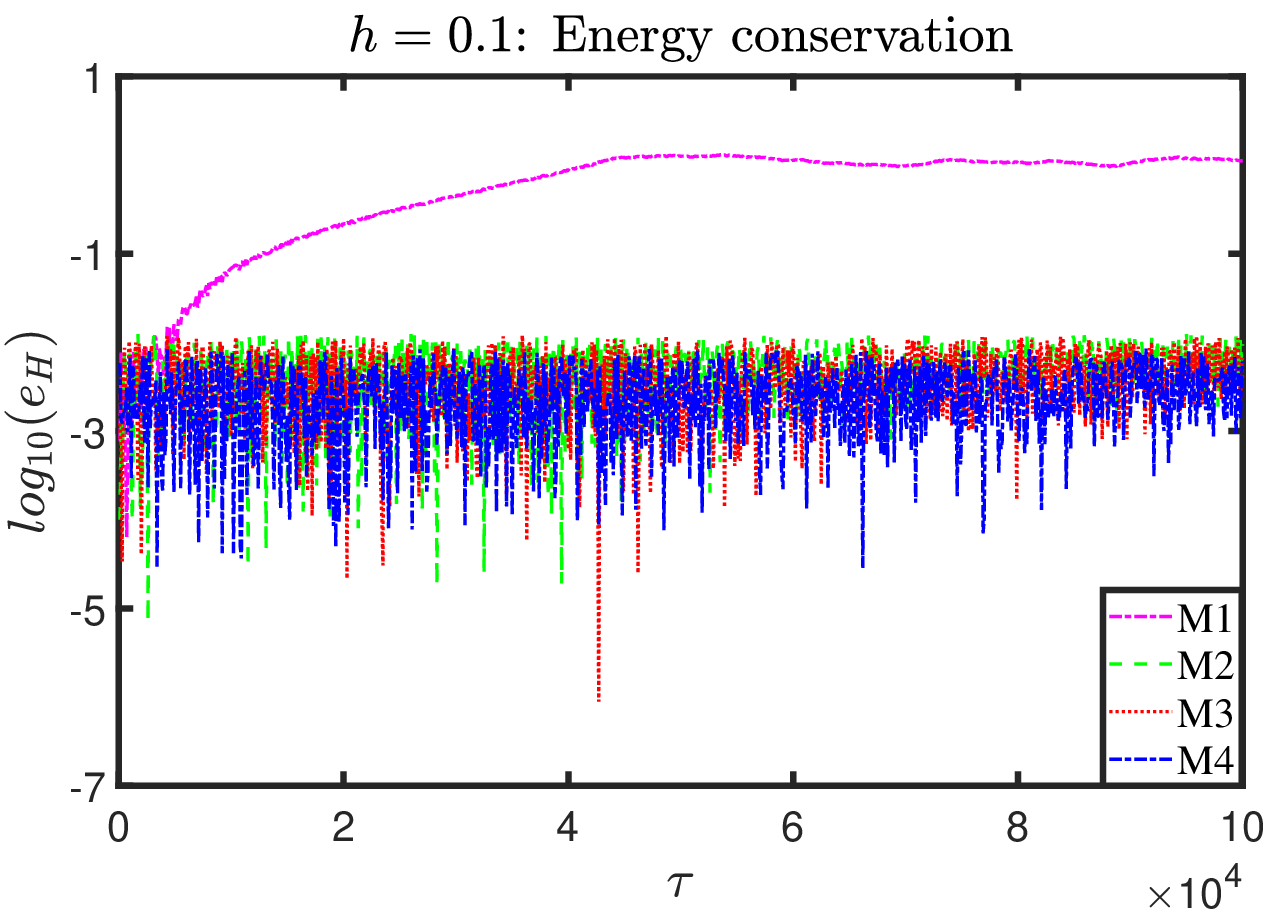}\qquad\quad
  \includegraphics[width=5.8cm,height=2.6cm]{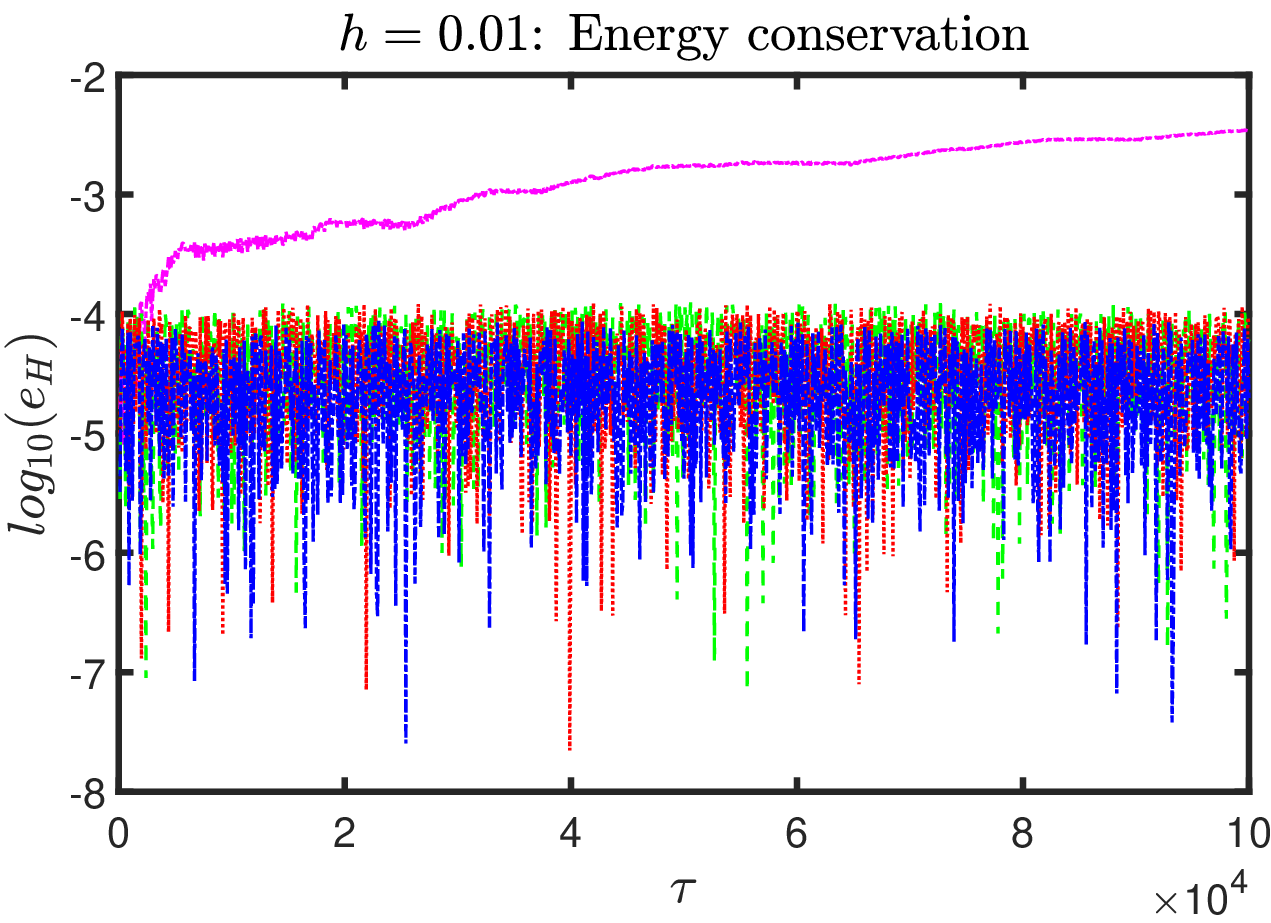}\\
  \includegraphics[width=5.8cm,height=2.6cm]{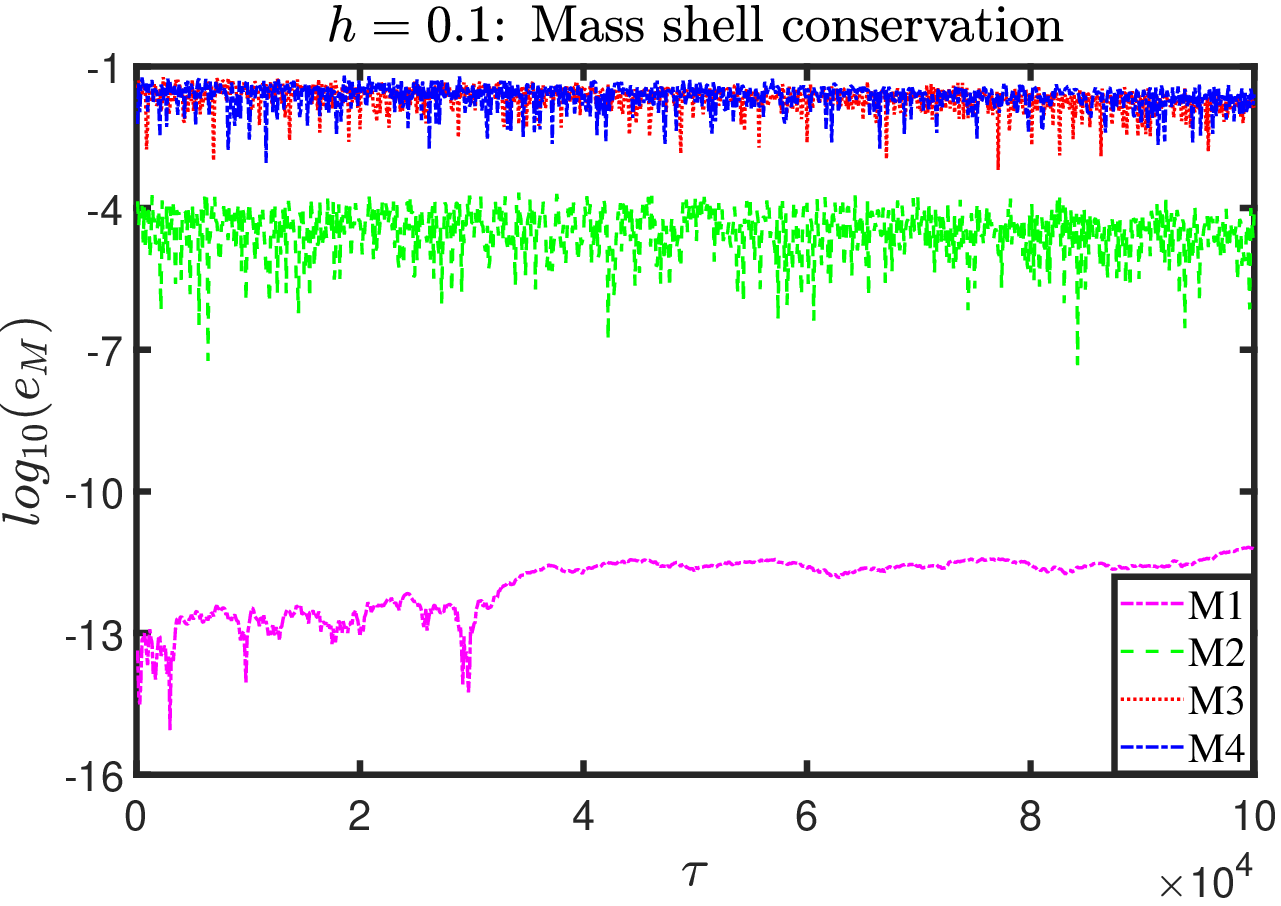}\qquad\quad
  \includegraphics[width=5.8cm,height=2.6cm]{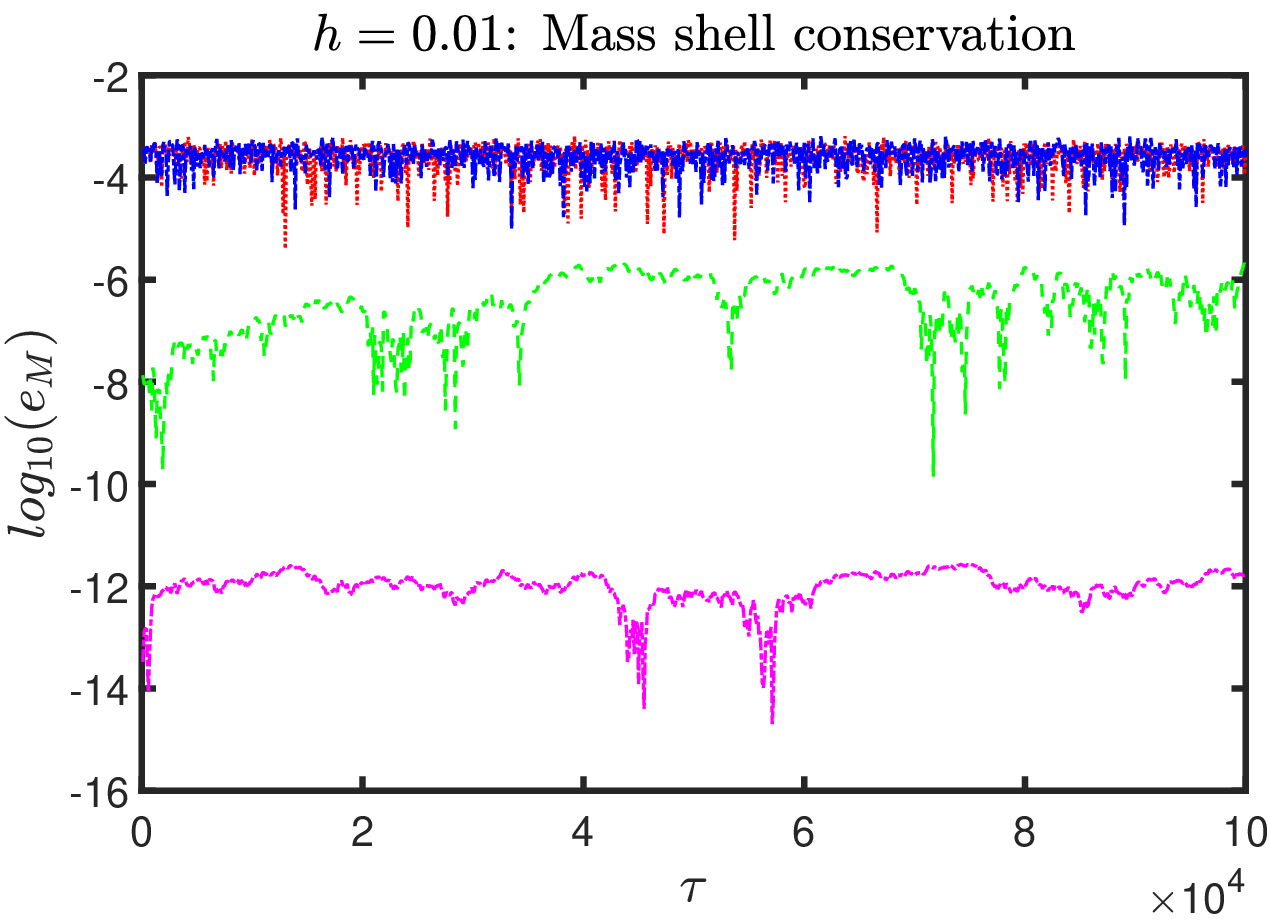}
 \end{tabular}
 \caption{Problem 2 (quadratic electric potential). Evolution of  $e_{H}$ and $e_{M}$ with different step sizes $h$.}
 \label{fig:problem2H}
 \end{figure}

\textbf{Problem 2 (Quadratic electric potential).}\label{ex1}
In the second numerical experiment, we consider
 the motion of a charged particle in the quadratic electric potential
$U(x)=x_1^2+2x_2^2+3x_3^2-x_1$, where the electric field is given by 
$E(x)=-\nabla U(x)$ and  the magnetic field is ${B}(x)=(\cos(x_2)-x_1,1+\sin(x_3),\cos(x_1)+x_3)^{\intercal}/2$. We choose the initial values as $(x(0); t(0))=(1/3,1/4,1/2,0)^{\intercal}$ 
and $(v(0); w(0))=(2/5,2/3,1/6,w(0))^{\intercal}$, where $w(0)=i\sqrt{1+v(0)^{\intercal}v(0)}$.
The global errors and the numerical conservation of  $H$ and $\mathcal{H}$ is presented in Figures \ref{fig:problem2Err} and \ref{fig:problem2H}, respectively.

\begin{figure}[H]
 \centering\tabcolsep=0.5mm
 \begin{tabular}
 [c]{cc}
 \includegraphics[width=4cm,height=2.8cm]{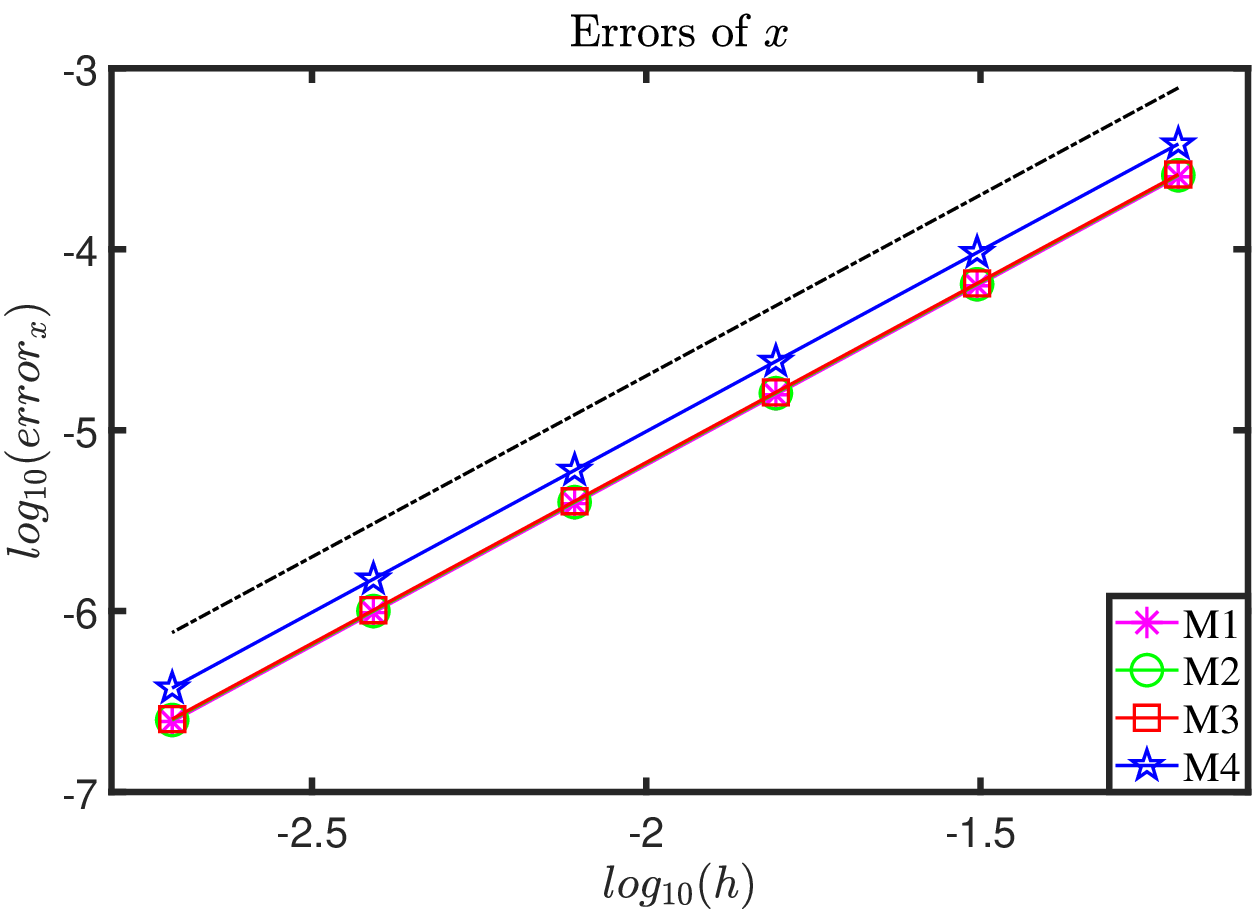}  \qquad\quad
 \includegraphics[width=4cm,height=2.8cm]{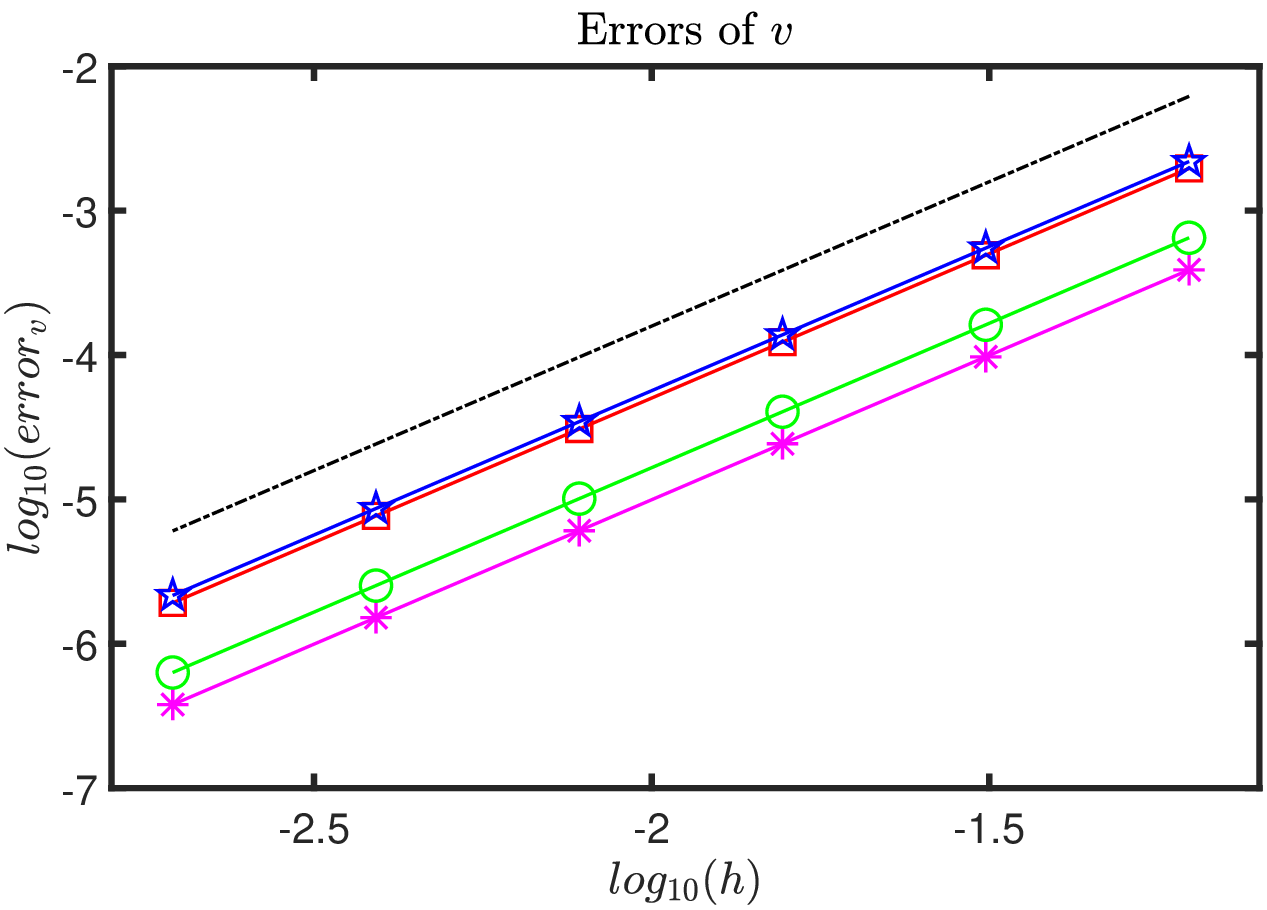}
 \end{tabular}
 \caption{Problem 3 (non-quadratic electric potential).
 		The global errors $error_{x}$ and $error_{v}$ with $\tau=1$ and $h=1/2^{k}$ for $k=4,\ldots,9$ (the dash-dot line is slope two).}
 \label{fig:problem3Err}
 \end{figure}

\begin{figure}[H]
	\centering\tabcolsep=0.5mm
	\begin{tabular}[c]{cc}
		\includegraphics[width=5.8cm,height=2.8cm]{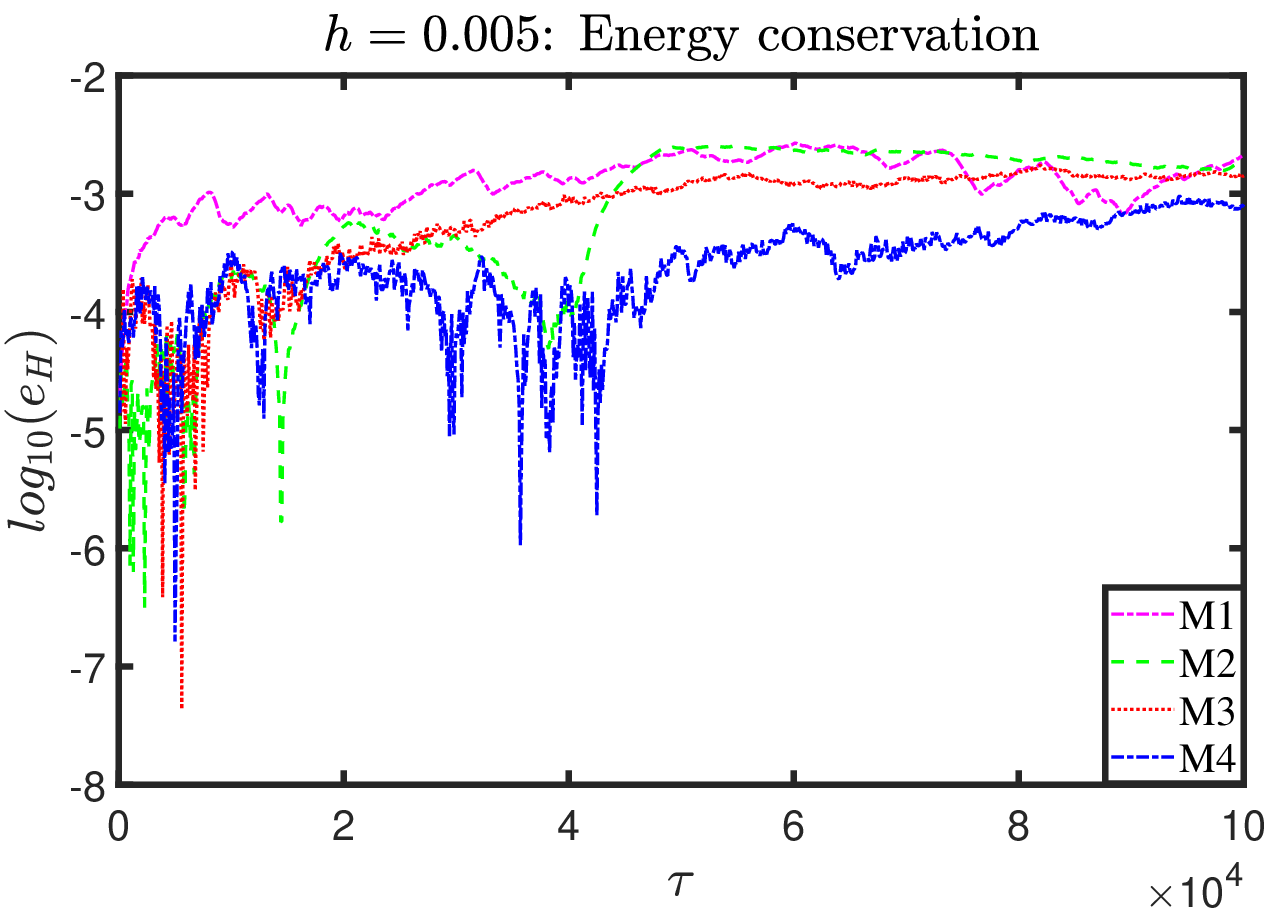}\qquad\quad \includegraphics[width=5.8cm,height=2.8cm]{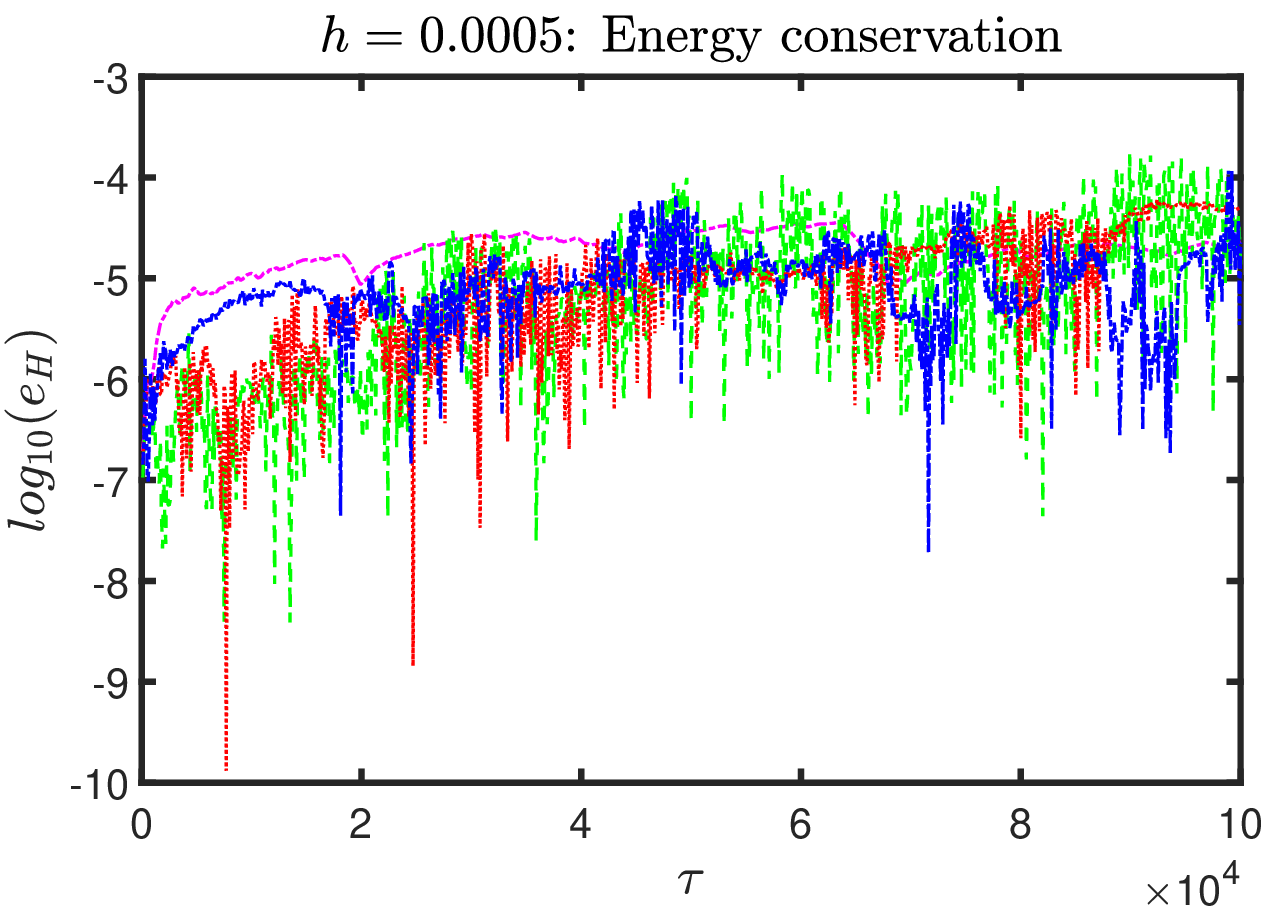}
	\end{tabular}
	\caption{Problem  3 (non-quadratic electric potential). Evolution of  $e_{H}$  with different step sizes $h$.}
	\label{fig:problem3H}
\end{figure}

\textbf{Problem 3 (Non-quadratic electric potential).}
The third numerical experiment  involves  the motion of a charged particle in a non-quadratic electric potential $U(x)=(x_1^3-x_2^3+x_1^4/5+x_2^4+x_3^4)/10$, with $E(x)=-\nabla U(x)$, and  the magnetic field  $B(x)=(\cos(x_2)-x_1,1+\sin(x_3),\cos(x_1)+x_3)^{\intercal}$. The initial values are given by  $(x(0); t(0))=(0,1,0.1,0)^{\intercal}$ 
and $(v(0); w(0))=(0.09,0.55,0.2,w(0))^{\intercal}$, where   $w(0)=i\sqrt{1+v(0)^{\intercal}v(0)}$.
Figures \ref{fig:problem3Err} and \ref{fig:problem3H} show the global errors and the numerical conservation of 
$H$ and $\mathcal{H}$, respectively.

 \begin{figure}[H]
 	\centering\tabcolsep=0.5mm
 	\begin{tabular}[c]{cc}
 			 \includegraphics[width=5.8cm,height=2.8cm]{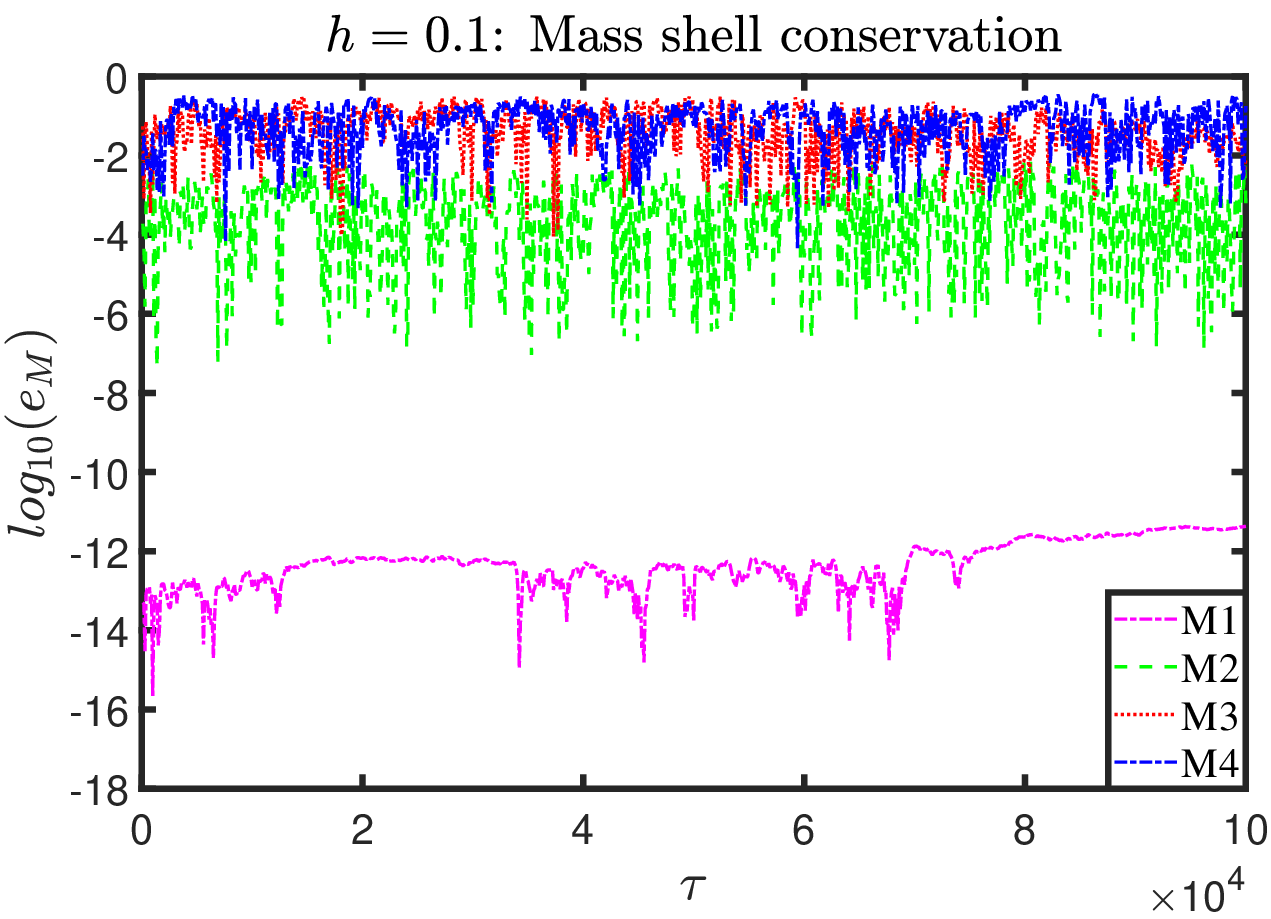}\qquad\quad
 			\includegraphics[width=5.8cm,height=2.8cm]{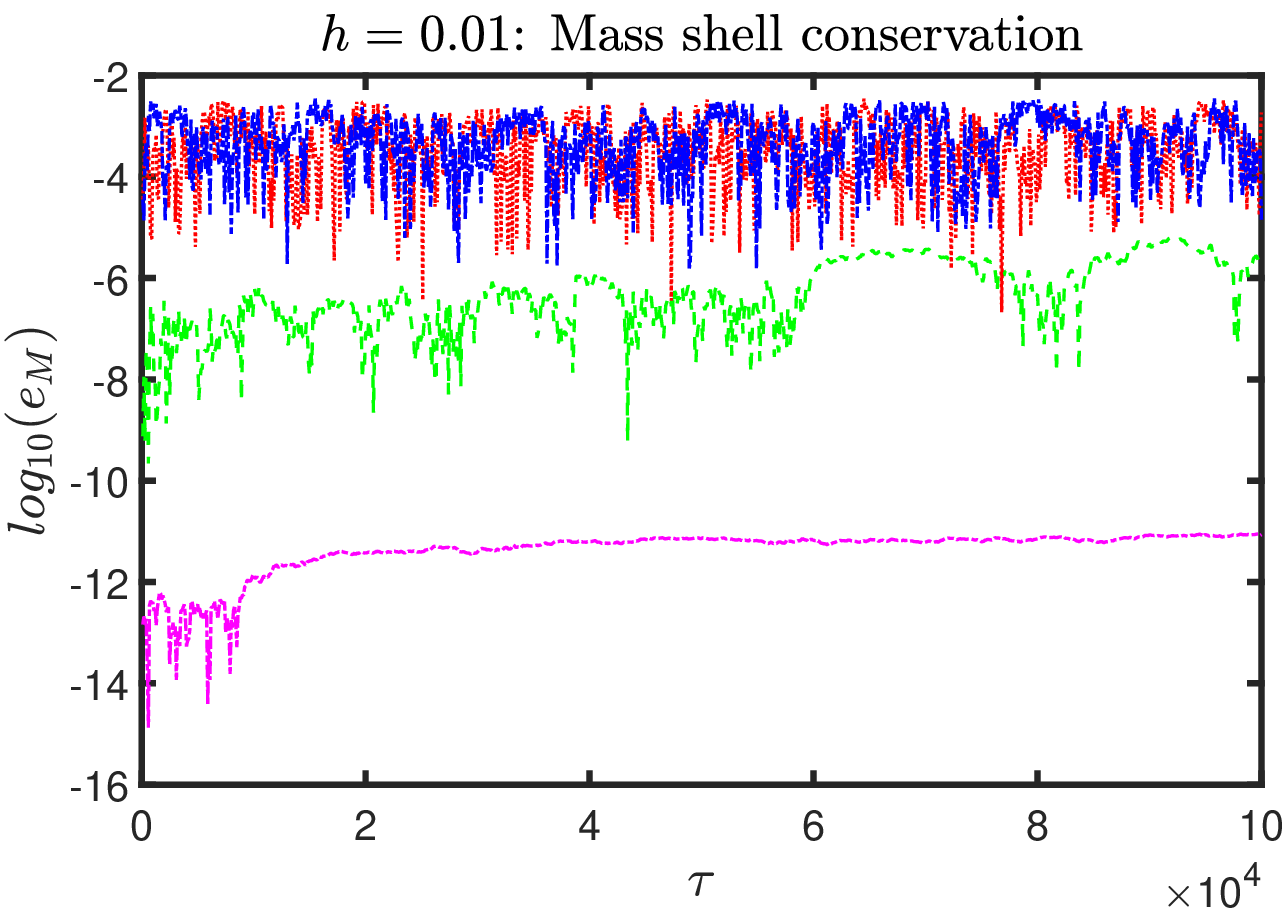}
 		\end{tabular}
 	\caption{Problem  3 (non-quadratic electric potential). Evolution of   $e_{M}$ with different step sizes $h$.}
 	\label{fig:problem3M}
 \end{figure}

\begin{figure}[h!]
 	\centering\tabcolsep=0.5mm
 	\begin{tabular}[c]{cc}
 \includegraphics[width=4cm,height=3cm]{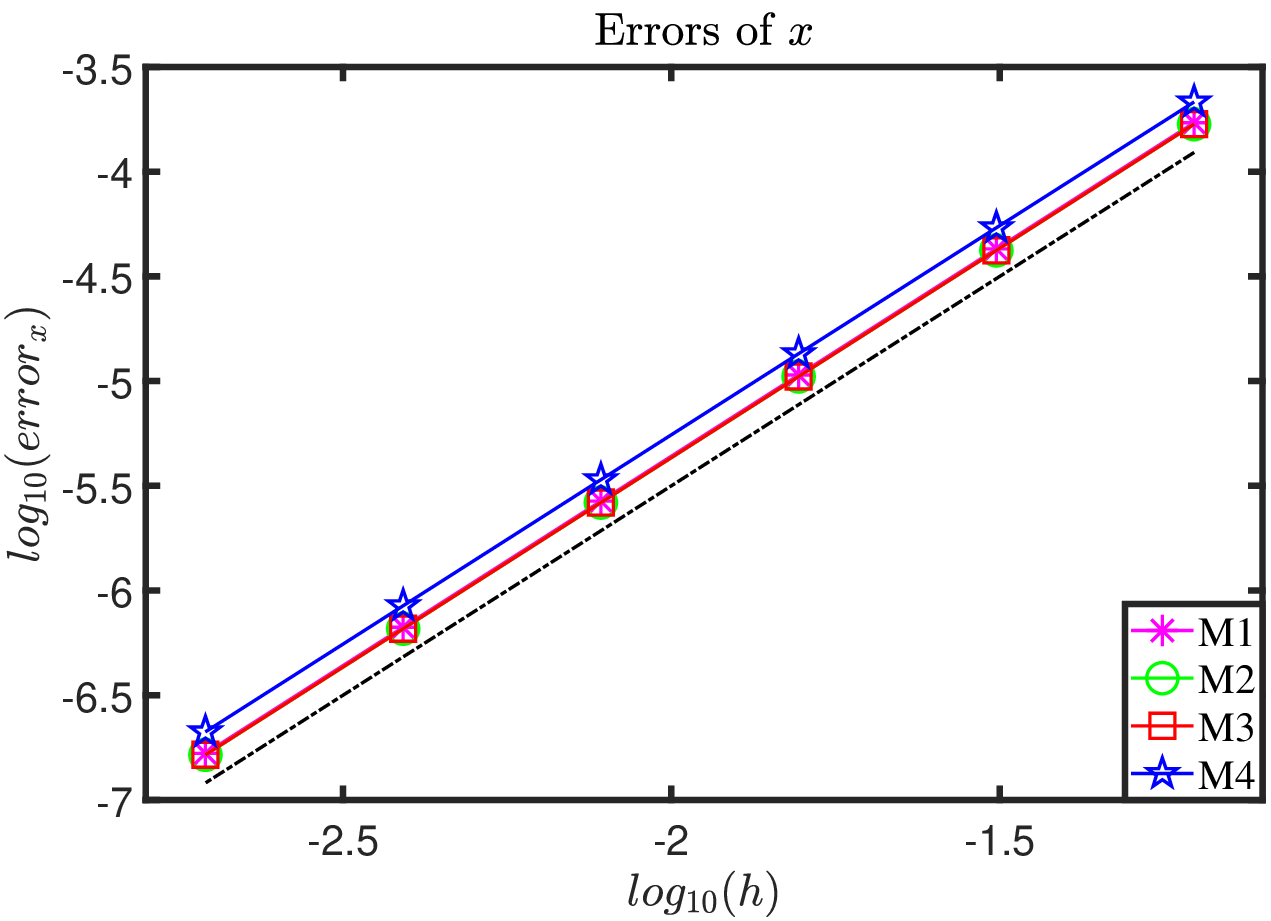}  \qquad\quad
 \includegraphics[width=4cm,height=3cm]{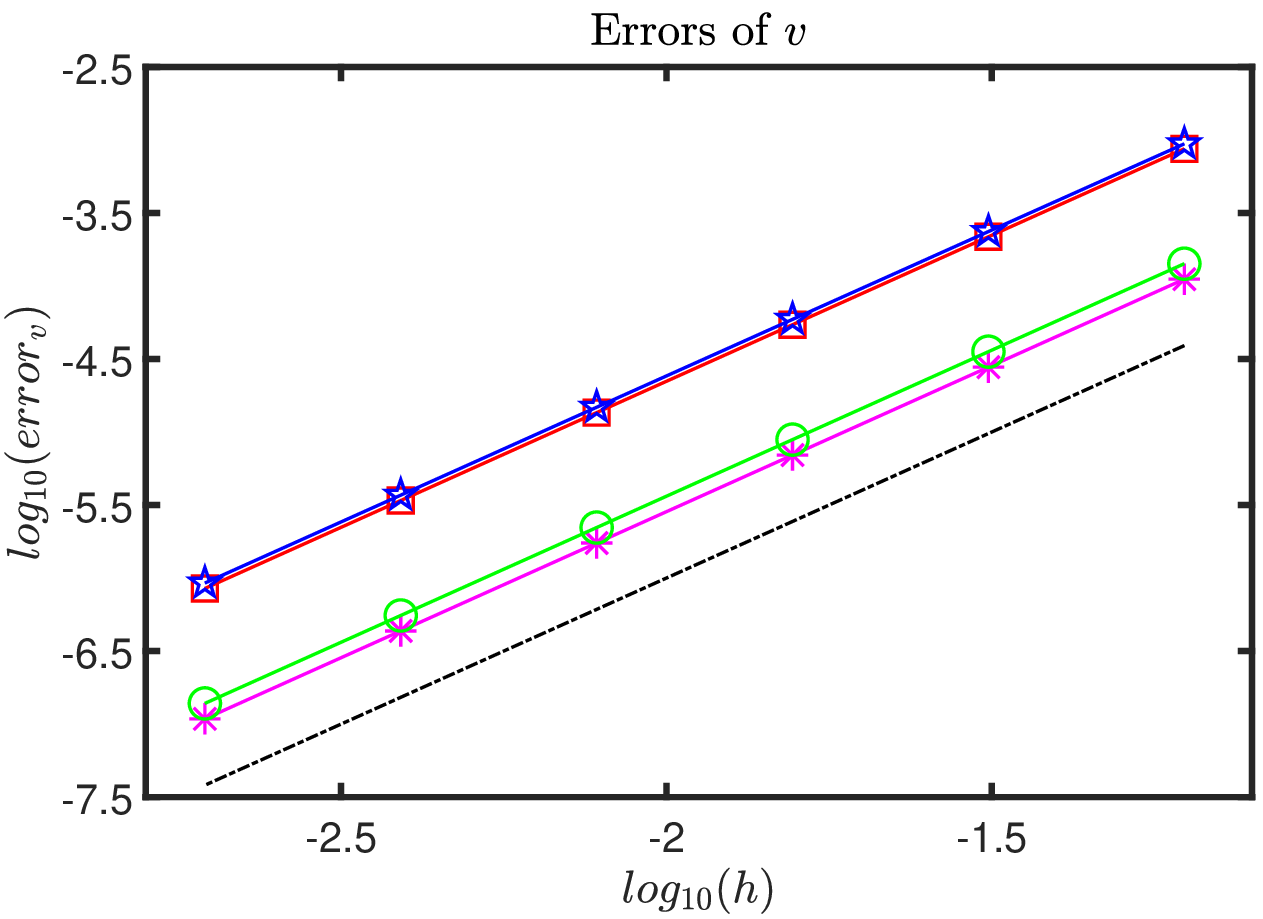}
 			\end{tabular}
 	\caption{Problem 4 (constant magnetic field).
 				The global errors $error_{x}$ and $error_{v}$ with $\tau=1$ and $h=1/2^{k}$ for $k=4,\ldots,9$ (the dash-dot line is slope two).}
 	\label{fig:problem4Err}
 \end{figure}

\begin{figure}[H]
 	\centering\tabcolsep=0.5mm
 	\begin{tabular}[c]{cc}
  \includegraphics[width=5.8cm,height=2.8cm]{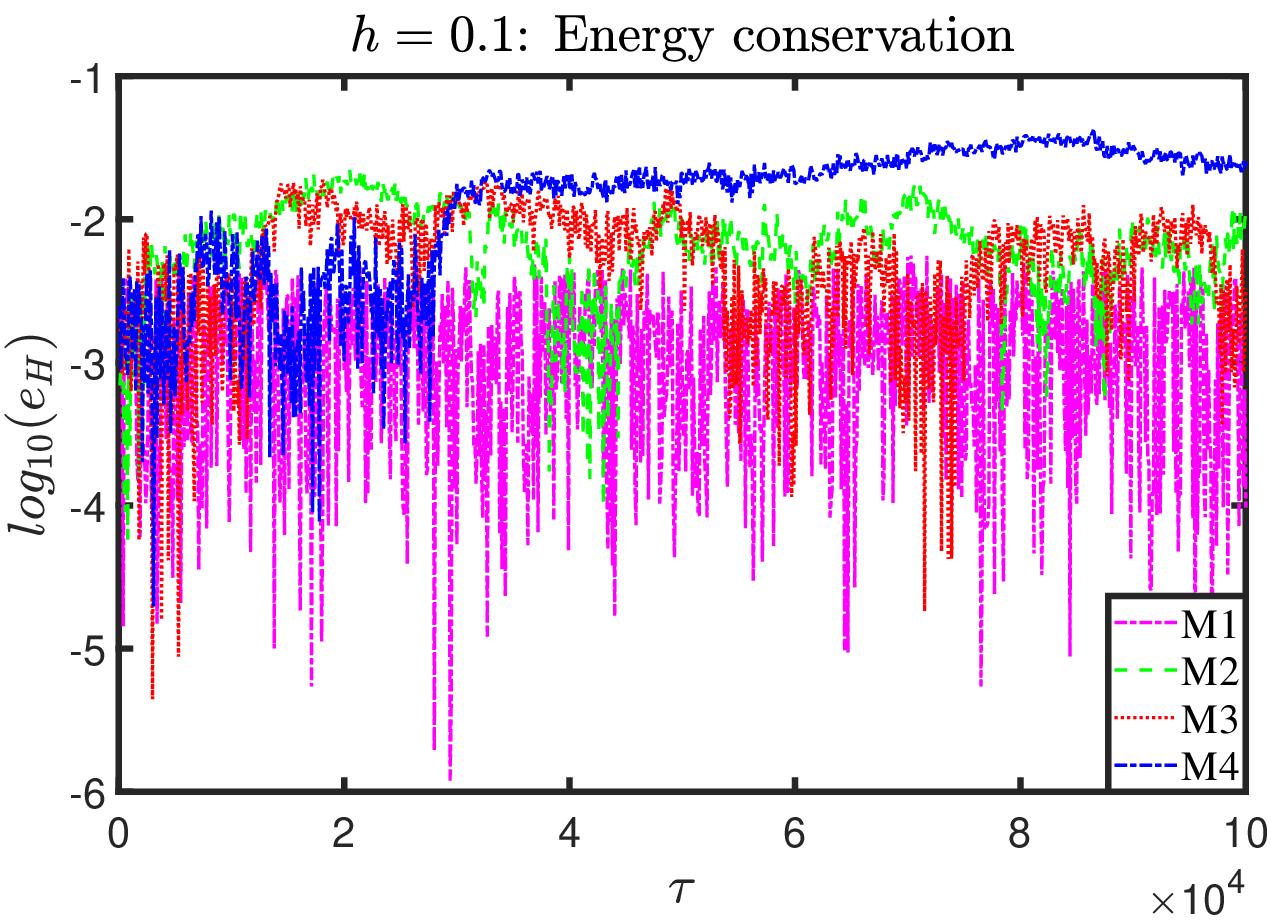}\qquad\quad \includegraphics[width=5.8cm,height=2.8cm]{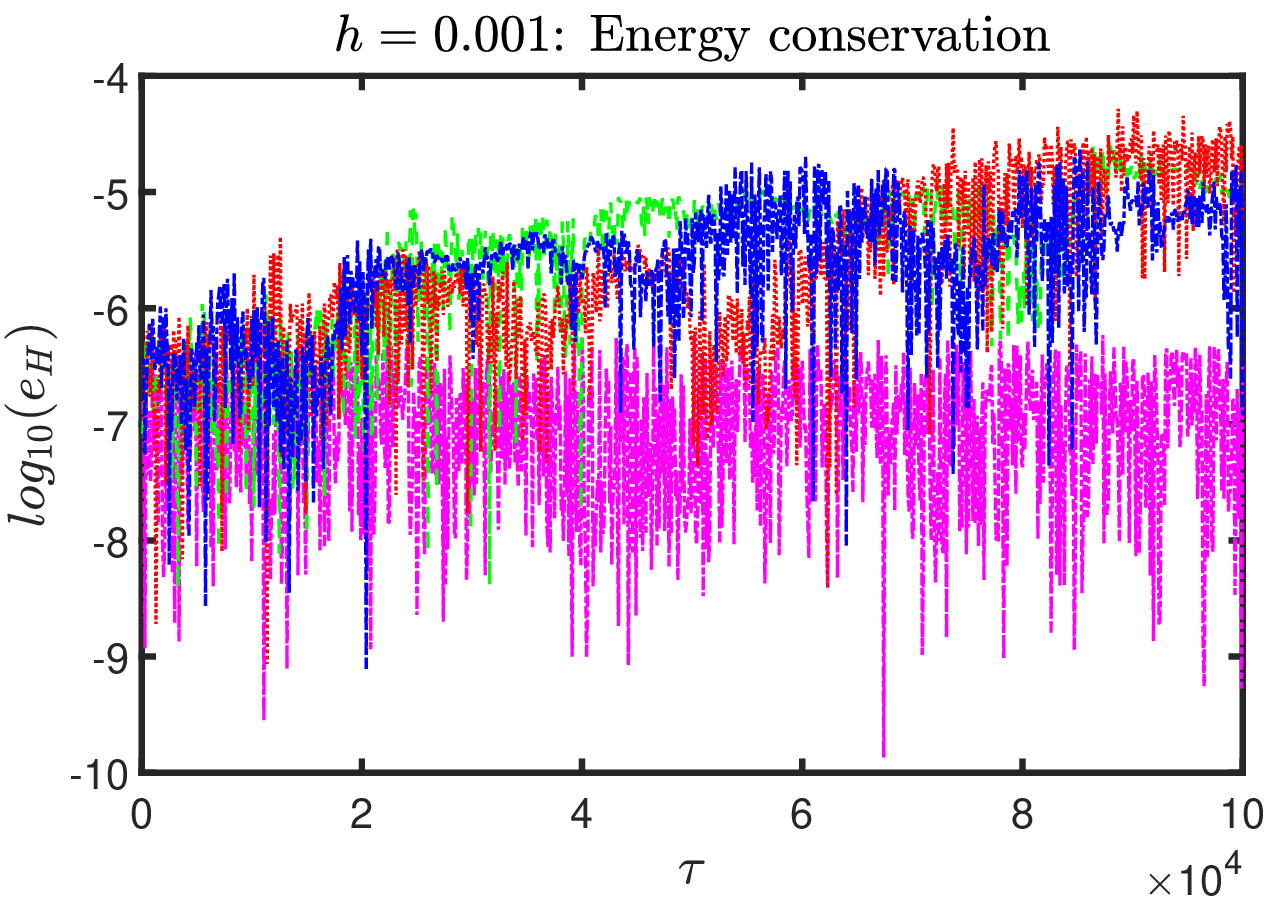}\\
  	\includegraphics[width=5.8cm,height=2.8cm]{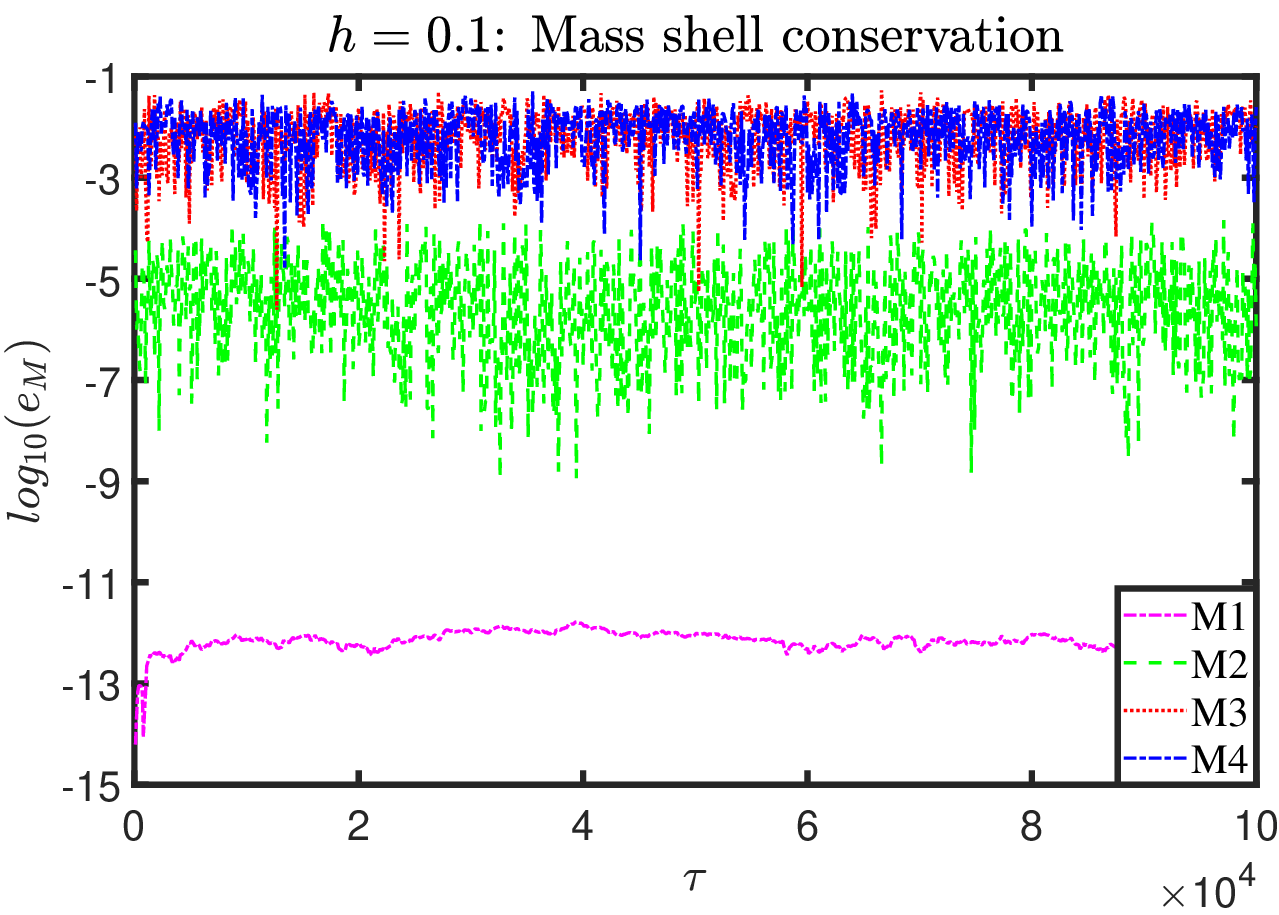}\qquad\quad
  \includegraphics[width=5.8cm,height=2.8cm]{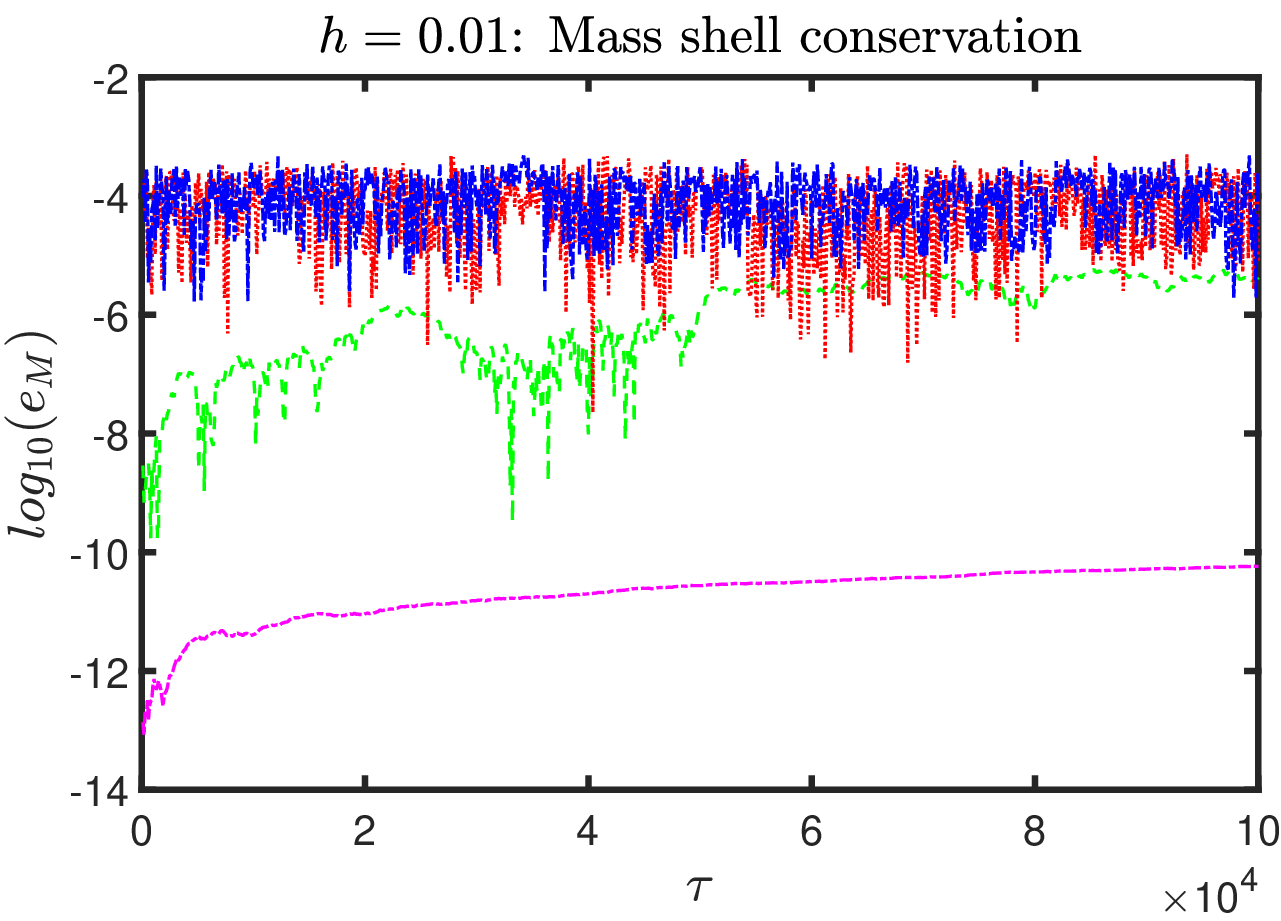}
 		\end{tabular}
 	\caption{Problem  4 (constant magnetic field). Evolution of  $e_{H}$ and $e_{M}$ with different step sizes $h$.}
 	\label{fig:problem4H}
 \end{figure}

\textbf{Problem 4 (Constant magnetic field).}
The fourth numerical experiment investigates a constant magnetic field $B=(0,0,1)^{\intercal}$ and the electric potential $U(x)=(x_1^3-x_2^3+x_1^4/5+x_2^4+x_3^4)/10$,
with $E(x)=-\nabla U(x)$. The initial values are  chosen as  $(x(0); t(0))=(0,1,0.1,0)^{\intercal}$ and  $(v(0); w(0))=(0.09,0.55,0.3,w(0))^{\intercal}$, where $w(0)=i\sqrt{1+v(0)^{\intercal}v(0)}$.
The global errors are shown in Figure \ref{fig:problem4Err}, and the numerical conservation of 
 $H$ and $\mathcal{H}$ is presented in Figure \ref{fig:problem4H}.

\begin{figure}[H]
 	\centering\tabcolsep=0.5mm
 	\begin{tabular}[c]{cc}
 \includegraphics[width=4cm,height=2.8cm]{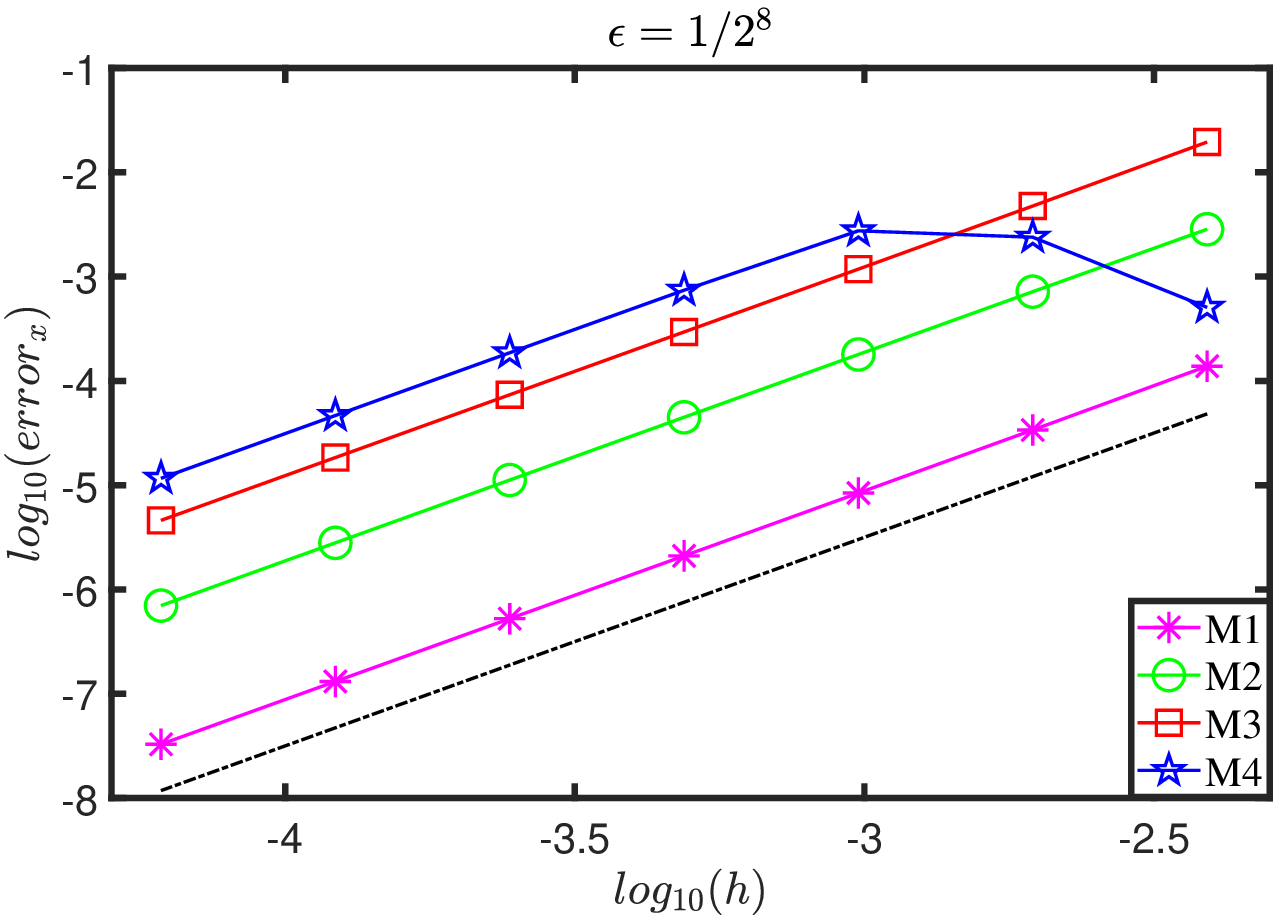}  \qquad\quad
 \includegraphics[width=4cm,height=2.8cm]{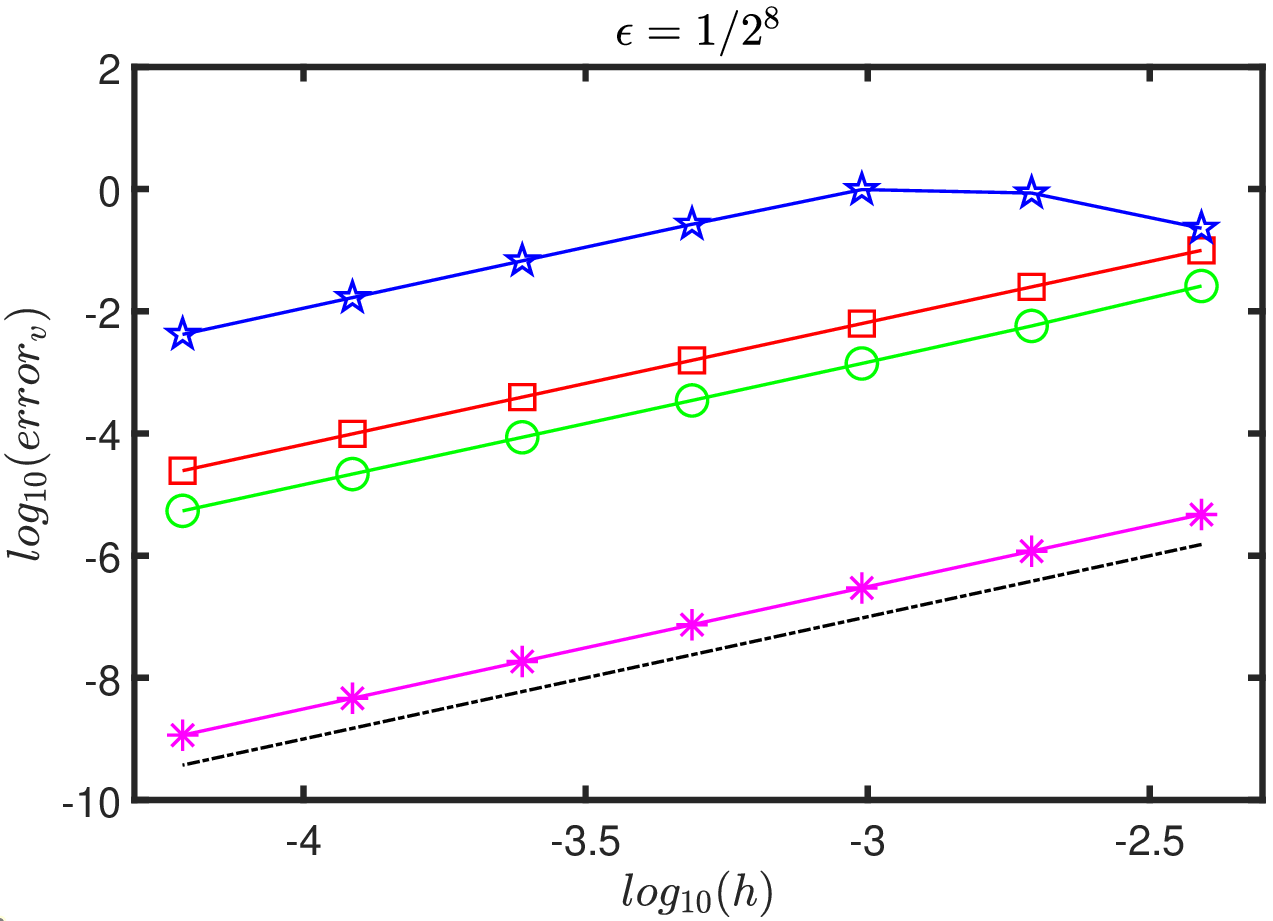}
 			\end{tabular}
 	\caption{Problem 5 (strong constant magnetic field).
 				The global errors $error_{x}$ and $error_{v}$ with $\tau=1$ and $h=1/2^{k}$ for $k=8,\ldots,14$ (the dash-dot line is slope two).}
 	\label{fig:problem5Err}
 \end{figure}

\begin{figure}[h!]
 	\centering\tabcolsep=0.5mm
 	\begin{tabular}[c]{cc}
  \includegraphics[width=5.8cm,height=2.8cm]{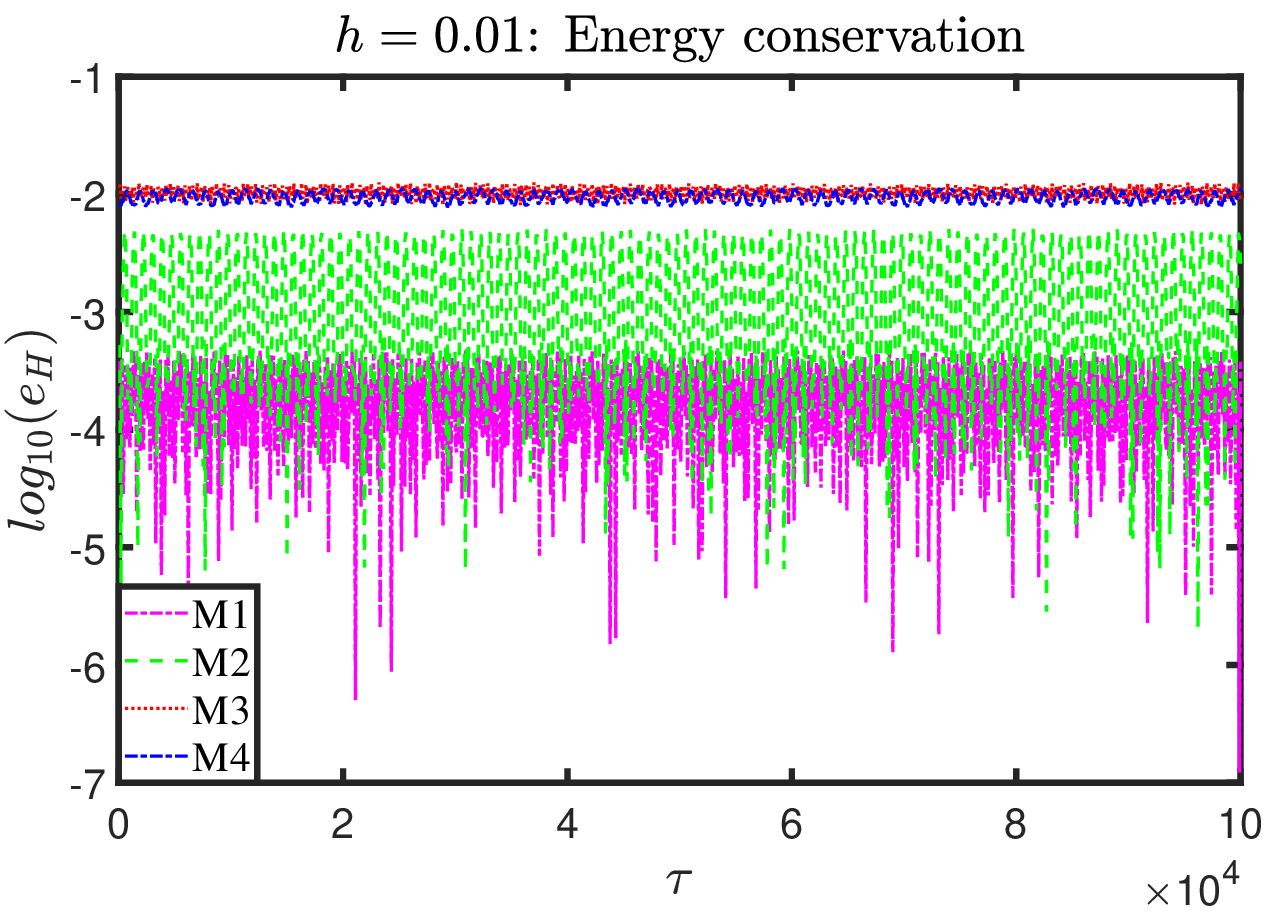}\qquad\quad \includegraphics[width=5.8cm,height=2.8cm]{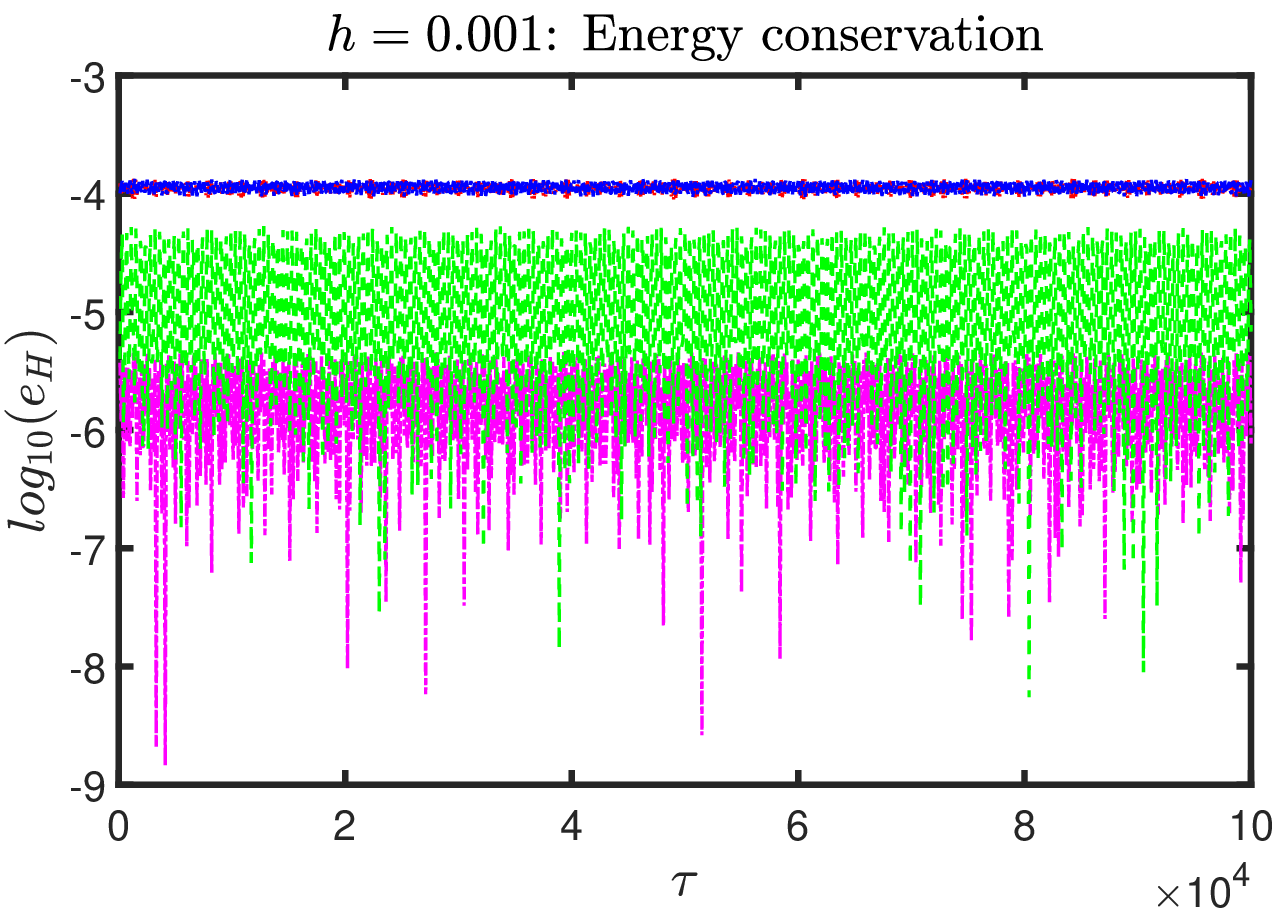}\\
  	\includegraphics[width=5.8cm,height=2.8cm]{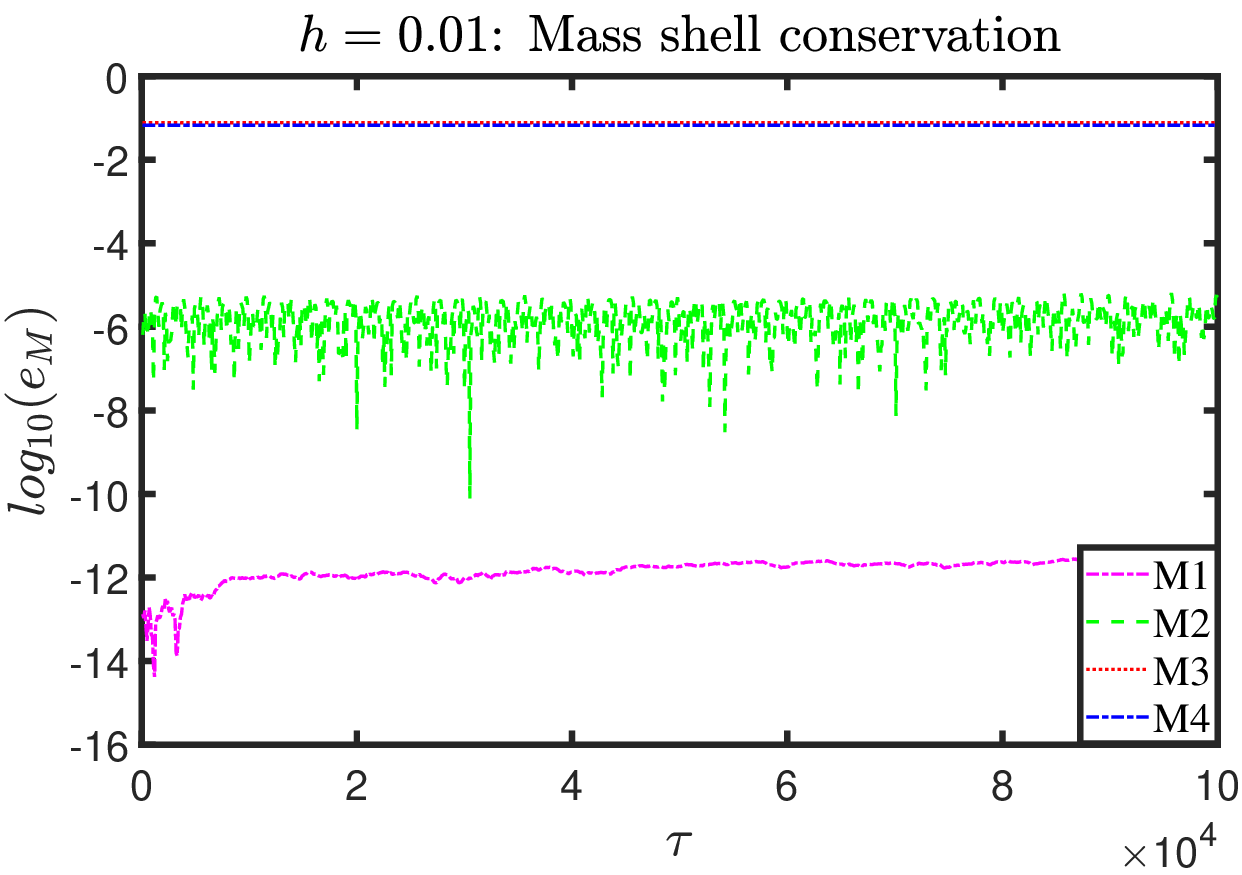}\qquad\quad
  \includegraphics[width=5.8cm,height=2.8cm]{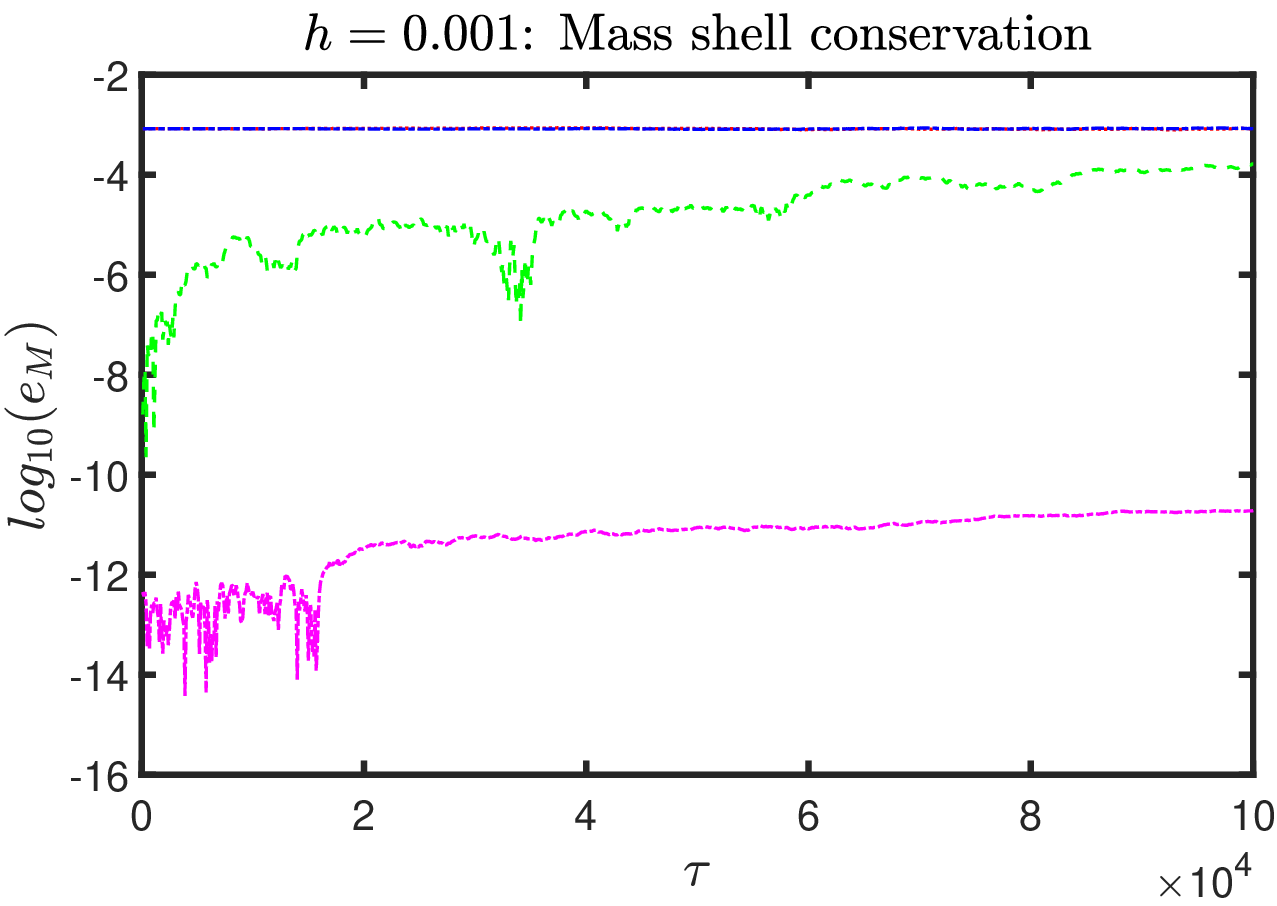}
 		\end{tabular}
 	\caption{Problem  5 (strong constant magnetic field).  Evolution of  $e_{H}$ and $e_{M}$ with step sizes $h=0.01$ (left) and $h=0.001$ (right).}
 	\label{fig:problem5H}
 \end{figure}

\textbf{Problem 5 (Strong constant magnetic field).}\label{Problem5}
The last numerical experiment considers a strong constant magnetic field $\frac{1}{\epsilon}B=\frac{1}{\epsilon}(0,0,1)^{\intercal}$, where $0<\epsilon\ll 1$ is a dimensionless parameter inversely proportional to the strength of the magnetic field. The potential is given by  $U(x)=x_1^2+2x_2^2+3x_3^2-x_1$, with electric field $E(x)=-\nabla U(x)$. The initial values are  chosen as  $(x(0); t(0))=(0,1,0.1,0)^{\intercal}$ and  $(v(0); w(0))=(0.09,0.55,0.3,w(0))^{\intercal}$, where $w(0)=i\sqrt{1+v(0)^{\intercal}v(0)}$. Figures \ref{fig:problem5Err} and \ref{fig:problem5H} illustrate the global errors and the numerical conservation of
$H$ and $\mathcal{H}$, respectively.

The numerical results shown in Figures \ref{fig:problem1Banana}--\ref{fig:problem4H} lead to the following observations:

\textbf{Order behavior.} The
global errors $error_{x}$ and $error_{v}$ at $\tau=1$ are displayed in Figures \ref{fig:problem1Err}, \ref{fig:problem2Err}, \ref{fig:problem3Err}, and  \ref{fig:problem4Err}.  These results indicate that the global error lines for M1--M4 are nearly parallel to a line with slope 2, confirming that all methods achieve second-order accuracy. Moreover, as shown in Figure \ref{fig:problem1Banana}, the projected particle orbit in the $(R, x_3)$ plane exhibits a banana-shaped trajectory. All proposed methods effectively preserve the geometric structure of the orbit.

\textbf{Energy conservation.}
In the Tokamak field case shown in Figure \ref{fig:problem1H} (top row), the numerical results show that M1 exactly preserves energy, which can be directly deduced from its formulation in the $(v;w)$ components. Meanwhile,  M2--M4 nearly preserve energy with high accuracy over long times.
As shown in Figure \ref{fig:problem2H} (top row), M2--M4 exhibit long-term near-conservation of energy under a quadratic electric potential, consistent with the result stated in Theorem \ref{conservation-H1}.
Figure \ref{fig:problem3H} illustrates a slow drift in the energy error for all methods under a non-quadratic electric potential.  Furthermore, as shown in Figure \ref{fig:problem4H} (top row), M1  demonstrates long-term near-conservation of energy under a constant magnetic field, even for a non-quadratic potential, whereas M2--M4  exhibit a noticeable energy drift under the same conditions. However, a rigorous theoretical justification for the long-term energy behavior of M1 is not yet available. This issue will be addressed in future work.

\textbf{Mass shell conservation.}   Figures \ref{fig:problem1H} (bottom row), \ref{fig:problem2H} (bottom row), \ref{fig:problem3M}, and \ref{fig:problem4H} (bottom row)  demonstrate the favorable long-term behavior of mass shell  along the numerical solutions obtained using M1--M4. The results clearly show that M1 exactly preserves the mass shell, while M2--M4 nearly preserve it over long times.

\textbf{Behavior under strong magnetic fields.} Beyond the moderate magnetic field ($\epsilon=1$) of Problem 1--4, we further examined the performance of the methods in the strong magnetic field regime ($0<\epsilon\ll 1$) in Problem 5. Although a rigorous theoretical justification is not yet available, the numerical results clearly indicate that the proposed methods M1--M4 still exhibit second-order accuracy and favorable long-term conservation properties in this regime. This will be addressed in future work.

\section{Error bounds and long-term analysis of M1--M4 (Proofs of Theorems \ref{error bound2}--\ref{conservation-vol})}\label{Error and long-term analysis}
In this section, we provide the proofs of the error bounds and long-term analysis for the two-step symmetric methods M1--M4. To demonstrate the long-term near-conservation of energy and mass shell, we employ the backward error analysis (see Chap. IX of \cite{Hairer2002}). By interpreting the numerical scheme as the exact solution of a modified differential equation, one obtains modified invariants approximating the original energy and mass shell, which remain nearly conserved over long times.

\subsection{Proof of Theorem \ref{error bound2}}\label{error analysis}

The proof is given only for the error bounds of M2, and the results for the other methods can be derived by a similar procedure.
We first derive the global error for the position vector $(x;t)$, and the result for the momentum vector $(v;w)$ can then be obtained in the same way.
		
The exact solution of the relativistic dynamics \eqref{rela-CPD-4d} can be expressed using the variation-of-constants formula as follows
		\begin{equation}\label{exact-solution}
			\begin{aligned}	
				{\bm{y}}(\tau_n+\sigma h)&=	{\bm{y}}(\tau_n) +\sigma h{\bm{u}}(\tau_n)+h^2	\int_{0}^{\sigma}(\sigma-z) {\bm{F}}(x(\tau_n+hz)){\bm{u}}(\tau_n+hz)dz,\\
				{\bm{u}}(\tau_n+\sigma h)&={\bm{u}}(\tau_n)+h
				\int_{0}^{\sigma} {\bm{F}}(x(\tau_n+hz)){\bm{u}}(\tau_n+hz)  dz.
			\end{aligned}
		\end{equation}
Setting $\sigma=-1$ and $\sigma=1$ in \eqref{exact-solution} and summing the two equations yields
		\begin{equation}\label{exact-two step}
			\begin{aligned}
				&x(\tau_n+h)-2x(\tau_n)+x(\tau_n-h)
				=h^2\int_{0}^{1}(1-z)\big( \tilde{B}(x(\tau_n+hz))v(\tau_n+hz)\\
			&\  + \tilde{B}(x(\tau_n-hz))v(\tau_n-hz)-iE(x(\tau_n+hz))w(\tau_n+hz)\\
			&\	-iE(x(\tau_n-hz))w(\tau_n-hz)\big)dz,\\
				&t(\tau_n+h)-2t(\tau_n)+t(\tau_n-h)=h^2\int_{0}^{1}(1-z)\big(
				iE(x(\tau_n+hz))^{\intercal}v(\tau_n+hz)\\
				&\ +iE(x(\tau_n-hz))^{\intercal}v(\tau_n-hz)
				\big)dz.
			\end{aligned}
		\end{equation}

		 \textbf{Local error of M2.} We begin by analyzing the local error of M2. Substituting the exact solution \eqref{exact-solution} into \eqref{M2} yields
		\begin{equation}\label{defect-two step}
			\begin{aligned}
				&x(\tau_{n+1})-2x(\tau_{n})+x(\tau_{n-1})
				=\big(S_1\tilde{B}( x(\tau_{n}))+S_2\tilde{E}( x(\tau_{n}))\big)\big(x(\tau_{n+1})-x(\tau_{n-1})\big)\\
				&\qquad\qquad\qquad\qquad\qquad\qquad\qquad-i\big(S_1E( x(\tau_{n})))-S_2B_0\big) \big(t(\tau_{n+1})-t(\tau_{n-1})\big)+\Delta_{x},\\[0.5em]
				&t(\tau_{n+1})-2t(\tau_{n})+t(\tau_{n-1})
				=i\big(S_1E( x(\tau_{n}))^{\intercal}-S_2B_0^{\intercal}\big)\big(x(\tau_{n+1})-x(\tau_{n-1})\big)+\Delta_{t},
			\end{aligned}
		\end{equation}	
where $\Delta_{x}$ and $\Delta_{t}$ denote the local discrepancies.	Combining the formula \eqref{exact-two step} with the Taylor expansions yields	
		\begin{equation*}\label{}
		\begin{aligned}
		\Delta_{x}
		&=h^2\int_{0}^{1}(1-z)\big( \tilde{B}(x(\tau_n+hz))v(\tau_n+hz)+\tilde{B}(x(\tau_n-hz))v(\tau_n-hz)\big)dz\\
		&\quad  -ih^2\int_{0}^{1}(1-z)\big(E(x(\tau_n+hz))w(\tau_n+hz)+E(x(\tau_n-hz))w(\tau_n-hz) \big)dz\\
		&\quad -\big(S_1\tilde{B}( x(\tau_{n}))+S_2\tilde{E}( x(\tau_{n}))\big)\big(x(\tau_{n+1})-x(\tau_{n-1})\big)-iS_2B_0\big(t(\tau_{n+1})-t(\tau_{n-1})\big) \\
		&\quad +iS_1E(x(\tau_{n}))\big(t(\tau_{n+1})-t(\tau_{n-1})\big)\\
		&=h^2\big(\tilde{B}( x(\tau_{n}))v(\tau_{n})-iE( x(\tau_{n}))w(\tau_{n})\big)+2ih\big(S_1E( x(\tau_{n})))-S_2B_0\big)\dot{t}(\tau_{n})+\mathcal{O}(h^3)	\\	
		&\quad-2h\big(S_1\tilde{B}( x(\tau_{n}))+S_2\tilde{E}( x(\tau_{n}))\big)\dot{x}(\tau_{n})\\
		&=h^2\big(\tilde{B}( x(\tau_{n}))v(\tau_{n})-iE( x(\tau_{n}))w(\tau_{n})\big)
		-h^2\tilde{B}( x(\tau_{n}))\dot{x}(\tau_{n})+\mathcal{O}(h^3)\\
		&\quad+ih^2E( x(\tau_{n}))\dot{t}(\tau_{n}).
    	\end{aligned}
		\end{equation*}
Here,  $\dot{(\cdot)}$ represents the derivative with respect to $\tau$,
and we use the Taylor expansions $$S_1=\frac{h}{2}-\frac{h^3}{24}r_1+\mathcal{O}(h^5),\ \ S_2=-\frac{h^3}{24}r_2+\frac{h^5}{240}r_1r_2+\mathcal{O}(h^7).$$ 
Thus, $\Delta_{x}=\mathcal{O}(h^3)$. Note that the Taylor expansion of $S_2$ begins at order $\mathcal{O}(h^3)$. The corresponding term $S_2 B_0$ is of the same order. Therefore, the modification in the $(x;t)$ components of M1 that leads to M2 does not affect the local error $\Delta_{x}$ or the subsequent analysis. Similarly, omitting the $S_2$ term has no impact on $\Delta_{x}$ or on the analysis, which ensures that both M3 and M4 also retain second-order accuracy. Furthermore, using the same argument, we obtain $\Delta_{t}=\mathcal{O}(h^3)$.
		
\textbf{Global error of M2.} Denote the global error of  position vector $(x; t)$ by $$e_{n}^x=x(\tau_n)-x_n,\qquad  e_{n}^{t}=t(\tau_n)-t_n.$$ Subtracting \eqref{M2} from \eqref{defect-two step} yields
\begin{equation*}
	\begin{aligned}
		&	e_{n+1}^x-2e_{n}^x+e_{n-1}^x
		=\big(S_1\tilde{B}( x(\tau_{n}))+S_2\tilde{E}( x(\tau_{n}))\big)\big(x(\tau_{n+1})-x(\tau_{n-1})\big)\\
		&\ \ -S_2\tilde{E}( x_{n})\big(x_{n+1}-x_{n-1}\big)
		-S_1\tilde{B}( x_{n})\big(x_{n+1}-x_{n-1}\big) +iS_1E( x_{n})\big(t_{n+1}-t_{n-1}\big)\\
		&\ \ -iS_2B_0\big(t_{n+1}-t_{n-1}\big)-i\big(S_1E( x(\tau_{n})))-S_2B_0\big)\big(t(\tau_{n+1})-t(\tau_{n-1})\big)+\Delta_{x}\\
		&=S_1\tilde{B}( x(\tau_{n}))\big(x(\tau_{n+1})-x_{n+1}\big)
		+S_1\big(\tilde{B}( x(\tau_{n})-\tilde{B}( x_{n}))x_{n+1}
		+S_1\tilde{B}( x_{n})x_{n-1}\\
		&\ \ -S_1\tilde{B}( x(\tau_{n}))x_{n-1}- S_1\tilde{B}( x(\tau_{n}))\big(x(\tau_{n-1})-x_{n-1}\big)+S_2\tilde{E}(x(\tau_{n}))\big(x(\tau_{n+1})-x_{n+1}\big)
		\\
		&\ \ +S_2\big(\tilde{E}(x(\tau_{n}))-\tilde{E}(x_{n})\big)x_{n+1}
		-S_2\tilde{E}(x(\tau_{n}))\big(x(\tau_{n-1})-x_{n-1}\big)+S_2\tilde{E}(x_{n})x_{n-1}
		\\
		&\ \ -S_2\tilde{E}(x(\tau_{n}))x_{n-1}
		-iS_1E( x_{n})\big(t(\tau_{n+1})-t_{n+1}\big) 
		-iS_1\big(E( x(\tau_{n}))-E( x_{n})\big)t(\tau_{n+1})
		\\
		&\ \ +iS_1E( x_{n})\big(t(\tau_{n-1})-t_{n-1}\big)
		+iS_1\big(E( x(\tau_{n}))-E( x_{n})\big)t(\tau_{n-1})+iS_2B_0t(\tau_{n+1})\\
		&\ \  -iS_2B_0t_{n+1}-iS_2B_0\big(t(\tau_{n-1})-t_{n-1}\big)+\Delta_{x}.
	\end{aligned}
\end{equation*} 
From the Taylor expansions of $S_1$ and $S_2$, it follows that
		\begin{equation*}
			\begin{aligned}
				&\abs{e_{n+1}^x-2e_{n}^x+e_{n-1}^x} 
				\leq Ch\big(\abs{e_{n+1}^x}+\abs{e_{n-1}^x}+\abs{e_{n+1}^{t}}+\abs{e_{n-1}^{t}}\big)\\
				&\quad +Ch^3\big(\abs{e_{n+1}^x}+\abs{e_{n-1}^x}+\abs{e_{n+1}^{t}}+\abs{e_{n-1}^{t}}\big)+\mathcal{O}(h^3).
			\end{aligned}
		\end{equation*}	
Using the recursion formula and the condition $e_{0}^x=0$, we obtain
$$
\abs{e_{n+1}^x-e_{n}^x}-\abs{e_{1}^x} \leq C(h+h^3)\sum_{m=1}^{n}\big(\abs{e_{m+1}^x}+\abs{e_{m-1}^x}+\abs{e_{m+1}^{t}}+\abs{e_{m-1}^{t}}\big)+\mathcal{O}(h^2).
$$
Applying the recursion formula once again yields
$$
\abs{e_{n+1}^x}-n\abs{e_{1}^x} \leq C(h+h^3)\sum_{m=1}^{n}\big(\abs{e_{m+1}^x}+\abs{e_{m-1}^x}+\abs{e_{m+1}^{t}}+\abs{e_{m-1}^{t}}\big)+\mathcal{O}(h).
$$
From the Gronwall's inequality, we conclude that if 
$e_{1}^x=\mathcal{O}(h^2)$, the global error bound $e_{n+1}^x=\mathcal{O}(h)$  follows immediately.
	
We now calculate the value of $e_{1}^x$ using the approach from the previous analysis. The starting value ${\bm{y}}_1=(x_1;t_1)$ is given by  
\begin{equation}\label{start_x1}
\begin{aligned}
{\bm{y}}_1=	{\bm{y}}_0+he^{\frac{h}{2}{\bm{F}}(x_{0})}{\bm{u}}_0.
\end{aligned}
\end{equation}	
The variation-of-constants formula \eqref{exact-solution} leads to
		\begin{equation}\label{eaxct_x1}
			{\bm{y}}(\tau_0+h)=	{\bm{y}}(\tau_0) + h{\bm{u}}(\tau_0)+h^2
			\int_{0}^{1}(1-z) {\bm{F}}(x(\tau_0+hz)){\bm{u}}(\tau_0+hz)
			dz.
		\end{equation}
		Inserting the exact solution \eqref{eaxct_x1} into \eqref{start_x1}, we obtain
		\begin{equation}\label{local_x1}
			\begin{aligned}
				{\bm{y}}(\tau_1)=	{\bm{y}}(\tau_0)+he^{\frac{h}{2}{\bm{F}}(x(\tau_0))}{\bm{u}}(\tau_0)
				+\Delta_1.
			\end{aligned}
		\end{equation}
	The above two equations imply that
		\begin{equation*}
			\begin{aligned}
				\Delta_1
				=h{\bm{u}}(\tau_0)+h^2
				\int_{0}^{1}(1-z) {\bm{F}}(x(\tau_0+hz)){\bm{u}}(\tau_0+hz)
				dz-he^{\frac{h}{2}{\bm{F}}(x(\tau_0))}{\bm{u}}(\tau_0)
				=\mathcal{O}(h^2),
			\end{aligned}
		\end{equation*}	
		which follows directly from the Taylor expansions.  Comparing  \eqref{local_x1} with \eqref{start_x1} yields
		\begin{equation*}
			\begin{aligned}
				e_{1}^{\bm{y}}={\bm{y}}(\tau_1)-{\bm{y}}_1={\bm{y}}(\tau_0)+he^{\frac{h}{2}{\bm{F}}(x(\tau_0))}{\bm{u}}(\tau_0)-{\bm{y}}_0-he^{\frac{h}{2}{\bm{F}}(x_{0})}{\bm{u}}_0
				+	\Delta_1
				=\mathcal{O}(h^2).
			\end{aligned}
		\end{equation*}	
Thus, $e_{1}^x=\mathcal{O}(h^2)$, which implies the global error in position satisfies $e_{n+1}^x=\mathcal{O}(h)$. Similarly, the global error in $t$ is  $e_{n+1}^{t}=\mathcal{O}(h)$.  Moreover, it is shown in \cite{Hairer2002}  that the order of a symmetric method is always even. Therefore, the second-order global error in the position vector $(x;t)$ of Theorem \ref{error bound2} is obtained. An analogous argument yields the corresponding estimate for the momentum vector $(v;w)$. 

\subsection{Proof of Theorem \ref{conservation-H1}}\label{energy analysis}
The analysis focuses on M2, and the results for M3 and M4 follow analogously.
Using the idea of the backward error analysis, we seek an 
$h$-dependent modified differential equation, given as a formal power series in $h$, whose solutions  $x(\tau)$ and $t(\tau)$  formally satisfy $x(nh)=x_n$ and $t(nh)=t_n$, where $x _n$ and $t_n$ are the numerical solutions obtained from M2. Applying the scheme \eqref{M2}, we obtain		
\begin{equation*}\label{modTSI}
			\begin{aligned}
&\frac{x(\tau+h)-2x(\tau)+x(\tau-h)}{h^2}
				=\frac{2}{h}\big(S_1\tilde{B}( x(\tau))+S_2\tilde{E}( x(\tau))\big)\frac{x(\tau+h)-x(\tau-h)}{2 h}\\
				&\qquad\qquad\qquad\qquad\qquad\qquad\qquad\ -\frac{2}{h}i\big(S_1E( {x(\tau)})-S_2B_{0}\big)\frac{t(\tau+h)-t(\tau-h)}{2 h},\\
				&\frac{{t(\tau+h)}-2{t(\tau)}+{t(\tau-h)}}{h^2}
				=\frac{2}{h}i\big(S_1E( {x(\tau)})^{\intercal}-S_2B_{0}^{\intercal}\big)\frac{x(\tau+h)-x(\tau-h)}{2 h}.
			\end{aligned}
		\end{equation*} 				
		From the Taylor expansions, it follows  that
		\begin{equation}\label{original-Mod}
			\begin{aligned}
				&	\ddot{x}+\frac{h^2}{12}\ddddot{x}+\cdots=\frac{2}{h}\big(S_1\tilde{B}( x)+S_2\tilde{E}( x)\big)\Big(\dot{x}+\frac{h^2}{6}\dddot{x}+\cdots \Big)\\
				&\qquad\qquad\qquad\qquad -\frac{2}{h}i\big(S_1E( x)-S_2B_{0}\big)\Big(\dot{t}+\frac{h^2}{6}\dddot{t}+\cdots \Big),\\
				&	\ddot{t}+\frac{h^2}{12}\ddddot{t}+\cdots
				=\frac{2}{h}i\big(S_1E( x)^{\intercal}-S_2B_{0}^{\intercal}\big)\Big(\dot{x}+\frac{h^2}{6}\dddot{x}+\cdots \Big),
			\end{aligned}
		\end{equation}
where the argument  $\tau$ is omitted for simplicity. 
Recalling the Taylor expansions $S_1=\frac{h}{2}-\frac{h^3}{24}r_1+\mathcal{O}(h^5)$ and $S_2=-\frac{h^3}{24}r_2+\frac{h^5}{240}r_1r_2+\mathcal{O}(h^7)$, it follows that \eqref{original-Mod} has a formal expansion in even powers of $h$. Thus, we have derived a modified differential equation that can formally be represented as a power series in $h$.

We now show that the modified differential equation possesses a formal first integral that approximates the energy $H$.  This derivation relies on the second formula of \eqref{original-Mod}. For the quadratic potential $U(x)=\frac{1}{2}x^{\intercal}Qx+q^{\intercal}x$, we have $E(x)=-\nabla U(x)=-(Qx+q)$, and the solution of the modified differential equation satisfies
$E(x)^{\intercal}\dot{x}=-\frac{d}{d\tau}U(x)$. Moreover, for any integer $l\geq 0$, one readily verifies that $x^{\intercal}Qx^{(2l+1)}$ is a total derivative, namely,
$$
x^{\intercal}Qx^{(2l+1)}=\frac{d}{dt}\big(x^{\intercal}Qx^{(2l)}-\cdots \pm\big({x^{(l-1)}}\big)^{\intercal}Qx^{(l+1)}\mp\frac{1}{2}\big({x^{(l)}}\big)^{\intercal}Qx^{(l)}\big).$$
Thus, by substituting the Taylor expansions of $S_1$ and $S_2$ into the second formula of \eqref{original-Mod}, the equation can be expressed as a total derivative. We emphasize that the above formal argument is valid on time intervals of length $\mathcal{O}(1)$. Consequently, there exist  $h$-independent functions $H_{2j}(x, t)$ such that the function
		\begin{equation}\label{modH1}
			H_{h}(x, t)=-i\dot{t}+U(x)+h^2H_{2}(x, t)
			+h^4H_{4}(x, t)+\cdots,
		\end{equation}
		truncated at the $\mathcal{O}(h^N)$ term, satisfies $$\frac{d}{d\tau}	H_{h}(x, t)=\mathcal{O}(h^N)$$
		along solutions of the modified differential equation.   
Using the $h$-dependent function      
\begin{equation*}
\begin{aligned}
i\gamma(\tau)&= iS_4E( x(\tau))^{\intercal}\tilde{B}( x(\tau))\big(x(\tau+h)-x(\tau-h)\big)/2h\\
&\quad \ +\big(S_3+S_4E( x(\tau))^{\intercal}E( x(\tau))\big)\big(t(\tau+h)-t(\tau-h)\big)/2h,
\end{aligned}
\end{equation*}
which yields
\begin{equation}\label{gamma}
\gamma= S_4E( x)^{\intercal}\tilde{B}( x)\Big(\dot{x}+\frac{h^2}{6}\dddot{x}+\cdots \Big)-i\big(S_3+S_4E( x)^{\intercal}E(x)\big)\Big(\dot{t}+\frac{h^2}{6}\dddot{t}+\cdots \Big).
\end{equation}
Substituting the Taylor expansions $S_3=1+ \frac{5h^4}{384} r_2^2 +\mathcal{O}(h^6)$ and $S_4 = -\frac{h^2}{8}+ \frac{5h^4}{384} r_1 +\mathcal{O}(h^6)$ into \eqref{gamma} yields 
\begin{equation}\label{texpression}
	\dot{t}=i\gamma+\mathcal{O}(h^2).
\end{equation}
Combining this result with \eqref{modH1} and the energy definition \eqref{E} leads to  $$H_h(x,t)=H(x,\gamma)+\mathcal{O}(h^2).$$

We prove that M2 nearly preserves the energy $H$ over long time intervals by patching together short-time estimates on subintervals of length $h$.
If the numerical solution remains in a compact set independent of $h$, a telescoping sum gives
\begin{equation*}
\begin{aligned}
&H(x(nh),\gamma(nh))-H(x(0),\gamma(0))
=H(x(nh),\gamma(nh))-H_{h}(x(nh),\gamma(nh))\\
&\quad+\sum_{j=1}^{n}\big(H_{h}(x(jh),\gamma(jh))-H_{h}(x((j-1)h),\gamma((j-1)h))\big)\\
&\quad+H_{h}(x(0),\gamma(0))-H(x(0),\gamma(0))=\mathcal{O}(h^2)+\mathcal{O}(nh^{N+1}).
\end{aligned}
\end{equation*}
Therefore, for any fixed integer $N \geq 3$, it follows that
	\begin{equation*}
	\abs{	H(x_{n},\gamma_n)-	H(x_{0},\gamma_0)} \leq Ch^2 \quad \ \textmd{for} \quad nh \leq h^{-N+2}.
\end{equation*}	
This establishes the statement of Theorem \ref{conservation-H1}.

\begin{remark}\label{remark3}
For a constant magnetic field $B(x)=B$, M1 nearly preserves energy over long times, independently of the form of the potential $U(x)$. This behavior is numerically illustrated for a non-quadratic potential in Problem~4 of Section \ref{experiments}. However, this cannot be established by the same analysis as in Theorem \ref{conservation-H1}, for reasons explained below. From the formulation of M1, we obtain 
	\begin{equation*}
		\ddot{t}+\frac{h^2}{12}\ddddot{t}+\cdots	=\frac{2}{h}i\big(K_1(x)E(x)^{\intercal}-K_2(x)B^{\intercal}\big)\Big(\dot{x}+\frac{h^2}{6}\dddot{x}+\cdots \Big),
	\end{equation*}
	where $B$ is a constant magnetic field.  Substituting $E( x)^{\intercal}(\dot{x}+\frac{h^2}{6}\dddot{x}+\cdots )$, obtained from the first relation of \eqref{M1}, into the above formula yields
	\begin{equation}\label{long-term-M1}
		\begin{aligned}
			\ddot{t}+\frac{h^2}{12}\ddddot{t}+\cdots
			=&i\Big(\dot{t}+\frac{h^2}{6}\dddot{t}+\cdots\Big)^{-1}\Big(\ddot{x}+\frac{h^2}{12}\ddddot{x}+\cdots-\frac{2}{h}iK_2(x)B\Big(\dot{t}+\frac{h^2}{6}\dddot{t}+\cdots \Big)\Big)^{\intercal}\\
			&\Big(\dot{x}+\frac{h^2}{6}\dddot{x}+\cdots \Big)
			-\frac{2}{h}iK_2(x)B^{\intercal}\Big(\dot{x}+\frac{h^2}{6}\dddot{x}+\cdots \Big).
		\end{aligned}
	\end{equation}
Based on the backward error analysis, we aim to find a function $H_{h}(x, \gamma)=H(x, \gamma)+h^2	H_{2}(x, \gamma)+h^4H_{4}(x, \gamma)+\cdots$ such that its derivative with respect to  $\tau$ satisfies $	\frac{d}{d\tau}	H_{h}(x, \gamma)=\mathcal{O}(h^N)$. However, such a function cannot currently be derived from \eqref{long-term-M1}. The difficulty arises from two main aspects. First, the momentum approximation $(v_{n+1};w_{n+1})$ in M1 does not readily provide an explicit or usable relation between $t$ and $\gamma$ easily, which is crucial for constructing $H_{h}(x, \gamma)$. Second, since $K_2(x)$ is related to $E(x)$, it becomes challenging to express the right-hand side of equation \eqref{long-term-M1} as a total differential. Therefore, the theoretical analysis of long-term near-conservation of energy for M1 is not achievable at present. We will address this issue in future work.
\end{remark}

\subsection{Proof of Theorem \ref{conservation-H2}}
\textbf{Proof of \eqref{MS0}.}
Substituting the approximate momentum variables $(v_{n+1};w_{n+1})$ from M1 into the mass shell \eqref{mass shell} yields
\begin{equation*}
	\begin{aligned}
		\mathcal{H}(v_{n+1},w_{n+1})
		&=({v_{n}};w_{n})^{\intercal}e^{\frac{h}{2}{\bm{F}}(x_{n})^{\intercal}}e^{\frac{h}{2}{\bm{F}}(x_{n+1})^{\intercal}}e^{\frac{h}{2}{\bm{F}}(x_{n+1})}e^{\frac{h}{2}{\bm{F}}(x_{n})}({v_{n}};w_{n})/2\\
        &=(v_{n}^2+w_{n}^2)/2
		=\mathcal{H}(v_{n},w_{n}),
	\end{aligned}
\end{equation*}
where we have used the skew-symmetry  of ${\bm{F}}(x_{n})$ and ${\bm{F}}(x_{n+1})$. This establishes the result of \eqref{MS0}.\\
\textbf{Proof of \eqref{MS2}.}
This proof also employs the backward error analysis from Theorem~\ref{conservation-H1}, with a focus on the differences in deriving the modified differential equation. We concentrate on M2, as the arguments for M3 and M4 are analogous. 

Multiplying the first equation of \eqref{original-Mod} by $\big(\dot{x}+\frac{h^2}{6}\dddot{x}+\cdots\big)^{\intercal}$, we have 
		\begin{equation*}\label{}
		\begin{aligned}
			&	\Big(\dot{x}+\frac{h^2}{6}\dddot{x}+\cdots\Big)^{\intercal}\Big(\ddot{x}+\frac{h^2}{12}\ddddot{x}+\cdots\Big)\\
			&	=\frac{2}{h}\Big(\dot{x}+\frac{h^2}{6}\dddot{x}+\cdots\Big)^{\intercal}\big(S_1\tilde{B}( x)+S_2\tilde{E}(x)\big)\Big(\dot{x}+\frac{h^2}{6}\dddot{x}+\cdots\Big)\\
			&\quad\
		 -\frac{2i}{h}\Big(\dot{x}+\frac{h^2}{6}\dddot{x}+\cdots\Big)^{\intercal}\big(S_1E( x)-S_2B_0\big)\Big(\dot{t}+\frac{h^2}{6}\dddot{t}+\cdots\Big).
			\end{aligned}
		\end{equation*}
By the skew-symmetry of $\tilde{B}$ and $\tilde{E}$, it  follows that
		\begin{equation}\label{xt}
			\begin{aligned}
				&	\Big(\dot{x}+\frac{h^2}{6}\dddot{x}+\cdots\Big)^{\intercal}\Big(\ddot{x}+\frac{h^2}{12}\ddddot{x}+\cdots\Big)\\
				&= -\frac{2i}{h}\Big(\dot{x}+\frac{h^2}{6}\dddot{x}+\cdots\Big)^{\intercal}\big(S_1E( x)-S_2B_0\big)\Big(\dot{t}+\frac{h^2}{6}\dddot{t}+\cdots\Big).
			\end{aligned}
		\end{equation}
		Multiplying the second formula of \eqref{original-Mod} by $\dot{t}+\frac{h^2}{6}\dddot{t}+\cdots$ yields
		\begin{equation*}\label{}
		\begin{aligned}		
				&	\Big(\dot{t}+\frac{h^2}{6}\dddot{t}+\cdots\Big)\Big(\ddot{t}+\frac{h^2}{12}\ddddot{t}+\cdots\Big)\\	&=\frac{2i}{h}\Big(\dot{t}+\frac{h^2}{6}\dddot{t}+\cdots\Big)\big(S_{1}E( x)^{\intercal}-S_2B_0^{\intercal}\big)\Big(\dot{x}+\frac{h^2}{6}\dddot{x}+\cdots \Big).
			\end{aligned}
		\end{equation*}
		Combining this with the equation \eqref{xt}, we obtain
		\begin{equation*}\label{ph}
		\begin{aligned}
				\Big(\dot{x}+\frac{h^2}{6}\dddot{x}+\cdots\Big)^{\intercal}\Big(\ddot{x}+\frac{h^2}{12}\ddddot{x}+\cdots\Big)=-\Big(\dot{t}+\frac{h^2}{6}\dddot{t}+\cdots\Big)\Big(\ddot{t}+\frac{h^2}{12}\ddddot{t}+\cdots\Big).
			\end{aligned}
		\end{equation*}
Note that the scalar product ${x^{(k)}}^{\intercal}{x^{(l)}}$ is a total differential whenever $k+l$ is odd (see \cite{Hairer2017}). Hence, there exist $h$-independent functions $\mathcal{H}_{2j}(x,t)$ such that the function
\begin{equation}\label{H_h}
\mathcal{H}_{h}(x,t)=\big(\dot{x}^{\intercal}\dot{x}+\dot{t}^2\big)/2+h^2\mathcal{H}_{2}(x,t)
+h^4\mathcal{H}_{4}(x,t)+\cdots,
\end{equation}
truncated at the $\mathcal{O}(h^N)$ term, satisfies $\frac{d}{d\tau}\mathcal{H}_{h}(x,t)=\mathcal{O}(h^N)$
along solutions of the modified differential equation. Using the expression for $v_n$,
we have
\begin{equation*}\label{}
\begin{aligned}
&v
=\big(S_3I+S_4\big(\tilde{B}( {x(\tau)})^2+E( {x(\tau)})E( {x(\tau)})^{\intercal}\big)\big)\Big(\dot{x}+\frac{h^2}{6}\dddot{x}+\cdots \Big)\\
&\qquad -iS_4\tilde{B}( {x(\tau)})E( {x(\tau)})\Big(\dot{t}+\frac{h^2}{6}\dddot{t}+\cdots \Big).\\
\end{aligned}
\end{equation*} 
By combining the Taylor expansions of $S_3$ and $S_4$ with the above formula, we obtain $v=\dot{x}+\mathcal{O}(h^2)$. Substituting this relation and \eqref{texpression} into \eqref{H_h}, and applying the definition of the mass shell \eqref{mass shell}, we obtain
$$\mathcal{H}_{h}(x,t) = \mathcal{H}(v,w) + \mathcal{O}(h^2).$$ The same argument as in the proof of Theorem \ref{conservation-H1} yields \eqref{MS2}. This completes the proof of Theorem \ref{conservation-H2}.

\subsection{Proof of Theorem \ref{conservation-vol}}
As in the original system \eqref{rela-CPD-4d}, one can check that the
right-hand side vector fields of the subsystems $Z_1$ and $Z_2$ are
divergence-free, so their exact flows $\Psi_h^i$ $(i=1,2)$ are
volume-preserving. Since each subflow is volume-preserving and this property is preserved under composition, the Strang splitting scheme
$
\Psi_h=\Psi_{h/2}^{1}\circ \Psi_{h/2}^{2}\circ \Psi_{h/2}^{1}
$
satisfies
\[
\operatorname{vol}(\Psi_h(\Omega))=\operatorname{vol}(\Omega)\qquad
\forall \Omega\subset\mathbb{R}^8,
\]
which establishes the volume conservation property stated in Theorem \ref{conservation-vol}.

\section{Conclusions}
\label{sec:conclusions}
In this paper, we proposed four two-step symmetric methods for solving relativistic charged-particle dynamics. First, we derived a two-step symmetric method based on a splitting scheme. This method exactly preserved both the mass shell and the phase-space volume of the system. Building on this, three other two-step symmetric methods were formulated by introducing appropriate modifications to the first. We proved second-order convergence for all methods and established long-time conservation or near-conservation of energy, mass shell, and phase-space volume. Numerical experiments demonstrated that all proposed methods achieved second-order accuracy and preserved or nearly preserved geometric structures over long times.

Last but not least, several issues remain to be addressed in future research. 
(a) The long-term near-conservation of energy by the M1 method under a constant magnetic field requires further investigation. 
(b) Error analysis and geometric conservation under strong magnetic fields will be analyzed using the modulated Fourier expansion technique \cite{Hairer2018, Lubich2022, Hairer2002, Lubich2020}. 
(c) Uniformly accurate integrators \cite{Chartier2019,Chartier2020,WangB2023} for relativistic charged-particle dynamics appears to be a promising direction for future research.

		\section*{Conflict of interest}
		The authors declare that they have no known competing financial interests or personal
		relationships that could have appeared to influence the work reported in this paper.
		
		\section*{Acknowledgments}
		
	The authors thank Professor Christian Lubich for his valuable discussions and helpful comments.	The research was supported by NSFC (12371403).

\end{document}